\documentclass[letterpaper, 10pt]{article}
\usepackage{preprint}

\title{Degree-of-Freedom and Optimization-Dynamic Effects on the Observability of Kuramoto--Sivashinsky Systems}

\author{
  Noah B. Frank\textsuperscript{1}, Joshua L. Pughe-Sanford\textsuperscript{2}, and Samuel J. Grauer\textsuperscript{1,}\thanks{Corresponding author: \href{mailto:sgrauer@psu.edu}{sgrauer@psu.edu}
  }\vspace*{.15em}\\
  {\small \textsuperscript{1}Department of Mechanical Engineering, Pennsylvania State University}\vspace*{-.25em}\\
  {\small \textsuperscript{2}Flatiron Institute, Simons Foundation}}
\date{}

\begin{document}

\maketitle
\setcounter{footnote}{2}
\vspace*{-2em}

% Abstract
\begin{abstract}
\noindent 
Simulations of chaotic systems can only produce high-fidelity trajectories if the initial and boundary conditions are well specified. When these conditions are unknown but measurements are available, variational state estimation can reconstruct a trajectory that is consistent with both the data and the governing equations. A key open question is how many measurements are required for accurate reconstruction, making the full system trajectory observable from sparse data. We establish observability criteria for variational state estimation applied to the Kuramoto--Sivashinsky equation by linking its observability to embedding theory for dissipative dynamical systems. For a system whose attractor lies on an inertial manifold of dimension $d_\mathcal{M}$, we show that $m \geq d_\mathcal{M}$ measurements ensures local observability from an arbitrarily good initial guess, and $m \geq 2d_\mathcal{M} + 1$ implies global observability using a gradient-based smoother since the only critical point on $\mathcal{M}$ is the global minimum. We also analyze optimization-dynamic limitations that persist even when these topological conditions are met, including drift off the manifold, degeneracy of the Hessian, negative curvature, and vanishing gradients. To address these issues, we introduce a robust reconstruction strategy that combines non-convex Newton updates with a novel pseudo-projection step. Numerical simulations of the Kuramoto--Sivashinsky equation validate our analysis and show practical limits of observability for chaotic systems with low-dimensional inertial manifolds.\par\vspace{.5em}

\end{abstract}
\vspace*{2em}

%%% Introduction %%%
\section{Introduction}
\label{sec: introduction}
Chaotic dynamics arise across natural settings and technological systems: emerging in the whorls of turbulent fluid flow at high Reynolds numbers, in the fiery coupling of heat release and acoustic resonance that drives combustion instability, in the tremor of flexible wings at the precipice of flutter, in the jitter of micro-electromechanical resonators, and in the fluctuating frequencies of power grids. Although these systems differ in size and setting, they share a common character: non-linear interactions that span many scales, a marked sensitivity to initial and boundary conditions, and aperiodic long-time dynamics. And yet, for all their apparent unruliness, the governing equations of many chaotic systems are well established, and modern numerical solvers can accurately integrate these systems forward in time from a prescribed starting point. Indeed, where computational resources suffice, scale-resolving simulations, such as direct numerical simulations of turbulent fluid flow, can map out dynamical structures. Such simulations reveal mechanisms of instability \cite{Smith2014} and energy transfer \cite{Marati2004, Bolotnov2010}, they trace out the pathways through which disturbances amplify or decay, and they capture coherent structures \cite{Robinson1991, Smits2011}, along with low-dimensional manifolds that organize long-time behavior. Fully resolved states from high-fidelity simulations also furnish the data from which reduced-order models are extracted \cite{Holmes2012, Linot2023} and from which closures for filtered simulations are calibrated \cite{Duraisamy2019, Argyropoulos2015}, thereby extending our predictive capability into regimes pertinent to engineering design and control.\par

Unfortunately, accurately simulating chaos comes at a steep computational cost. Direct simulations at scales relevant to engineering devices are often prohibitively expensive \cite{Slotnick2014, Cary2021}, and reduced-order or closure models require calibration with data that may not exist in practical regimes. One compromise is to restrict simulations to smaller subdomains \cite{Gronskis2013}, for example, resolving a shock wave--boundary layer interaction on a test article without modeling the entire wind tunnel environment. Although less expensive, such computations are strongly influenced by uncertain and usually unsteady inflow and outflow conditions, heat fluxes, material responses, and so forth. Experimental measurements, by contrast, provide access to the true dynamics of chaotic systems, free from assumptions and compute limitations, but the data are most often sparse, noisy, and indirectly related to the quantities of interest. Fortunately, data assimilation (DA) can blend partial observations from sensors with a system's governing equations, using a numerical solver to produce high-fidelity trajectories that are anchored to real data \cite{Asch2016, Hayase2015, Zaki2025a, Zaki2025b, He2025}. In this paper, we focus on \emph{state estimation}, where the goal is to reconstruct the full system evolution within an observation window (not necessarily to forecast future behavior). With this approach, DA can deliver accurate and dynamically consistent approximations of chaos in regimes where direct simulation is hindered by uncertain initial conditions, boundary conditions, or system parameters.\par

This perspective brings us to a central question: under what circumstances do available measurements provide enough information to permit accurate reconstruction of a chaotic trajectory in state space? In other words, when is the system \emph{observable} via state estimation? For our purposes, a system must satisfy two criteria to be observable. (1)~The data uniquely determine the underlying trajectory. (2)~The reconstruction method can recover that trajectory from those data. Taken together, observability requires that the inverse problem be numerically well posed for the chosen solver.\par

% State estimation
\subsection{Data assimilation methods for state estimation}
\label{sec: introduction: state estimation}
Observability depends on the dynamics of the target system and the available data, including their density, fidelity, and relation to the system state. As indicated above, it also depends on the chosen reconstruction scheme. Broadly, there are three families of DA methods for state estimation: filters, nudging, and smoothers. Filters, such as the Kalman filter \cite{Welch1995} and extensions thereof \cite{Evensen1994}, evolve the governing equations exactly between observation times but introduce discontinuous updates during the analysis steps, so the resultant trajectory does not satisfy the dynamics across the entire assimilation window. Nudging and synchronization observers add feedback terms that drive the modeled system toward the measurements. This can be effective, but these solvers also perturb the true system dynamics since the feedback is not physical \cite{Clark2020, Vela2021, Lalescu2013}. Variational smoothers, by contrast, reconstruct the entire trajectory at once by minimizing a loss functional defined over the full measurement window. They assimilate all the data simultaneously, thereby producing a trajectory that is dynamically consistent (or can be, depending on the solver) while sustaining some discrepancies with the data. Hence, these methods ``smooth out'' said discrepancies by imposing dynamical constraints.\par

Some smoothers enforce the governing equations in an approximate manner via soft penalty terms, while others impose hard constraints embedded in the solver. Physics-informed neural networks (PINNs) \cite{Raissi2019}, for instance, use soft constraints (for the most part). PINNs represent the system's full trajectory with a global model that comprises one or more neural networks. The network parameters are tuned to minimize both measurement error and residuals of the governing equations, and the dynamics are weakly constrained by minimizing these residuals. In experimental fluid mechanics, related DA strategies approximate trajectories of flow states using B-splines \cite{Gesemann2016}, radial basis functions \cite{Casa2013}, or empirical modal bases \cite{Wang2015}, sometimes enforcing linear constraints like mass continuity as exact conditions.\par

A second category of smoothers enforces the discretized system dynamics: the system is parameterized solely by initial and boundary conditions, for instance, and a high-fidelity solver is used to propagate the state forward in time. The loss functional, which compares predicted and experimental observations, is differentiated with respect to the unknown conditions. These conditions are then tuned to minimize the loss via gradient-based optimization. One way to compute these gradients is by solving adjoint equations; adjoint--variational DA \cite{Le1986, Wang2019} is often referred to as ``4DVar'' for unsteady systems in three spatial dimensions. 4DVar solves an adjoint equation that propagates measurement residuals backward in time to yield the gradient. An alternative is ensemble--variational DA, which avoids adjoint equations by estimating gradients statistically from an ensemble of forward model realizations \cite{Liu2008}. For high-dimensional chaotic systems, variational DA has been shown to be accurate and efficient \cite{Mons2016, Zaki2025b}, since it provides exact gradients and its computational cost does not scale with the dimension of the control vector---in this case, initial and boundary conditions---which becomes very large in 4DVar state estimation problems.\par

Our goal in this work is to assess the fundamental limits of observability, so we seek to minimize solver-induced biases and avoid data--physics trade-offs, where possible. Since variational state estimation can enforce the governing equations exactly, it provides a natural framework to probe observability limits.\par

% Embedding theory
\subsection{Application of embedding theory to observability in state estimation}
\label{sec: introduction: embedding theory}
Although many dissipative chaotic systems, ranging from the relatively simple Kuramoto--Sivashinsky (KS) equation to fluid turbulence governed by the three-dimensional (3D) Navier--Stokes equations, formally evolve in an infinite-dimensional state space, their long-time dynamics are expected to collapse onto a finite-dimensional invariant subset of state space know as the global attractor, which we denote by $\mathcal{A}$ \cite{Zelik2014}. The box counting dimension of this attractor, $d_\mathcal{A}$, quantifies the system's \emph{effective degrees of freedom} and provides a measure of its complexity. It naturally follows that $d_\mathcal{A}$ should influence the number of measurements required to uniquely determine the system state.\par

Embedding theory allows us to make this notion precise. It considers mappings of the form $\Phi : \mathcal{A} \to \mathbb{R}^m$, which takes a point on the attractor to an $m$-dimensional vector of observations $\boldsymbol{y} \in \mathbb{R}^m$. When the observations are sufficiently rich, $\Phi$ becomes an \emph{embedding}, meaning that it is a smooth, one-to-one, and invertible mapping from points on the attractor to its image in $\mathbb{R}^m$, and it has a smooth inverse. If $\Phi $ is an embedding, then $\Phi(\mathcal{A})$ is topologically equivalent to the attractor, so the geometry of the dynamics can be unfolded in measurement space without self-intersections \cite{Kugiumtzis1996}. Takens' embedding theorem and its extensions \cite{Deyle2011} formalize this principle, showing that if the observation dimension $m$ exceeds twice the attractor dimension ($m > 2d_\mathcal{A}$), then $\Phi$ is almost always an embedding. (These theorems holds under assumptions that we shall discuss later.) These results underpin the field of \emph{state space reconstruction},\footnote{State space reconstruction recovers an attractor by assembling delay-coordinate vectors $\boldsymbol{y}(t), \boldsymbol{y}(t-\tau), \boldsymbol{y}(t-2\tau), \dots$ from one or more observed time series, typically to examine the attractor geometry or to estimate dynamical invariants of the system, like its Lyapunov exponents or fractal dimension, without explicitly enforcing the governing equations. State estimation in DA, by contrast, seeks to reconstruct the system's trajectory in the original state space. This work is concerned with the latter.} which leverages the topological equivalence of $\mathcal{A}$ and $\Phi(\mathcal{A})$ to determine features of chaotic dynamics directly from measured data.\par

Although embedding theory has been developed primarily for state space reconstruction, we suggest that its insights can likewise inform gradient-based state estimation. An embedding, by definition, is a smooth injective map with a smooth inverse $\Phi^{-1}$, which implies that the initial state of the system, and thus its full trajectory over an observation window, can be uniquely identified from the data in $\boldsymbol{y}$. Embedding theorems specify when such a correspondence exists in principle, linking the number of measurements $m$ to the attractor dimension $d_\mathcal{A}$. This is the first component of observability that we introduced above. What these theorems do not provide is a practical means of recovering $\Phi^{-1}$, nor do they address the complications introduced by noise, limited observation windows, sensor placement, or model error \cite{Casdagli1991}. These gaps motivate our investigation. We ask under what measurement conditions, and with which DA schemes, the observability conditions suggested by embedding theory can be achieved in practice.\par

% Roadmap
\subsection{Roadmap to the paper}
\label{sec: introduction: roadmap}
To address the observability question posed above, we require a model system that is both chaotic and computationally accessible. The KS equation with periodic boundary conditions provides such a testbed. It is a non-linear partial differential equation, often regarded as a minimal model of spatiotemporal chaos \cite{Wittenberg1998}, that transitions from intermittent disorder to fully developed chaos with increasing domain lengths. Despite its simplicity, the KS equation exhibits many hallmarks of more complex systems of engineering interest, including multiscale interactions and an energy cascade. Its long-term dynamics are well characterized in the literature \cite{Ding2016, Linot2020, Linot2022, Zeng2024, Yang2009, Takeuchi2011}, which provides benchmarks for testing predictions from embedding theory. At the same time, the modest computational cost of KS simulations enables systematic studies of measurement configurations and optimization strategies that would be prohibitive in higher-dimensional systems such as 3D fluid turbulence.\par

We apply embedding theory to the problem of observability in DA-based state estimation, with the goal of testing whether the reconstruction limits suggested by theory can be realized in practice. Our aim is not to recover invariants of chaotic attractors, per se, as in state space reconstruction, but to reconstruct full system trajectories that resolve the underlying fields, thereby enabling physical interpretation and modeling. To this end, we adopt a variational state estimation method in which the state is parameterized by its initial condition and marched forward by a high-fidelity solver. Because the KS equation is one-dimensional (1D) in space and evolves in time, we refer to the approach as ``2DVar.'' We use 2DVar reconstructions to examine how accuracy varies with the density, spacing, and repetition rate of ``sensors,'' as well as on the duration of the observation window. First, we analyze our results in relation to predictions from embedding theory, expressed through a manifold dimension that bounds $d_\mathcal{A}$. We refer to dependencies on the attractor dimension as \emph{degree-of-freedom effects}. After that, we examine how the topology of the loss landscape and the behavior of gradient-based optimizers influence the outcome of state estimation. We call this influence \emph{optimization-dynamic effects}.\par

In the remainder of this paper, Sec.~\ref{sec: state estimation} introduces the KS equation as well as our numerical framework for forward and inverse computations. Section~\ref{sec: cases} lays out our procedure for generating test cases and presents sample reconstructions. Section~\ref{sec: DoF} investigates how reconstruction accuracy depends on the sensor network, relating these results to the system's degrees of freedom and embedding-based observability criteria. Section~\ref{sec: optimization} examines the role of optimization dynamics, showing how key features of the loss landscape influence convergence, and demonstrating that most difficulties arise from vanishing gradients and negative curvature rather than spurious local minima. Finally, Sec.~\ref{sec: conclusions} presents our conclusions and discusses some broader implications of our work for higher-dimensional systems. Pertinent details on numerical methods and supporting derivations are provided in the appendices.\par

%%% KSE_DA %%%
\section{Variational state estimation for Kuramoto--Sivashinsky systems}
\label{sec: state estimation}
We begin by introducing the model system and reconstruction framework used throughout the paper, starting with the KS equation, its key properties, and our numerical solver. Next, we present our 2DVar formulation and examine how chaos affects gradient computations. We then review first- and second-order optimization strategies and introduce a novel stabilization method for adjoint marching termed \emph{pseudo-projection}.\par

% KSE overview
\subsection{Formulation and dynamics of the Kuramoto--Sivashinsky equation}
\label{sec: state estimation: KS}
In its derivative-form, the KS equation reads
\begin{equation}
\label{equ: KSE}
    \frac{\partial u}{\partial t}
    = -\underbrace{\frac{\partial^2 u}{\partial x^2}}_\text{(I)}
    - \underbrace{\frac{\partial^4 u}{\partial x^4}}_\text{(II)}
    - \underbrace{u \frac{\partial u}{\partial x}}_\text{(III)},
\end{equation}
for positions $x \in [-L/2, L/2]$, where $L$ is the domain length, and times $t \in [0,\infty)$. In practice, we consider finite observation windows $[t_0, t_0 + T]$, where $t_0$ is an arbitrary start time and $T$ is the window length. Observation times are relative to $t_0$, with $t \in [0, T]$. Periodic boundary conditions are imposed,
\begin{equation}
    \left(\frac{\partial^a u}{\partial x^a}\right)_{x=-L/2} =
    \left(\frac{\partial^a u}{\partial x^a}\right)_{x=L/2}, \quad
    \forall \ a \in \mathbb{N}_0,
\end{equation}
and the dynamics are fully specified by the initial condition $u(x,0)$.\par

Terms (I)–(III) can be analyzed in Fourier space, where periodic solutions admit the expansion
\begin{equation}
    \label{equ: Fourier}
    u(x, t) = \sum_{j  \in \mathbb{Z}} \widehat{u}_j(t) \,e^{\mathrm{i} k_j x},
\end{equation}
with wavenumbers $k_j = 2\pi j / L$ and Fourier coefficients $\widehat{u}_j$. Substituting this series into Eq.~\eqref{equ: KSE} yields
\begin{equation}
    \label{equ: KSE Fourier}
    \frac{\partial \widehat{u}_j}{\partial t} = (k_j^2 - k_j^4) \,\widehat{u}_j - \frac{\mathrm{i} k_j}{2}\sum_{i \in \mathbb{Z}} \widehat{u}_i \widehat{u}_{j-i}.
\end{equation}
Here, $(k_j^2 - k_j^4) \,\widehat{u}_n$ combines the energy-producing anti-diffusion and the energy-dissipating hyper-diffusion terms, i.e., (I) and (II). It amplifies low-wavenumber modes ($|k_j| < 1$) and damps high-wavenumber ones ($|k_j| > 1$), with maximum growth near $k_\mathrm{crit} = \pm 1/\sqrt{2}$. The conservative convective term (III) redistributes energy across wavenumbers, coupling the stable and unstable modes and orchestrating the balance between the production, dissipation, and redistribution of energy that gives rise to chaos.\par

For all domain lengths $L$, solutions of the KS equation have been proven to rapidly approach a smooth, finite-dimensional manifold $\mathcal{M}$, called the inertial manifold (IM) \cite{Foias1986}, which contains the global attractor $\mathcal{A} \subset \mathcal{M}$. The manifold is invariant under the system dynamics, so long-time trajectories of $u$ are embedded in $\mathcal{M}$ \cite{Jones1996}. The dimension of this manifold, $d_\mathcal{M}$, has been extensively studied and provides an upper bound for the attractor's box counting dimension, $d_\mathcal{A} \leq d_\mathcal{M}$, where computing $d_\mathcal{A}$ is often intractable. Thus, $d_\mathcal{M}$ serves as a rigorous measure of the system's effective degrees of freedom, and embedding theory allows us to connect it to the number of measurements needed for reliable state estimation.\par

\begin{figure}[htb!]
    \center
    \includegraphics[width=0.65\textwidth]{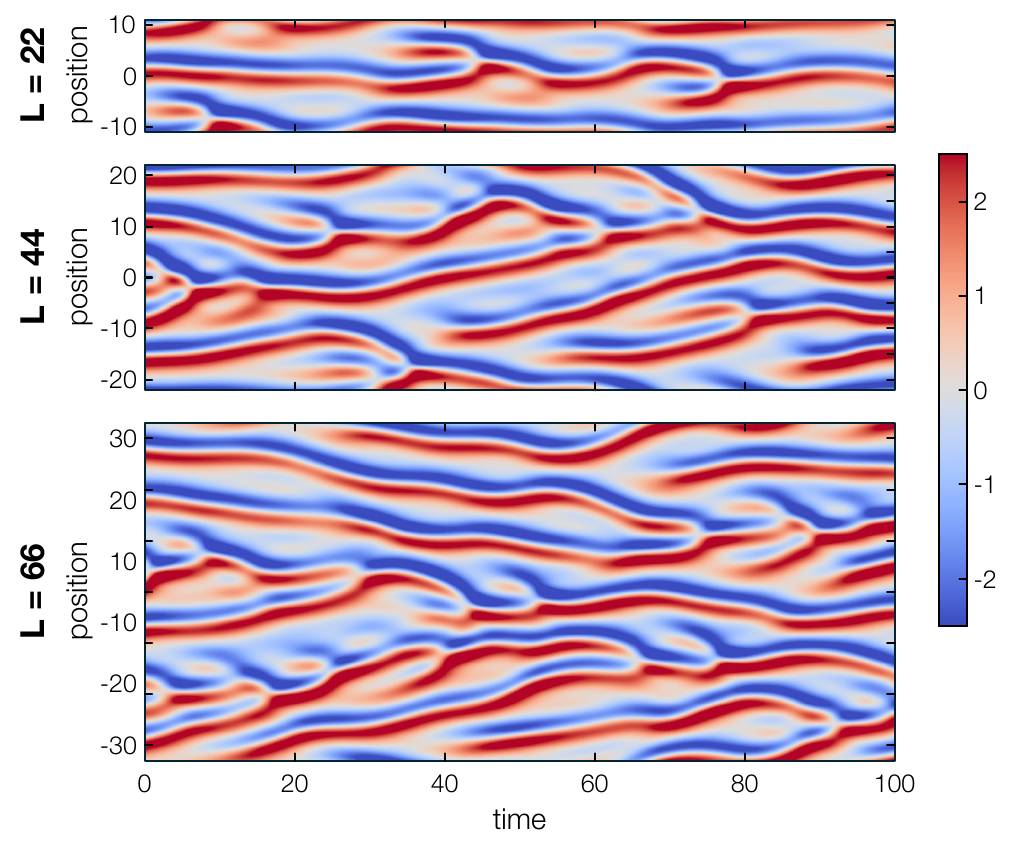}
    \caption{Representative trajectories of Kuramoto--Sivashinsky systems for domain lengths $L = 22$ (top), $L = 44$ (middle), and $L = 66$ (bottom). All cases exhibit meandering streaks whose spatial and temporal complexity increases with $L$, reflecting the additional active degrees of freedom.}
    \label{fig: sample solutions}
\end{figure}

Forward solutions to the KS equation are obtained with a custom solver implemented in JAX, with full details provided in Appendix~\ref{app: numerics}. The solver employs a uniform grid with a Fourier pseudo-spectral discretization in space and a fourth-order Runge--Kutta exponential time-differencing scheme \cite{Cox2002} to handle the stiff linear terms \cite{Cvitanovic2010}. JAX also provides efficient gradient computations via automatic differentiation (AD), which we exploit in our 2DVar formulation. For $L \in \{22, 44\}$ we use $64$ grid points with a time step of $0.1$, while for $L = 66$ we use $72$ points with a time step of $0.05$. These discretizations are consistent with prior studies on KS systems \cite{Linot2020, Cvitanovic2010}, and the Kaplan--Yorke dimension computed with our solver agrees with the results of Edson et al. \cite{Edson2019}.\par

\begin{table}[htb!]
    \caption{Inertial manifold dimension and leading Lyapunov exponent for different domain lengths.}
    \renewcommand{\arraystretch}{1.25}
    \centering\vspace*{.3em}
    \begin{tabular}{c c c}
        \hline\hline
        \multicolumn{1}{c}{\bf Domain Length $\boldsymbol{L}$} &
        \multicolumn{1}{c}{\bf IM Dimension $\boldsymbol{d_\mathcal{M}}$} &
        \multicolumn{1}{c}{\bf Leading Lyapynov Exponent $\boldsymbol{\ell_1}$} \\
        \hline
        22 & 8 & 0.05 \\
        44 & 18 & 0.083 \\
        66 & 28 & 0.087 \\
        \hline\hline
    \end{tabular}
    \label{tab: invariants}
\end{table}

Numerical DA experiments are carried out in domains of length $L \in \{22, 44, 66\}$, each of which sustains chaotic dynamics \cite{Cvitanovic2010}. As $L$ grows larger, KS systems display behavior of increasing complexity: $L = 22$ lies just beyond the onset of structurally stable chaos \cite{Cvitanovic2010}, while $L = 66$ exhibits strongly chaotic behavior. Representative trajectories are shown in Fig.~\ref{fig: sample solutions}, characterized by undulating waveforms that randomly appear, drift and mingle across the domain, and merge together.\par

\begin{figure}[htb!]
    \center\includegraphics[height=6.5cm]{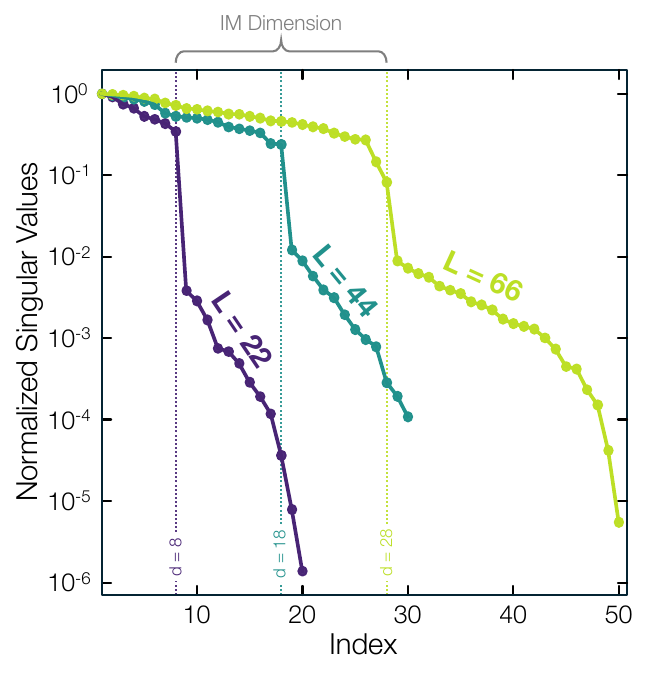}
    \caption{Autoencoder-based estimates of $d_\mathcal{M}$ for $L \in \{22, 44, 66\}$. Singular values of the centered latent state data matrix are shown, with vertical dashed lines marking the inferred IM dimension identified by the sharp drop in eigenvalues.}
    \label{fig: IM dimension}
\end{figure}

To verify our solver, we benchmark our solutions against known chaotic invariants. Table~\ref{tab: invariants} reports the IM dimension and leading Lyapunov exponent for each domain length. Lyapunov spectra are computed from long-time simulations using the QR method \cite{Benettin1980}, with additional details provided in Appendix~\ref{app: Lyapunov}. Estimates of $d_\mathcal{M}$ are obtained using the autoencoder-based approach of Zeng et al. \cite{Zeng2024}. An autoencoder couples an encoder $\mathsf{E} : \mathcal{M} \to \mathcal{L}$ to a decoder $\mathsf{D} : \mathcal{L} \to \mathcal{M}$, where the composite map $\mathsf{A} = \mathsf{D} \circ \mathsf{E}$ learns an identity function on $\mathcal{M}$ and $\mathcal{L}$ is a low-dimensional latent space of dimension $d_\mathcal{L}$. The encoder and decoder are jointly trained to minimize reconstruction error, such that $\mathsf{A}$ learns to represent states on $\mathcal{M}$ in the compressed latent space $\mathcal{L}$. We set $d_\mathcal{L}$ conservatively so that $d_\mathcal{M} \leq d_\mathcal{L} < n$, where $n$ is the dimension of the discrete state vector. The architecture of $\mathsf{A}$ is designed to promote models which only use a low-dimensional subset of the  $d_\mathcal{L}$-dimensional latent space. After training, the autoencoder provides mappings to and from a low-dimensional embedding of KS states. A principal component analysis (PCA) of the encoded states in $\mathcal{L}$ yields a covariance matrix whose effective rank is taken as an estimate of $d_\mathcal{M}$. Figure~\ref{fig: IM dimension} shows singular values of the centered data matrix and inferred dimensions for the $L = 22,\ 44,\ \text{and}\ 66$ domains. These estimates of $d_\mathcal{M}$, computed with our solver and autoencoders, are consistent with rigorous analyses \cite{Ding2016}, physical-mode counts \cite{Yang2009, Takeuchi2011}, and previous autoencoder-based studies \cite{Linot2020, Linot2022, Zeng2024}. Further information on our autoencoder architectures and training procedures is provided in Appendix~\ref{app: AE}.\par

% Overview of adjoint DA
\subsection{Variational state estimation}
\label{sec: state estimation: adjoint}
We consider state estimation for a KS system from sparse spatio-temporal point measurements. Our formulation broadly follows that of Protas et al.~\cite{Protas2004}; interested readers are also directed to Jardak et al.~\cite{Jardak2010} for a discussion of ensemble smoothers applied to KS systems. We present the state estimation problem in discrete form, retaining only the definitions and equations needed to specify the inverse problem and describe its gradient-based optimization. Derivations of adjoint systems are deferred to Appendix~\ref{app: adjoint}.\par

The discrete state is represented by a vector $\boldsymbol{u}_k \in \mathbb{R}^{n}$ containing the solution $u$ at $n$ uniformly spaced spatial nodes and at time $k \Delta t$. The time index satisfies $k \in \mathcal{K} = \{0, \dots, K\}$, where $K = T/\Delta t$ is the final step of the rollout. The system is advanced by a numerical solver,
\begin{equation}
    \boldsymbol{u}_{k+1} = \mathsf{f}(\boldsymbol{u}_k),
\end{equation}
where $\mathsf{f}$ represents one time step of the forward solver described in Appendix~\ref{app: numerics}.\par

We consider a limited-data problem in which $m$ scalar observations are available, collected in the vector $\boldsymbol{y} \in \mathbb{R}^m$, with $m \ll nK$. Each entry of $\boldsymbol{y}$ corresponds to a single observation of the system at a particular point in space and time. More generally, the $i$th observation is represented by a smooth observation operator $\mathsf{h}_i(\boldsymbol{u}_0) : \mathbb{R}^n \to \mathbb{R}$ of the form
\begin{equation}
    y_i = \mathsf{h}_i(\boldsymbol{u}_0) =
    \mathsf{g}_i\!
    \left[
        \mathsf{f}^{\tau_i}(\boldsymbol{u}_0)
    \right] \!,
\end{equation}
where $\mathsf{g}_i$ is a smooth scalar-valued function, $\mathsf{f}^{\tau_i}$ denotes $\tau_i$ successive applications of the discrete flow map, and $\tau_i$ is the discrete time index associated with the $i$th observation. An observation is said to be ``unlagged'' when $\tau_i = 0$, in which case $\mathsf{f}^{\tau_i}$ is the identity map.\par

In this paper, we restrict our attention to point measurements, for which
\begin{equation}
    y_i = \boldsymbol{e}_{\ell_i}^\top \boldsymbol{u}_{\tau_i}
    = \boldsymbol{e}_{\ell_i}^\top \mathsf{f}^{\tau_i}(\boldsymbol{u}_0),
\end{equation}
where $\boldsymbol{e}_{\ell_i}$ is the $\ell_i$th standard basis vector. Spatial sensor locations are selected as $x \in \mathcal{X}$, and measurement times are selected as $t \in \mathcal{T}$. A distinct observation coordinate is denoted by $(x, t)_i \in \mathcal{X} \times \mathcal{T}$ and is identified by the index $i \in \mathcal{I}$, where $m = |\mathcal{I}|$. These coordinates correspond to discrete indices $\ell_i = x_i/\Delta x$ and $\tau_i = t_i/\Delta t$. We assume that all observation points coincide with grid nodes and time steps so that $\ell_i$ and $\tau_i$ are integers.\par

\begin{figure}[htb!]
    \centering
    \includegraphics[height=3.75cm]{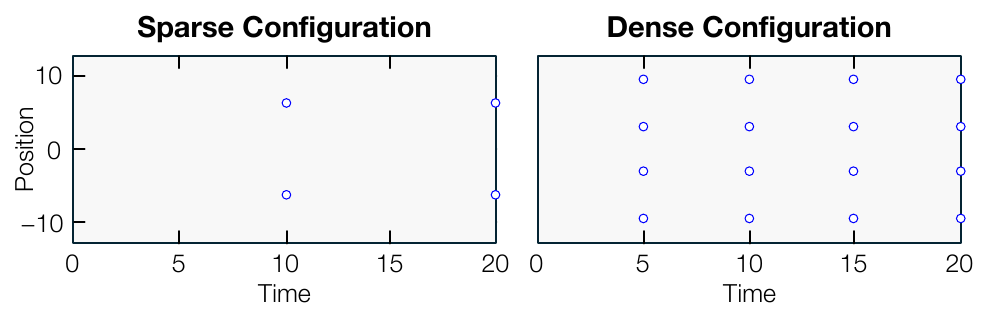}
    \caption{Exemplary measurement configurations in the $L = 22$ domain. The left panel shows a sparse layout with two spatial sensors and two observation times; the right panel shows a denser configuration with four spatial sensors and four observation times.}
    \label{fig: sensor layout}
\end{figure}

In our numerical experiments, the observation set is defined by placing sensors uniformly in space and sampling them at regular time intervals. The spatial locations are evenly distributed, with $m_x$ sensors separated by $\Delta x = L/m_x$ and centered within the domain. Each sensor records $m_t$ samples at a constant rate, with observations beginning at $t_0 + \Delta t$ and including the final time $t_0 + T$, where $\Delta t = T/m_t$. Thus, the initial state is always excluded, the final time is always observed, and the total number of measurements is $m = m_x \,m_t$. Sample configurations are shown in Fig.~\ref{fig: sensor layout} for the $L = 22$ domain and a time horizon of $T = 20$. The left panel shows a sparse layout with $m = 4$, and the right panel shows a denser layout with $m = 16$.\par

With these $m$ measurements, the objective is to reconstruct the initial condition that produced the observed trajectory. To this end, we parameterize the observer system's initial condition using $p$ Fourier coefficients $\boldsymbol{\theta} \in \mathbb{C}^p$, such that
\begin{equation}
    \boldsymbol{u}_{\boldsymbol{\theta}} = \mathsf{F}^{-1}(\boldsymbol{\theta})
    \quad\Longleftrightarrow\quad
    \boldsymbol{\theta} = \mathsf{F}(\boldsymbol{u}_{\boldsymbol{\theta}}),
\end{equation}
where $\mathsf{F}$ and $\mathsf{F}^{-1}$ are the discrete Fourier transform and its inverse. We use 15, 20, and 25 Fourier modes for the $L = 22$, $44$, and $66$ domains, respectively. The predicted observation at the $i$th measurement coordinate is thus
\begin{equation}
    \mathsf{h}_i(\boldsymbol{\theta})
    = \boldsymbol{e}_{\ell_i}^\top \,\mathsf{f}^{\,\tau_i}\!
    \left[
        \mathsf{F}^{-1}(\boldsymbol{\theta})
    \right] \!.
\end{equation}
In the absence of measurement noise or prior information about the initial condition, the variational DA problem reduces to minimization of the mean squared error between the reference and observer measurements. This corresponds to minimization of the loss functional
\begin{equation}
    \label{equ: discrete loss}
    \mathscr{J}(\boldsymbol{\theta}) = \frac{1}{m} \sum_{i \in \mathcal{I}}
    \left[\mathsf{h}_i(\boldsymbol{\theta}) - y_i\right]^2.
\end{equation}
The initial state is reconstructed by solving
\begin{equation}
    \widetilde{\boldsymbol{\theta}} = \arg\min_{\boldsymbol{\theta}} \mathscr{J}(\boldsymbol{\theta}),
\end{equation}
with the corresponding reconstructed initial state given by $\boldsymbol{u}_{\widetilde{\boldsymbol{\theta}}, 0} = \mathsf{F}^{-1}(\widetilde{\boldsymbol{\theta}})$. For notational convenience, we suppress the time index and write the initial observer state $\boldsymbol{u}_{\boldsymbol{\theta}, 0}$ as $\boldsymbol{u}_{\boldsymbol{\theta}}$.\par

We solve this optimization problem using gradient-based methods. Our solver is implemented in JAX, which enables AD of full rollouts to compute gradients and Hessians of $\mathscr{J}$. Although we also derived and implemented the discrete adjoint equations, we use AD for the results reported in this work because it is straightforward to implement and computationally efficient for the present problem sizes. When AD is unavailable, or when memory constraints limit its use, the same derivatives can be obtained from the discrete adjoint systems given in Appendix~\ref{app: adjoint}.\par

% Adjoint chaos
\subsubsection{Implications of chaos for the inverse problem}
\label{sec: state estimation: adjoint: chaos}
Numerical simulations of chaotic systems necessarily diverge from the true system trajectory over long time horizons due to the exponential growth of errors. The adjoint system inherits the Lyapunov spectrum of the forward dynamics \cite{Chandramoorthy2019, Zaki2021}, so gradients of $\mathscr{J}$ are vulnerable to the amplification of measurement noise and numerical errors accumulated during backward integration. Consequently, optimizing a single trajectory over long time horizons is not feasible \cite{Zaki2025b}. A common strategy is to restrict the assimilation window to be on the order of the Lyapunov timescale, $T_\ell = 1/\ell_1$, where $\ell_1$ is the leading Lyapunov exponent \cite{Li2020, Chandramouli2020}. For measurements spanning a total duration $T > T_\ell$, the problem is divided into $\mathrm{ceil}(T/T_\ell)$ segments of length $T_\ell$ or less \cite{Chandramouli2020}. In order to speed up convergence in a multi-window reconstruction, the terminal condition from one segment may be used as the initial guess for the next \cite{Wang2025}. In this paper, we follow this practice and take the assimilation windows to be of duration $T_\ell$.\par

% Optimization methods
\subsection{Optimization methods}
\label{sec: state estimation: optimization}
Given exact gradients and Hessians of the loss functional, i.e., computed by AD or adjoint equations, the performance of gradient-based state estimation depends on the optimization algorithm. The choice of optimizer is pivotal in chaotic systems, where ill-conditioning and instability are endemic.\par

It is often assumed that optimization is difficult because descent methods become trapped in spurious local minima \cite{Dauphin2014, Paternain2019}. However, in many high dimensional problems, such as dictionary learning \cite{Sun2016}, tensor decomposition \cite{Ge2015}, matrix completion \cite{Ge2016}, and training certain (restrictive) classes of neural networks \cite{Kawaguchi2016}, \emph{all} minima are global minima with the same loss. Theory for random high-dimensional error surfaces, which are good surrogates for practical loss landscapes, likewise suggests that nearly all critical points of high loss are saddle points rather than minima \cite{Bray2007, Fyodorov2007}. Ergo, the main challenge in most high-dimensional, gradient-based state estimation problems is escaping saddle points \cite{Dauphin2014, Baldi1989, Saxe2013}. As a corollary, any local minimum can be regarded as a satisfactory solution. Similar challenges arise in variational state estimation, namely, effective handling of negative curvature and achieving convergence despite an inherently ill-conditioned Hessian.\par

We consider four minimizers in this work: vanilla gradient descent, Newton's method, a quasi-Newton algorithm, and a regularized variant of Newton's method. The first three illustrate common failure modes in variational state estimation, whilst the last mitigates these pathologies. Some other methods, not tested in this work, are briefly discussed in Sec.~\ref{sec: cases: optimizers}. Throughout our discussion, the control vector at iteration $k$ is denoted $\boldsymbol{\theta}_k$; it is initialized at $\boldsymbol{\theta}_0$ and iteratively updated with increasing $k$. The gradient and Hessian of the loss functional are
\begin{equation}
    \boldsymbol{g} = \left(\frac{\partial \mathscr{J}}{\partial \boldsymbol{\theta}}\right)^\top
    \quad\text{and}\quad
    \boldsymbol{H} = \frac{\partial^2 \mathscr{J}}{\partial \boldsymbol{\theta}^2},
\end{equation}
where $\boldsymbol{g} \in \mathbb{R}^n$ and $\boldsymbol{H} \in \mathbb{R}^{n \times n}$ is symmetric. Adjoint expressions for $\boldsymbol{g}$ and $\boldsymbol{H}$ are provided in Eqs.~\eqref{equ: loss gradient} and \eqref{equ: adjoint Hessian}, respectively.\par

% Gradient descent
\subsubsection{Vanilla gradient descent}
\label{sec: state estimation: optimization: gradient descent}
Gradient descent updates the control vector along the steepest slope of the loss landscape,
\begin{equation}
    \boldsymbol{\theta}_{k+1} = \boldsymbol{\theta}_k - \eta \boldsymbol{g},
\end{equation}
where $\eta > 0$ is the step size, which may vary with $k$. Unfortunately, progress can slow to a crawl when the Hessian is ill-conditioned \cite{Paternain2019}. To see this, consider a quadratic loss with ordered Hessian eigenvalues $\lambda_1 \geq \dots \geq \lambda_n > 0$. The condition number of $\boldsymbol{H}$ is $\kappa_{\boldsymbol{H}} = \lambda_1/\lambda_n$. Using the optimal step size, given by $\eta = 2/(\lambda_1 + \lambda_n)$, progress near a critical point $\boldsymbol{\theta}^*$ satisfies
\begin{equation}
    \left\lVert \boldsymbol{\theta}_{k+1} - \boldsymbol{\theta}^* \right\rVert_2
    < \frac{\kappa_{\boldsymbol{H}} - 1}{\kappa_{\boldsymbol{H}} + 1}
    \left\lVert \boldsymbol{\theta}_k - \boldsymbol{\theta}^* \right\rVert_2,
\end{equation}
where $\lVert \cdot \rVert_2$ is the Euclidean norm \cite{Saarinen1993}. Thus, when $\kappa_{\boldsymbol{H}} \gg 1$, the contraction factor approaches unity, updates shrink extremely slowly, and the optimizer lingers near $\boldsymbol{\theta}^*$.\par

% Newton
\subsubsection{Newton's method}
\label{sec: state estimation: optimization: Newton}
Gradients can be rescaled using local curvature information to deal with ill-conditioned Hessians. Newton's method, which achieves quadratic local convergence when $\boldsymbol{H}$ is full rank \cite{Saarinen1993}, modifies gradient descent by taking smaller steps along directions of high curvature and larger steps along directions of low curvature,
\begin{equation}
    \boldsymbol{\theta}_{k+1} = \boldsymbol{\theta}_k - \eta \boldsymbol{H}^{-1} \boldsymbol{g}.
\end{equation}
When $\boldsymbol{H}$ has negative eigenvalues, however, the Newton step can move uphill, converging to a nearby saddle point or local maximum. Moreover, if $\boldsymbol{H}$ is singular or nearly so, regularization or other modifications are required to approximate $\boldsymbol{H}^{-1}$.\par

% Quasi-Newton
\subsubsection{Quasi-Newton methods}
\label{sec: state estimation: optimization: BFGS}
Quasi-Newton methods are employed when computing, storing, or inverting the full Hessian is too costly. Instead, an approximation to the inverse Hessian $\boldsymbol{B}_k^{-1} \approx \boldsymbol{H}^{-1}$ is built from past gradients and iterates, so that the update step is
\begin{equation}
    \boldsymbol{\theta}_{k+1} = \boldsymbol{\theta}_k - \eta \boldsymbol{B}_k^{-1}\boldsymbol{g}.
\end{equation}
Like in Newton's method, $\boldsymbol{B}_k^{-1}$ is meant to accelerate convergence in locally convex regions by rescaling gradients according to the estimated curvature at $\boldsymbol{\theta}_k$. However, the theoretical foundations of quasi-Newton methods break down in regions with strong negative curvature \cite{Dauphin2014}. The most widely used variant, i.e., the Broyden--Fletcher--Goldfarb--Shanno (BFGS) algorithm, builds a positive definite approximation $\boldsymbol{B}_k^{-1}$ so long as the curvature condition is satisfied. At points of negative curvature this condition can fail, producing a non-–positive definite or poorly conditioned update. A common remedy is to reset $\boldsymbol{B}_k^{-1} = \boldsymbol{I}$, thereby reverting to a gradient descent step at iteration $k$. As a result, BFGS can perform poorly in regions where negative curvature is prevalent. We use BFGS as our representative quasi-Newton method.\par

% Regularized Newton
\subsubsection{Regularized Newton methods}
\label{sec: state estimation: optimization: NCN}
Several modifications to Newton's method have been proposed to address ill-conditioning of the Hessian and negative eigenvalues. A common strategy is to take the absolute value and\slash or threshold the eigenvalues of $\boldsymbol{H}$ prior to inversion \cite{Greenstadt1967, Gould1998, Dauphin2014, Paternain2019}. We adopt one such technique called the non-convex Newton (NCN) method \cite{Paternain2019}, which conditions the gradient using the ``positive definite truncated inverse Hessian,'' denoted $|\boldsymbol{H}|^{-1}$. NCN steps are given by
\begin{equation}
    \boldsymbol{\theta}_{k+1} = \boldsymbol{\theta}_k - \eta
    |\boldsymbol{H}|^{-1} \boldsymbol{g}.
\end{equation}
To compute $|\boldsymbol{H}|^{-1}$, we perform an eigendecomposition of the Hessian,
\begin{equation}
    \boldsymbol{H} = \boldsymbol{Q\Lambda Q}^\top,
\end{equation}
where $\boldsymbol{\Lambda} = \mathrm{diag}(\lambda_1, \dots, \lambda_n)$ contains the ordered eigenvalues and $\boldsymbol{Q}$ is an orthonormal matrix that comprises the corresponding eigenvectors. The eigenvalues are regularized:
\begin{equation}
    \lambda_i^\prime = \max(|\lambda_i|, \delta),
\end{equation}
where $\delta \in (0,\lambda_1]$ is a threshold, which yields the diagonal matrix $|\boldsymbol{\Lambda}| = \mathrm{diag}(\lambda_1', \dots, \lambda_n')$. Note that $\delta > 0$ ensures that $|\boldsymbol{H}|$ is full rank. The regularized inverse is
\begin{equation}
    \label{equ: NCN Hessian}
    |\boldsymbol{H}|^{-1} =
    \boldsymbol{Q} |\boldsymbol{\Lambda}|^{-1} \boldsymbol{Q}^\top.
\end{equation}
The threshold $\delta$ controls the condition number of $|\boldsymbol{H}|^{-1}$, which equals $\lambda_1 / \delta$.\par

Unlike Newton's method, NCN can descend along directions of negative curvature and is proven to escape saddle points exponentially fast \cite{Paternain2019}. Furthermore, because $|\boldsymbol{H}|^{-1}$ is inherently positive definite, NCN can always decrease the loss function for a sufficiently small step, provided that $\boldsymbol{g} \neq \mathbf{0}$. Since NCN can handle ill-conditioned Hessians and negative concavity, and because the Hessian is easy to compute for our KS state estimation problem, NCN is the default optimizer in this work.\par

Selection of the threshold $\delta$ and step size $\eta$ are key considerations when using NCN. We first set $\delta = \lambda_1 / \kappa_{|\boldsymbol{H}|}$, where $\kappa_{|\boldsymbol{H}|}$ is the targeted condition number, and we then determine $\eta$ through a standard backtracking line search that satisfies the Armijo condition. Setting $\kappa_{|\boldsymbol{H}|}$ to 1 removes curvature-based scaling from the step direction, i.e., $|\boldsymbol{H}|^{-1} \boldsymbol{g} = \lambda_1^{-1} \boldsymbol{g}$, which is equivalent to gradient descent (i.e., robust to noise in $\boldsymbol{g}$ but slow and prone to stagnation near critical points). Conversely, as $\delta \to 0$, we have $\kappa_{|\boldsymbol{H}|} \to \kappa_{\boldsymbol{H}}$ and all components of the Hessian are retained. However, when the Hessian is ill-conditioned, NCN steps amplify gradient components aligned with directions of low curvature, which we observe are typically associated with high-wavenumber components in $\boldsymbol{u}_0$ (see Sec.~\ref{sec: DoF}). Such components are usually non-physical because the hyper-diffusion term (II) in Eq.~\eqref{equ: KSE} rapidly damps high-wavenumber content and the true state lies near the attractor (by assumption), which primarily contains low-wavenumber modes. Selection of $\kappa_{|\boldsymbol{H}|}$ is thus a balancing act: it must be large enough to accelerate convergence by exploiting curvature information but small enough to avoid amplifying non-physical gradient components. In practice, we constrain $\kappa_{|\boldsymbol{H}|}$ to $\{10^3, 10^5\}$ to simplify the regularization procedure.\par

We set $\kappa_{|\boldsymbol{H}|} = 10^3$ at the start of each optimization, which biases updates toward states dominated by low-wavenumber modes. However, if $\boldsymbol{g}^\top |\boldsymbol{H}|^{-1} \boldsymbol{g}$ becomes much smaller than the objective loss, optimization effectively stalls (see Sec.~\ref{sec: optimization: convergence}). To address this, whenever $\mathscr{J}^{-1} (\boldsymbol{g}^\top |\boldsymbol{H}|^{-1} \boldsymbol{g}) < 0.01$, we increase $\kappa_{|\boldsymbol{H}|}$ to $10^5$, enabling larger steps along ``higher-wavenumber directions.'' This adjustment typically gets the optimizer unstuck and reduces the loss, though it also tends to increase early-time reconstruction errors.\par

% Pseudo-projection
\subsection{Pseudo-projection}
\label{sec: state estimation: projection}
Most initial conditions in state space do not lie on the attractor, and many points off the attractor converge toward trajectories that are nearly indistinguishable from those on it. As a result, when optimizing the initial condition, the estimate can drift away from the attractor along directions that yield similar trajectories but contain erroneous early transients. This tendency is exacerbated when $\kappa_{|\boldsymbol{H}|}$ is large, since curvature scaling amplifies gradient components aligned with directions of low curvature (typically associated with high-wavenumber content in $\boldsymbol{u}_0$, as discussed above). We refer to this behavior as a blow-up of the initial condition, since values of $\boldsymbol{u}_{\boldsymbol{\theta}}$ can become unphysically large along these directions.\par

Because all trajectories approach the attractor exponentially fast, it is reasonable to assume that the true initial condition lies on the attractor or else very near it. To incorporate this knowledge into our state estimation algorithm and to prevent blow-up, we introduce \emph{pseudo-projection}. This operation projects a state toward the attractor by integrating the governing equations forward in time for a short duration. Pseudo-projection is given by
\begin{equation}
    \boldsymbol{\theta}_{k+1} = \mathsf{F}\! \left\{ \mathsf{f}^{\tau}\!
    \left[
        \smash{\underbrace{\mathsf{F}^{-1}\! \left(\boldsymbol{\theta}_k\right)}_{\boldsymbol{u}_{\boldsymbol{\theta}}}}
    \right] \right\} \!,
    \vphantom{\underbrace{1_1^{1}}_{1_1}}
\end{equation}
where the time index $\tau$ is chosen such that $\tau \Delta t \ll T$. By keeping the rollout short, pseudo-projection acts as a dynamical filter that damps non-physical, high-wavenumber components while leaving the long-time trajectory essentially unchanged. Hence, pseudo-projection regularizes an underdetermined inverse problem by incorporating the prior information that admissible states should lie on the attractor, thereby reducing the multiplicity of feasible solutions. This idea is conceptually related to the ``Bayesian--variational cyclic'' method of Gejadze et al. \cite{Gejadze2023}, in which periodic Bayesian updates of lower-dimensional latent variables are used to help combat non-uniqueness in variational problems. It also shares the motivation of preconditioning methods \cite{Haben2011, Ke2026}.\par

We use pseudo-projection in conjunction with NCN. The optimization is performed for 350 iterations, with pseudo-projection applied at steps 50, 100, and 150. We set $\kappa_{|\boldsymbol{H}|} = 10^3$ after the final application of pseudo-projection, allowing the optimizer to refine physical modes.\par

%%% Cases and Sample Results %%%
\section{Test cases and sample reconstructions}
\label{sec: cases}
This section presents representative test cases to illustrate typical outcomes of variational state estimation. We begin by describing the procedure used to generate the dataset of trajectories and initial guesses employed throughout our study, followed by the error metrics used to assess reconstruction quality. We then report representative reconstructions for the $L = 22$ domain to contextualize these metrics, compare optimizer performance, and demonstrate the effects of pseudo-projection, in that order.\par

% Case generation
\subsection{Generation of cases}
\label{sec: cases: generation}
Variational state estimation for KS systems is a highly non-convex problem that depends strongly on both the reference trajectory and the initial guess for the observer trajectory. To marginalize these dependencies and obtain representative reconstruction statistics, we perform reconstructions across a large ensemble of ground truth trajectories and guesses. For each domain size, $L \in \{22, 44, 66\}$, we generate a collection of states on the attractor by integrating a single system forward for 10~000 time units. The first 1000 time units are discarded to ensure convergence to $\mathcal{A}$, and the remaining 9000 time units are retained at intervals of $\Delta t = 1$. The center and radius of the attractor are approximated as
\begin{equation}
    \boldsymbol{u}_{\mathcal{A}}
    \approx \frac{1}{9000} \sum_{k=1001}^{10~000} \boldsymbol{u}_k
    \quad\text{and}\quad
    R_{\mathcal{A}} 
    \approx \frac{1}{9000} \sum_{k=1001}^{10~000}
    \left\| \boldsymbol{u}_k - \boldsymbol{u}_{\mathcal{A}} \right\|_2,
\end{equation}
where $k$ indicates time units. The radius serves as a characteristic scale in state space, and we use it to normalize errors and sample initial guesses at prescribed distances from the true initial condition.\par

For each domain size, we define test cases using 20 random initial conditions and 400 random initial guesses per initial condition, yielding a total of 8000 cases per domain. All cases are reconstructed using data from multiple sensor configurations, with $m_x \in \{2, \dots, 16\}$ spatial sensors and $m_t \in \{2, \dots, 8\}$ measurement times. Both the reference initial conditions and initial guesses (observer systems) are drawn from the long-time rollout for the corresponding domain. The $i$th reference initial state is denoted $\sdx{\boldsymbol{u}_0}[(i)]$, and the $j$th initial guess for that system is $\sdx{\boldsymbol{u}_{\boldsymbol{\theta}, 0}}[(i,j)]$. When generating guesses for a given condition, we compute the distances
\begin{equation}
    D_{ij} = \|\sdx{\boldsymbol{u}_{\boldsymbol{\theta}, 0}}[(i,j)] - \sdx{\boldsymbol{u}_0}[(i)]\|_2.
\end{equation}
We sample states with distances $D_{ij} \in [0.01 R_\mathcal{A}, R_\mathcal{A}]$ to ensure a mixture of good guesses and poor ones. To do this, random target distances are drawn from $[0.01 R_\mathcal{A}, R_\mathcal{A}]$ with uniform probability, and the state for which $D_{ij}$ most closely matches the sample is selected. Duplicates are redrawn until we have 400 unique starting points for our observer system. \par

To simplify notation, we henceforth drop the subscript $0$ when referring to the initial condition of the observer system. We also omit the $(i,j)$ superscript when considering a single observer--reference pair, since no ambiguity arises between different observer systems or reference systems. Thus, we write $\sdx{\boldsymbol{u}_{\boldsymbol{\theta},0}}[(i,j)]$ as $\boldsymbol{u}_{\boldsymbol{\theta}}$.
For observer states at later time indices, with $k > 0$, we write $\boldsymbol{u}_{\boldsymbol{\theta}, k}$.\par

% Metrics
\subsection{Error metrics}
\label{sec: cases: metrics}
Reconstruction accuracy is primarily evaluated using two metrics: a normalized Euclidean distance between initial conditions of the observer and reference systems as well as the cosine similarity between the full trajectories. The initial condition error is
\begin{equation}
    e_{\boldsymbol{u}} = R_\mathcal{A}^{-1} \,\|\boldsymbol{u}_0 - \boldsymbol{u}_{\boldsymbol{\theta}}\|_2,
\end{equation}
where $\boldsymbol{u}_0$ denotes the initial condition of the reference system, and $\boldsymbol{u}_{\boldsymbol{\theta}}$ denotes the initial condition of the observer system. During optimization, $\boldsymbol{\theta}$ is updated, and both $\boldsymbol{u}_{\boldsymbol{\theta}}$ and $e_{\boldsymbol{u}}$ evolve accordingly. Because distinct initial conditions can yield nearly indistinguishable trajectories on the attractor, we also quantify accuracy at the trajectory level via the cosine similarity,
\begin{equation}
    \mathrm{CS}_{\boldsymbol{U}} = 
    \frac{\boldsymbol{U}^\top \boldsymbol{U}_{\boldsymbol{\theta}}}
    {\|\boldsymbol{U}\|_2 
    \|\boldsymbol{U}_{\boldsymbol{\theta}}\|_2},
\end{equation}
where $\boldsymbol{U} = (\boldsymbol{u}_0; \dots; \boldsymbol{u}_K)$ and $\boldsymbol{U}_{\boldsymbol{\theta}} = (\boldsymbol{u}_{\boldsymbol{\theta}}; \dots; \boldsymbol{u}_{\boldsymbol{\theta},K})$ are the ground truth and reconstructed trajectories in $\mathbb{R}^{nK}$.\par

We also assess whether an embedding is well conditioned by computing the largest loss below which trajectory estimates are accurate with high probability. Specifically, we define
\begin{equation}
    \label{equ: sup loss metric}
    \varepsilon^* = \sup 
    \left\{ \varepsilon \mid
    \smash{\underbrace{p\! \left( \smash{\overbrace{\mathrm{CS}_{\boldsymbol{U}} \geq \tau}^\text{accurate est.}} \mid
    \smash{\overbrace{\mathscr{J} < \varepsilon}^\text{of low loss}} \right) \geq 
    1 - \delta}_\text{with high probability}} \right\}
    \vphantom{\underbrace{\overbrace{1}^1}_{1_1}},
\end{equation}
where $\tau \approx 1$ indicates an accurate trajectory and $0 < \delta \ll 1$. Thus, any loss below $\varepsilon^*$ almost certainly corresponds to an accurate reconstruction. Equivalently, let $\boldsymbol{U}_{\mathrm{A}}$ and $\boldsymbol{U}_{\mathrm{B}}$ denote trajectories initialized at $\boldsymbol{u}_{\mathrm{A}}$ and $\boldsymbol{u}_{\mathrm{B}}$, respectively, with measurements $\boldsymbol{y}_{\mathrm{A}} = \Phi(\boldsymbol{u}_{\mathrm{A}})$ and $\boldsymbol{y}_{\mathrm{B}} = \Phi(\boldsymbol{u}_{\mathrm{B}})$. Here, $\Phi$, defined in Eq.~\eqref{eq: phi def}, is the composition of the system flow map with the measurement operator, mapping initial states to a vector containing all available measurements. Well-conditioned embeddings require that trajectories which are close in measurement space are also necessarily close in state space. In particular, accurate estimation occurs whenever
\begin{equation*}
    \frac{1}{m}\,
    \|\boldsymbol{y}_{\mathrm{A}} - \boldsymbol{y}_{\mathrm{B}}\|_2
    < \varepsilon^*.
\end{equation*}
For state estimation to be well posed, $\varepsilon^*$ should exist and it should be reasonably large.\par

% Classification
\subsection{Representative reconstructions}
\label{sec: cases: samples}
\begin{figure}[htb!]
    \center\includegraphics[width=0.95\textwidth]{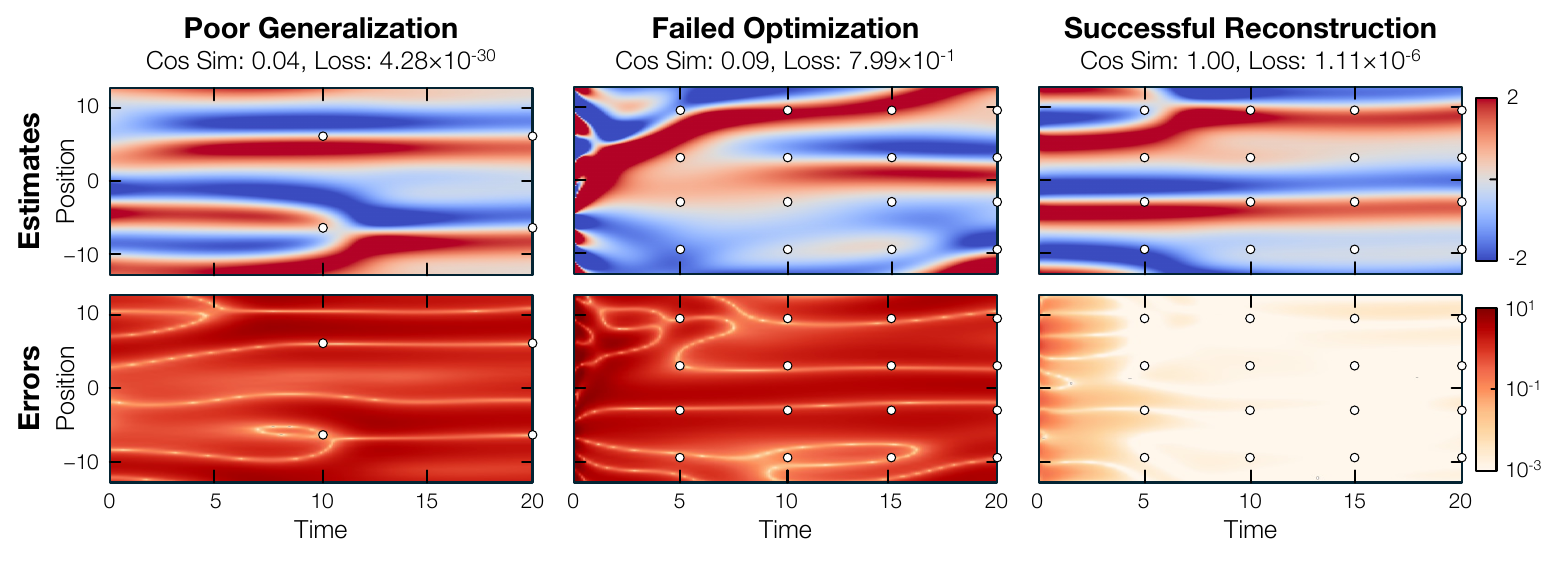}
    \caption{Sample reconstructions for $L = 22$. The top row shows reconstructed trajectories with sensor locations superimposed; the bottom row shows absolute error fields. Final loss values and cosine similarities are reported above each column. From left to right: poor generalization (low loss, high error), failed optimization (high loss, high error), and successful reconstruction (low loss, low error).}
    \label{fig: sample reconstructions}
\end{figure}

State estimation has three characteristic outcomes: poor generalizations, failed optimizations, and successful reconstructions. Figure~\ref{fig: sample reconstructions} presents examples of all three for the $L = 22$ domain. The top row shows reconstructed trajectories with sensor positions superimposed on the estimates and the bottom row shows absolute error fields. All reconstructions were computed using the same reference state, $\boldsymbol{u}_0$, and the same optimizer configuration, namely, our default sequence of 350 NCN iterations with pseudo-projection applied at iterations 50, 100, and 150. The only differences among the examples in Fig.~\ref{fig: sample reconstructions} are the number of observations and the initial guess. For the cases shown (left to right), the initial distances $D_{ij}$ are 0.58, 0.80, and 0.89. The values of $\mathscr{J}$ and $\mathrm{CS}_{\boldsymbol{U}}$ reported in this subsection correspond to these reconstructions.\par

Poor generalization (left) occurs when the observations are sparse ($m_x = 2$ and $m_t = 2$) and the initial guess lies far from the truth. Although the optimizer drives the loss to an extremely low value ($4.3 \times 10^{-30}$), the reconstructed trajectory differs markedly from the reference, yielding a cosine similarity of 0.04. The problem is underdetermined: many trajectories can reproduce this sparse set of observations, so a low loss does not imply an accurate reconstruction. Indeed, as shown in Sec.~\ref{sec: DoF: critical points: local}, even for very good guesses, as $D_{ij} \to 0$, one generally requires $m \geq d_\mathcal{M}$ measurements for the initial state to be observable, where $d_\mathcal{M} = 8$ for the $L = 22$ domain. In this case, the error field exhibits two shallow valleys of extremely low loss centered on the sensor positions, which is a common feature of low-sensor-count reconstructions.\par

Even with a denser set of observations, the optimization can still fail, as shown in the middle panel. Here, the measurement density ($m = 16$) is close to the embedding criterion $m \geq 2d_\mathcal{M} + 1$ discussed in Sec.~\ref{sec: DoF: critical points: global}, yet the optimizer converges to a spurious solution with a high loss ($7.9 \times 10^{-1}$) and a low cosine similarity (0.09). Once again, the error field contains valleys of low loss, though not nearly as deep as in the first reconstruction, and the reconstructed field bears little resemblance to the true system. Such failures arise when the optimizer becomes trapped in high plateaus on the loss landscape, stalling convergence.\par

Lastly, the right panel shows a successful reconstruction for the same sensor configuration as the middle panel, where the optimizer converges to a physically consistent solution with low loss ($1.1 \times 10^{-6}$) and a cosine similarity near unity. Notably, the initial guess in this case was \emph{further} from $\boldsymbol{u}_0$ than in the second example ($D_{ij} = 0.89$ as compared to 0.80). The relationship between the initial separation and the conditions required for accurate reconstruction is discussed in detail in Secs.~\ref{sec: DoF} and \ref{sec: optimization}.\par

% Optimizers
\subsection{Characteristic optimizer behavior}
\label{sec: cases: optimizers}
\begin{figure}[htb!]
    \center\includegraphics[height=6.5cm]{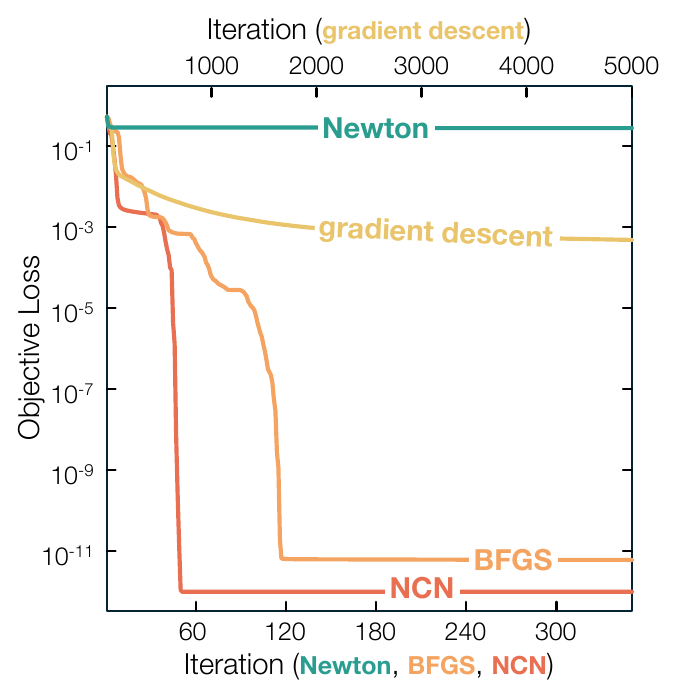} 
    \caption{Optimization loss versus iteration for gradient descent, modified Newton, BFGS, and NCN applied to the same $L = 22$ case with $m_x = 4$ and $m_t = 4$. The lower axis (0--350 iterations) corresponds to Newton, BFGS, and NCN; the upper axis (0--5000 iterations) corresponds to gradient descent.}
    \label{fig: opt conv comp}
\end{figure}

To illustrate the behavior of different optimizers, Fig.~\ref{fig: opt conv comp} shows representative loss traces for a case with $m_x = 4$, $m_t = 4$, and $L = 22$. All methods are initialized from the same initial guess. We compare gradient descent, BFGS, and NCN, as described in Sec.~\ref{sec: state estimation: optimization}, but we do not include a bona fide Newton method. Due to dissipative dynamics and measurement sparsity, the true Hessian for these cases is either very ill-conditioned or degenerate, per Sec.~\ref{sec: optimization}, so exact Newton steps are not well defined. Even with minimal regularization---i.e., retaining only non-zero eigenvalues in Eq.~\eqref{equ: NCN Hessian} (setting $\lambda_i^{-1} = 0$ when $\lambda_i=0$) without enforcing positivity or applying a threshold---Newton iterations quickly diverge. To isolate the effects of negative curvature in our comparison, therefore, we include a modified Newton scheme that applies the same cutoff as NCN but does not enforce positivity of $\boldsymbol{\Lambda}$. Pseudo-projections are omitted from these tests to highlight the intrinsic behavior of each optimizer.\par

Due to severe ill-conditioning of the loss landscape, gradient descent converges at a glacial pace: even after 5000 iterations, its loss remains orders of magnitude higher than those achieved by BFGS and NCN in far fewer steps (note the separate $x$-axes). The modified Newton method also plateaus, despite having access to exact curvature information. Because the Hessian can become indefinite, Newton steps point uphill whenever the gradient overlaps with directions of negative curvature. Such steps would increase the loss, so the backtracking line search sets the step size to zero, causing the optimizer to get stuck. By contrast, BFGS preconditions the gradient with a positive definite matrix $\boldsymbol{B}_k^{-1}$, ensuring that the loss necessarily decreases provided the step is sufficiently small. When the curvature condition fails, we reset $\boldsymbol{B}_k$ to $\boldsymbol{I}$, which is positive definite and thus allows descent to continue. Finally, NCN exhibits the fastest convergence and the lowest final loss, in line with its relative performance across all the test cases we examined.\par

It should be noted that classical and non-convex Newton methods are generally impractical for high-dimensional problems due to the cost of forming, factorizing, and inverting the Hessian. Moreover, the standard Newton step is only guaranteed to be a descent direction when the Hessian is positive definite; otherwise, regularization or modification is required. We present these methods here as diagnostic tools to probe the topology of the loss landscape and explore its implications for state estimation. For higher-dimensional problems, truncated Newton methods \cite{Dembo1983} may provide a suitable alternative. These techniques compute search directions by approximately solving the Newton system using Hessian--vector products with a Krylov subspace method, most commonly via conjugate gradients. The inner iteration is terminated once a prescribed tolerance is met or when negative curvature is detected. In the latter case, if a Krylov direction $\boldsymbol{p}$ satisfies $\boldsymbol{p}^{\top} \boldsymbol{Hp} < 0$, the method terminates the solve and constructs a model-decreasing step from the current Krylov iterate or the detected negative-curvature direction. Thus, truncated Newton methods can handle indefiniteness without explicitly decomposing $\boldsymbol{H}$, and they recover Newton-like behavior when the Hessian is locally positive definite.\par

% Pseudo-projection
\subsection{Reconstructions with pseudo-projection}
\label{sec: cases: projection}
\begin{figure}[htb!]
    \centering
    \includegraphics[width=0.9\textwidth]{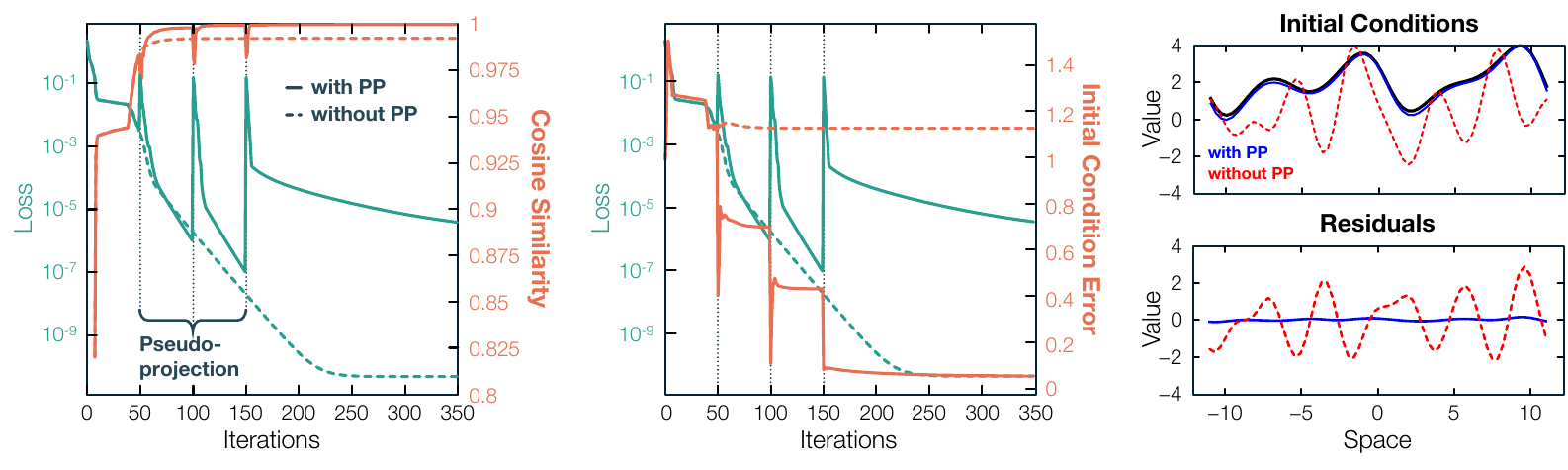} 
    \caption{Effect of pseudo-projections for a case with $m_x = 4$, $m_t = 4$, and $L = 22$. Shown are the loss and cosine similarity (left), the loss and initial-condition error (middle), true and estimated initial conditions (top right), and initial condition residuals (bottom right). Vertical dashed lines indicate pseudo-projection steps, which reduce the initial-state error by nudging the estimate back toward the attractor.}
    \label{fig: PP}
\end{figure}

Recall that pseudo-projection involves a short forward integration of the system dynamics that is meant to bring $\boldsymbol{u}_{\boldsymbol{\theta}}$ closer to $\mathcal{M}$. Figure~\ref{fig: PP} illustrates its effect for a representative case with $m_x = 4$, $m_t = 4$, and $L = 22$. As throughout this paper, we use 350 NCN iterations with pseudo-projection applied at iterations 50, 100, and 150. The plots compare two otherwise identical DA runs: one with pseudo-projection (solid lines) and the other without (dashed lines). The left panel shows traces of the loss $\mathscr{J}$ and cosine similarity $\mathrm{CS}_{\boldsymbol{U}}$, the middle panel shows the loss and the initial condition error $e_{\boldsymbol{u}}$, and the right panel displays initial conditions of the reconstructed and reference systems (top right) as well as the residuals (bottom right).\par

Pseudo-projection events are indicated by vertical dotted lines. At each instance, the loss spikes up and the cosine similarity dips down. Both effects are expected because the action of the system dynamics per se does not account for observations of the reference system. By contrast, the initial condition error consistently decreases with pseudo-projection, meaning that the dynamics do indeed pull $\boldsymbol{u}_{\boldsymbol{\theta}}$ toward $\mathcal{M}$. Notably, the loss always remains below its initial value after projection, and the cosine similarity stays relatively high. Hence, pseudo-projection moves $\boldsymbol{u}_{\boldsymbol{\theta}}$ closer to the IM without fundamentally degrading the estimated trajectory.\par

Between pseudo-projections, the optimizer makes limited progress in reducing $e_{\boldsymbol{u}}$, mainly due to the small number of measurements and the moderately high NCN threshold used in this work. The case without pseudo-projection clearly highlights this limitation: although the optimization achieves a low loss, the final $e_{\boldsymbol{u}}$ is worse than at initialization, and the cosine similarity is lower than in the pseudo-projection case. Thus, for this example, pseudo-projection yields a more accurate trajectory even though the final measurement match is slightly worse. More broadly, pseudo-projection almost always makes the optimization harder for a few steps, i.e., because the spike in loss must be brought back down, but it also introduces perturbations that help the optimizer to escape plateaus or shallow valleys in $\mathscr{J}$. Across all of our tests, we find that pseudo-projection is the primary mechanism for reducing $e_{\boldsymbol{u}}$.\par

\begin{figure}[htb!]
    \centering
    \includegraphics[width=0.9\textwidth]{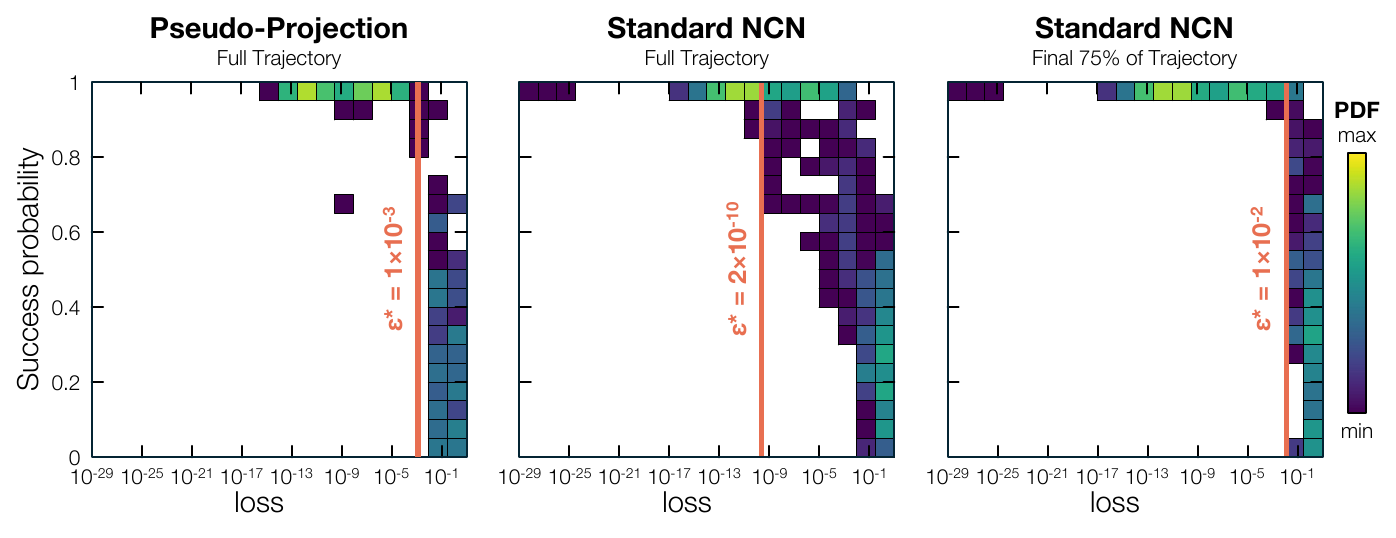}
    \caption{State estimation with and without pseudo-projections for $L = 22$. Each panel shows $p(\mathrm{CS}_{\boldsymbol{U}} = \tau, \mathscr{J})$ in the $(\mathscr{J}, \tau)$ plane. From left to right: with pseudo-projection, without pseudo-projection, and without pseudo-projection but computing $\mathrm{CS}_{\boldsymbol{U}}$ using only the latter 75\% of the trajectory. Vertical lines indicate $\varepsilon^*$.}
    \label{fig: global PP}
\end{figure}

Next, we examine the global effect of pseudo-projection using all 8000 trials of the $m_x = 4,\ = m_t = 4,\ L = 22$ case. Figure~\ref{fig: global PP} shows joint probability density functions (PDFs) $p(\mathrm{CS}_{\boldsymbol{U}} \geq \tau, \mathscr{J})$, evaluated for $\tau = 0.95$. The left and middle panels compare optimizations performed with pseudo-projection (left) and without it (middle). The right panel shows the same analysis, but we restrict the cosine similarity metric to the latter 75\% of the trajectory, thereby excluding early-time transients that take place before the first measurement time in $\mathcal{T}$. For each plot, we also indicate the supremum threshold $\varepsilon^*$ from Eq.~\eqref{equ: sup loss metric}, computed using $\delta = 0.001$.\par

Large values of $\varepsilon^*$ indicate that, for a given optimization scheme, low-loss solutions correspond to accurate reconstructions with high probability. Naturally, small values of $\varepsilon^*$ suggest the opposite, where low-loss solutions can arise from distinct trajectories that nearly match the reference observations. Such cases can occur when trajectories begin off the IM but closely shadow trajectories on it. This behavior is apparent in the middle panel of Fig.~\ref{fig: global PP}, wherein many runs produce a low loss yet have a low probability of accurate reconstruction. With pseudo-projection (left), these spurious \emph{low-loss--low-accuracy} cases are greatly reduced, and $\varepsilon^*$ is much higher. We interpret this increase as improved numerical robustness due to pseudo-projection, since the underlying problem is unchanged across these cases. The right panel confirms that the difficulty originates in early-time reconstruction errors. When the cosine similarity is only computed for the latter 75\% of the trajectory, the bulk of the problematic low-loss--low-accuracy region vanishes. This suggests that pseudo-projection primarily improves observability of the initial condition, as opposed to the full trajectory, by pulling it closer to the IM.\par

%%% Degree-of-Freedom Effects %%%
\section{Degree-of-freedom effects}
\label{sec: DoF}
In this section, we assume the existence of a smooth compact manifold in state space that contains the global attractor with minimal possible dimension and on which the flow map is a diffeomorphism. In other words, we assume an \emph{inertial manifold}, which is known to exist for KS systems. It is also reasonable to posit its existence for more complex dissipative systems such as Navier--Stokes flows, where an IM has not yet been rigorously proven but is expected to exist \cite{Hopf1948}. The results of this section therefore have potential applicability to such systems. A detailed discussion of the existence and properties of IMs in dissipative systems is provided by Zelik~\cite{Zelik2014}.\par

Kuramoto--Sivashinsky dynamics on the IM can be expressed as a system of $d_\mathcal{M}$ ordinary differential equations. Trajectories are determined by state vectors in $\mathbb{R}^{d_\mathcal{M}}$ that specify initial positions on $\mathcal{M}$. Hence, $d_\mathcal{M}$ provides a natural measure of the information necessary to define the system state. Given sufficient knowledge of the system dynamics, $d_\mathcal{M}$ should correspond to the number of measurements $m$ needed for state estimation. Embedding theory formalizes this connection by relating $d_\mathcal{M}$ to the number of measurements $m$ required for a smooth, invertible mapping $\Phi : \mathcal{M} \to \mathbb{R}^m$ to exist, which holds when $m \geq 2d_\mathcal{M} + 1$. While embedding theory has been widely applied to state space reconstruction, we employ it here for the first time to analyze the well-posedness of variational state estimation. Because $d_\mathcal{M}$ increases with the domain length, the number of measurements required for reconstruction likewise grows. We refer to this dependence as a \emph{degree-of-freedom effect on observability}.\par

\begin{figure}[htb!]
    \centering
    \includegraphics[width=0.9\textwidth]{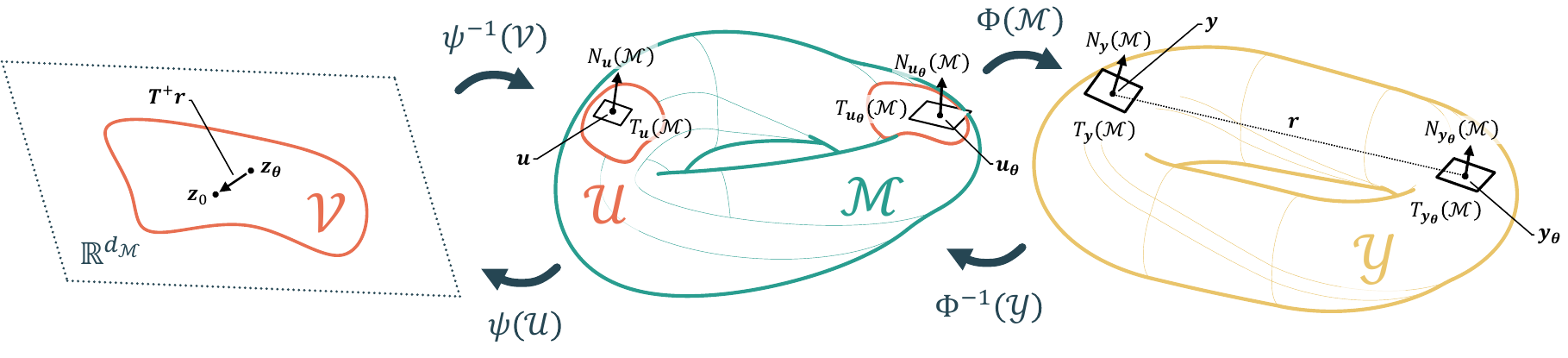}
    \caption{Schematic illustrating the relationships among the inertial manifold $\mathcal{M}$, the measurement manifold $\mathcal{Y}$, and a local Euclidean parametrization $\mathcal{V}$. The observation map $\Phi$ takes states on $\mathcal{M}$ to measurements on $\mathcal{Y}$, and the chart $\psi$ provides local coordinates on $\mathcal{U} \subset \mathcal{M}$. Tangent and normal spaces are shown for representative states and observations.}
    \label{fig: DoF summary}
\end{figure}

Figure~\ref{fig: DoF summary} provides a graphical summary of the spaces relevant to state estimation and their relationships to one another. At the center lies the inertial manifold $\mathcal{M}$, which contains the system's long-time dynamics and is assumed to include both the reference trajectory starting at $\boldsymbol{u}_0$ and the observer trajectory starting at $\boldsymbol{u}_{\boldsymbol{\theta}}$. To the right appears the measurement manifold $\mathcal{Y} = \Phi(\mathcal{M})$, where the observation operator $\Phi$ maps states $\boldsymbol{u} \in \mathcal{M}$ to measurements $\boldsymbol{y} \in \mathcal{Y}$. When $\Phi$ is an embedding, this mapping is smooth and invertible, so $\mathcal{Y}$ and $\mathcal{M}$ are topologically equivalent. On the left is a local embedding: a chart $\psi$ defined on an open patch $\mathcal{U} \subset \mathcal{M}$ that maps $\boldsymbol{u} \in \mathcal{U}$ to manifold coordinates $\boldsymbol{z} \in \mathcal{V} = \psi(\mathcal{U}) \subset \mathbb{R}^{d_\mathcal{M}}$. By definition, such patches form an open cover of $\mathcal{M}$. The remainder of this section develops these spaces and mappings and substantiates their role in variational state estimation.\par

% States to data
\subsection{Mappings from states to measurements}
\label{sec: DoF: mappings}
We begin by introducing the notation used in this section and by briefly reviewing some relevant aspects of embedding theory. Consider the mapping
\begin{equation}
    \Phi : \mathbb{R}^n \to \mathbb{R}^m,
\end{equation}
which takes an initial condition in state space to a vector of $m$ observations,
\begin{equation}
\label{eq: phi def}
    \Phi(\boldsymbol{u}) =
    [\mathsf{h}_1(\boldsymbol{u});
    \mathsf{h}_2(\boldsymbol{u});
    \dots;
    \mathsf{h}_m(\boldsymbol{u})] =
    \boldsymbol{y}.
\end{equation}
To analyze the properties of $\Phi$, we restrict its domain to states on $\mathcal{M}$, such that $\mathcal{Y} = \Phi(\mathcal{M})$ defines the corresponding shadow manifold. We denote by $\mathcal{P}^m$ the space of all such mappings with output dimension $m$. When the domain is further restricted to a subset $\mathcal{B} \subset \mathcal{M}$, we write
\begin{equation}
    \Phi_\mathcal{B}: \mathcal{B} \to \mathcal{Y}.
\end{equation}
While the forward problem $\boldsymbol{y} = \Phi(\boldsymbol{u})$ is well posed, the inverse problem ${\boldsymbol{u}} = \Phi^{-1}(\boldsymbol{y})$ may not be, since $\Phi^{-1}$ may not exist and is typically unavailable in closed form, regardless. Variational state estimation implicitly approximates this inverse through a constrained optimization. The problem can be well posed when each $\boldsymbol{y}$ corresponds to a unique initial condition, which holds for all $\boldsymbol{u} \in \mathcal{M}$ when $\Phi$ is an embedding, ensuring that $\Phi^{-1}$ does indeed exist.\par

Numerous results in the literature on state space reconstruction establish bounds on the number of measurements required for an embedding to exist. Takens' pioneering work showed that if $\boldsymbol{y}$ is a scalar time series obtained from $\boldsymbol{u} \in \mathcal{M}$,\footnote{A vector of delay coordinates contains measurements of $u$ at a fixed spatial position $x$ and at times $t, t + \tau, t + 2\tau, \dots$} then the delay-coordinate map $\Phi : \boldsymbol{u} \mapsto \boldsymbol{y}$ is a diffeomorphic embedding when the number of delays satisfies $m \geq 2d_\mathcal{M} + 1$ \cite{Takens1981, Noakes1991}. Sauer, Yorke, and Casdagli \cite{Sauer1991} later extended this result to strange attractors, showing that a generic observation function yields an embedding when $m > 2d_\mathcal{A}$. Deyle and Sugihara \cite{Deyle2011} further generalized these results to multivariate time series, demonstrating that if $m \geq 2d_\mathcal{M} + 1$, then $\Phi \in \mathcal{P}^m$ is generically an embedding for sufficiently smooth measurement functions, under mild assumptions about periodic points. They also showed that in a probabilistic formulation, any $\Phi_\mathcal{B} \in \mathcal{P}^m$ is almost surely an embedding when $m > 2d_\mathcal{B}$, where $d_\mathcal{B} < d_\mathcal{M}$ is the box-counting dimension of a compact subset $\mathcal{B} \subset \mathcal{M}$. Finally, if $\Phi$ is an embedding, then $\mathcal{M}$ and $\mathcal{Y}$ are topologically equivalent, so that $d_\mathcal{Y} = d_\mathcal{M}$, which suggests a lower bound $m \geq d_\mathcal{M}$ on the number of measurements required for state estimation. The implications of $m \geq d_\mathcal{M}$ and $m \geq 2d_\mathcal{M}+1$ for variational state estimation are derived in Sec.~\ref{sec: DoF: critical points}. Lastly, we note that these theorems establish when $\Phi$ is \emph{almost always} an immersion or an embedding, but they do not provide universal guarantees.\par

% State space reconstruction
\subsection{State space reconstruction}
\label{sec: DoF: SSR}
In state space reconstruction, the goal is to determine invariants of a dynamical system from sparse measurements and to predict their evolution \cite{Kugiumtzis1994}. If $\Phi$ is an embedding, then one can define dynamics of $\boldsymbol{y}$ on $\mathcal{Y}$ that are equivalent to the dynamics of $\boldsymbol{u}$ on $\mathcal{M}$ such that both systems share the same invariants. The measurement space dynamics are written as
\begin{equation}
    \label{equ: measurement dynamics}
    \boldsymbol{y}_{k+1} = \mathsf{g}_{\Delta t}(\boldsymbol{y}_k),
\end{equation}
where $\mathsf{g}_{\Delta t}$ is simply
\begin{equation}
    \mathsf{g}_{\Delta t} = \Phi \circ \mathsf{f} \circ \Phi^{-1}.
\end{equation}
Hence, the system dynamics can be examined entirely in measurement space when $\Phi$ is an embedding. The existence of $\Phi^{-1}$ is sufficient for this purpose, which stands in contrast to state estimation, where we seek a functional approximation to $\Phi^{-1}$.\par

% Sensor setup
\subsubsection{Sensor placement and repetition rate}
\label{sec: DoF: SSR: setup}
Although embedding theorems specify how many measurements are needed to establish a diffeomorphic mapping from $\mathcal{M}$ to $\mathcal{Y}$, they offer no guidance on where to place sensors or how rapidly they should record observations of $\boldsymbol{u}$. These are critical considerations in practice, especially for noisy measurements \cite{Kugiumtzis1996}. A standard approach for selecting the time lag $\tau$ between measurements is to analyze the average mutual information between measurements at times $t$ and $t + \tau$. The lag is often chosen as the first minimum of mutual information with increasing $\tau$ \cite{Fraser1986}, minimizing redundancy while ensuring that successive measurements are still correlated \cite{Abarbanel1993, Rosenstein1994}. Alternatively, one can fix the measurement time horizon $T$, from which $\tau$ is determined by the number of measurements as $\tau = T/m$. Rosenstein et al. \cite{Rosenstein1994} showed that the optimal lag scales with $m$ such that $T$ remains approximately constant. By fixing $T$, one can keep the earliest and latest measurements within a window where their dynamical relationship remains computable. A myriad of methods have been proposed to optimize $T$ or $\tau$ \cite{Caputo1986, Gibson1992, Rosenstein1994}, all aiming to strike a reasonable compromise between \emph{redundancy} (short $\tau$, strongly correlated measurements) and \emph{irrelevance} (long $\tau$, decorrelated measurements that convey little information about the initial state) \cite{Casdagli1991}.\par

For variational state estimation, prior studies recommended restricting the assimilation window to the Lyapunov time $T_\ell$ \cite{Li2020, Chandramouli2020}. Beyond this scale, the exponential sensitivity to initial conditions is assumed to cause gradient calculations to rapidly deteriorate. When $K \Delta t \gg T_\ell$, the \emph{computed} probability density $p(\boldsymbol{u}_K \mid \boldsymbol{u}_0)$ approaches the unconditional distribution $p(\boldsymbol{u}_0)$, and as a corollary, we have $p(\boldsymbol{u}_0 \mid \boldsymbol{u}_K) \to p(\boldsymbol{u}_0)$. Therefore, setting $T = T_\ell$ is a pragmatic choice for the assimilation window, independent of $m$, and we follow this convention throughout the present work. Because the literature on state space reconstruction primarily concerns 1D time series, however, there is little precedent for spatial sensor placement. We thus adopt uniform spatial coverage under the assumption that all spatial locations are equally informative.\par

% Critical points
\subsection{Critical points on \texorpdfstring{$\mathscr{J}$}{the loss functional}}
\label{sec: DoF: critical points}
In Sec.~\ref{sec: state estimation: optimization}, we discuss evidence that critical points of high loss are rare in high-dimensional non-convex optimization problems, while critical points of low loss are typically saddle points or global minima (possibly with multiple minima of equal loss). Here, we examine the conditions under which critical points can arise in the loss landscape on $\mathcal{M}$. We show that, under suitable assumptions, the global minimum is the only critical point on manifold, which holds locally for $m \geq d_\mathcal{M}$ and globally for $m \geq 2d_\mathcal{M} + 1$.\par

% Definitions
\subsubsection{Some definitions}
\label{sec: DoF: critical points: definitions}
Several geometric quantities must be defined to assess critical points on the IM. The mapping $\Phi$ must be an \emph{immersion} to qualify as an embedding, and an \emph{atlas of charts} is needed to parameterize the loss landscape. An immersion is simply a local embedding: for every state $\boldsymbol{u} \in \mathcal{M}$, there exists a neighborhood around it such that $\boldsymbol{u} \in \mathcal{U} \subset \mathcal{M}$, wherein the restricted mapping $\Phi_\mathcal{U} : \mathcal{U} \to \mathcal{Y}$ is an embedding \cite{bishop2012}. Immersions are known to exist generically when the number of measurements satisfies $m \geq d_\mathcal{M}$ \cite{Deyle2011}.\par

Equivalently, immersions can be characterized using the \emph{tangent spaces} of $\mathcal{M}$ and $\mathcal{Y}$, which are depicted in Fig.~\ref{fig: DoF summary}. The tangent space of a smooth manifold $\mathcal{B}$ at $\boldsymbol{x}$ is denoted $T_{\boldsymbol{x}}(\mathcal{B})$, with a Euclidean dimension that necessarily equals the manifold dimension $d_\mathcal{B}$, and the normal space is $N_{\boldsymbol{x}}(\mathcal{B})$. For states $\boldsymbol{u}$ and measurements $\boldsymbol{y}$, the Jacobian of $\Phi$ with respect to $\boldsymbol{u}$ must relate both the tangent and normal spaces of $\mathcal{M}$ and $\mathcal{Y}$,
\begin{equation}
    \label{equ: immersion}
    \frac{\partial \Phi}{\partial \boldsymbol{u}} :
    \underbrace{T_{\boldsymbol{u}}(\mathcal{M}) \oplus N_{\boldsymbol{u}}(\mathcal{M})}_{\mathbb{R}^n} \to
    \underbrace{T_{\boldsymbol{y}}(\mathcal{Y}) \oplus N_{\boldsymbol{y}}(\mathcal{Y})}_{\mathbb{R}^m}.
\end{equation}
The immersion property pertains solely to the restricted mapping
\begin{equation}
    \label{equ: restricted immersion}
    \left(\frac{\partial \Phi}{\partial \boldsymbol{u}}\right)_{T_{\boldsymbol{u}}(\mathcal{M})} :
    T_{\boldsymbol{u}}(\mathcal{M}) \to T_{\boldsymbol{y}}(\mathcal{Y}).
\end{equation}
If $\Phi$ is an immersion, this mapping is \emph{injective}, which implies that $d_\mathcal{M} \leq d_\mathcal{Y}$. Since the measurement manifold is the image of $\mathcal{M}$, i.e., $\mathcal{Y} = \Phi(\mathcal{M})$, we also have $d_\mathcal{Y} \leq d_\mathcal{M}$. Hence, $d_\mathcal{M}$ must equal $d_\mathcal{Y}$ and the tangent map is a bijection. A first-order Taylor series expansion of $\Phi$ gives
\begin{equation}
    \label{equ: immersion Taylor}
    \delta \boldsymbol{y} \approx
    \frac{\partial \Phi}{\partial \boldsymbol{u}}\,\delta \boldsymbol{u},
\end{equation}
so if $\Phi$ is an immersion, then any non-zero perturbation $\delta \boldsymbol{u} \in T_{\boldsymbol{u}}(\mathcal{M})$ produces a non-zero measurement perturbation $\delta \boldsymbol{y} \in T_{\boldsymbol{y}}(\mathcal{Y})$ and vice versa. In other words, as an immersion, $\Phi$ resolves all the intrinsic directions at all points on the inertial and shadow manifolds.\par

Next, we define an atlas, which allows us to map from a state $\boldsymbol{u} \in \mathcal{M}$ to a vector of manifold coordinates $\boldsymbol{z} \in \mathbb{R}^{d_\mathcal{M}}$ corresponding to the system's intrinsic degrees of freedom. An atlas is a collection of charts whose domains $\mathcal{U} \subset \mathcal{M}$ form an open cover of $\mathcal{M}$, where $\mathcal{U}$ is one of many such domains. Each chart
\begin{equation}
    \psi : \mathcal{U} \to \mathcal{V} \subset \mathbb{R}^{d_\mathcal{M}}
\end{equation}
is a diffeomporhism from $\mathcal{U}$ onto an open subset of $\mathbb{R}^{d_\mathcal{M}}$, such that $\boldsymbol{u} = \psi^{-1}(\boldsymbol{z})$ for $\boldsymbol{z} \in \mathcal{V}$. The tangent space of a manifold can also be obtained by differentiating the inverse of a chart. Specifically, 
\begin{equation}
    \label{equ: chart to tangent}
    T_{\boldsymbol{u}}(\mathcal{M})
    = \operatorname{span}\! \left(\frac{\partial \psi^{-1}}{\partial \boldsymbol{z}}\right) \!,
\end{equation}
where $\boldsymbol{u} \in \mathcal{U}$, $\boldsymbol{z} \in \mathcal{V}$, and the rank of the Jacobian is $d_\mathcal{M}$. An  example of this mapping is illustrated in Fig.~\ref{fig: DoF summary}.\par

% Local results
\subsubsection{Local behavior}
\label{sec: DoF: critical points: local}
To begin, we show that if $\Phi$ is an immersion, and if the initial conditions of the observer and reference systems, i.e., $\boldsymbol{u}_{\boldsymbol{\theta}}$ and $\boldsymbol{u}_0$, are confined to a sufficiently small region on a chart domain $\mathcal{U}$, then there exists a single critical point in this region at $\boldsymbol{u}_{\boldsymbol{\theta}} = \boldsymbol{u}_0$. The proximity of $\boldsymbol{u}_{\boldsymbol{\theta}}$ and $\boldsymbol{u}_0$ is required to justify a first-order Taylor expansion. In what follows, $\boldsymbol{z}_{\boldsymbol{\theta}} = \psi^{-1}(\boldsymbol{u}_{\boldsymbol{\theta}})$ and $\boldsymbol{z}_0 = \psi^{-1}(\boldsymbol{u}_0)$ are representations of $\boldsymbol{u}_{\boldsymbol{\theta}}$ and $\boldsymbol{u}_0$ in manifold coordinates, with $\boldsymbol{z}_{\boldsymbol{\theta}}, \boldsymbol{z}_0 \in \mathcal{V} \subset \mathbb{R}^{d_\mathcal{M}}$. We thus define the measurements as a function of $\boldsymbol{z}_0$,
\begin{equation}
    \boldsymbol{y} = \Phi \circ \psi^{-1}(\boldsymbol{z}_0),
\end{equation}
and so too for $\boldsymbol{y}_{\boldsymbol{\theta}}$ and $\boldsymbol{z}_{\boldsymbol{\theta}}$. The first-order Taylor series expansion about $\boldsymbol{z}_{\boldsymbol{\theta}}$ gives
\begin{equation}
    \label{equ: Taylor series}
    \boldsymbol{y} = \boldsymbol{y}_{\boldsymbol{\theta}} + \boldsymbol{T} (\boldsymbol{z}_0 - \boldsymbol{z}_{\boldsymbol{\theta}}),
\end{equation}
where
\begin{equation}
    \boldsymbol{T} =
    \frac{\partial \Phi}{\partial \boldsymbol{u}_{\boldsymbol{\theta}}}
    \frac{\partial \psi^{-1}}{\partial \boldsymbol{z}_{\boldsymbol{\theta}}}.
\end{equation}
The row space of ${\partial \Phi}/{\partial \boldsymbol{u}_{\boldsymbol{\theta}}} \in \mathbb{R}^{m \times n}$ contains $T_{\boldsymbol{u}_{\boldsymbol{\theta}}}(\mathcal{M})$ since $\Phi$ is an immersion and thus has a rank greater than or equal to $d_\mathcal{M}$. The column space of ${\partial \psi^{-1}}/{\partial \boldsymbol{z}}_{\boldsymbol{\theta}} \in \mathbb{R}^{n \times d_\mathcal{M}}$, whose rank is exactly $d_\mathcal{M}$, is identically $T_{\boldsymbol{u}_{\boldsymbol{\theta}}}(\mathcal{M})$. Therefore, the rank of $\boldsymbol{T} \in \mathbb{R}^{m \times d_\mathcal{M}}$ is $d_\mathcal{M}$.\par

We see this scenario on the left side of Fig.~\ref{fig: DoF summary}. Two nearby states that fall within the same chart domain $\mathcal{U}$ are mapped into $\mathcal{V} \subset \mathbb{R}^{d_\mathcal{M}}$, yielding manifold coordinates $\boldsymbol{z}_0$ and $\boldsymbol{z}_{\boldsymbol{\theta}}$. From Eq.~\eqref{equ: Taylor series}, the points are separated by
\begin{equation*}
    \underbrace{\left( \boldsymbol{T}^\top \boldsymbol{T} \right)^{-1} \boldsymbol{T}^\top}_{\boldsymbol{T}^+}
    \underbrace{\left( \boldsymbol{y} - \boldsymbol{y}_{\boldsymbol{\theta}} \right)}_{\boldsymbol{r}} = \boldsymbol{z}_0 - \boldsymbol{z}_{\boldsymbol{\theta}},
\end{equation*}
where $\boldsymbol{T}^+$ is the pseudoinverse of $\boldsymbol{T}$ and $\boldsymbol{r}$ is the measurement residual.\par

The loss functional may be written as
\begin{equation}
    \label{equ: z_loss_fn}
    \mathscr{J} = \frac{1}{2} (\boldsymbol{y} - \boldsymbol{y}_{\boldsymbol{\theta}})^\top
    (\boldsymbol{y} - \boldsymbol{y}_{\boldsymbol{\theta}})
    = \frac{1}{2} \boldsymbol{r}^\top \boldsymbol{r},
\end{equation}
which is equivalent to Eq.~\eqref{equ: discrete loss} up to a constant. Differentiating it with respect to ${\boldsymbol{z}}_{\boldsymbol{\theta}}$ gives
\begin{equation}
    \label{equ: z_loss_fn_grad}
    \frac{\partial \mathscr{J}}{\partial {\boldsymbol{z}}_{\boldsymbol{\theta}}}
    = \boldsymbol{r}^\top \boldsymbol{T}.
\end{equation}
Substituting the first-order expansion from Eq.~\eqref{equ: Taylor series} yields
\begin{equation}
   \frac{\partial \mathscr{J}}{\partial \boldsymbol{z}_{\boldsymbol{\theta}}}
    = \left(\boldsymbol{z}_0 - \boldsymbol{z}_{\boldsymbol{\theta}}\right)^\top \boldsymbol{T}^\top \boldsymbol{T}.
\end{equation}
Note that $\mathrm{rank}(\boldsymbol{T}^\top \boldsymbol{T}) = \mathrm{rank}(\boldsymbol{T}) = d_\mathcal{M}$, so the gradient only vanishes when $\boldsymbol{z}_{\boldsymbol{\theta}} = \boldsymbol{z}_0$. Hence, so long as $\Phi$ is an immersion, the only critical point local to the global minimum $\boldsymbol{u}_0$ is in fact $\boldsymbol{u}_0$ itself.\par

This result establishes a lower bound on the number of measurements required for local state estimation. If $m < d_\mathcal{M}$, then $\mathrm{rank}(\boldsymbol{T}^\top \boldsymbol{T}) \leq m < d_\mathcal{M}$ and the quadratic approximation to $\mathscr{J}$ in manifold coordinates is degenerate at optimality. In this case, there exist nonzero tangent perturbations that do not change the measurements to first order and, as a result, the initial state is not locally observable from $\boldsymbol{y}$. Conversely, if $m \geq d_\mathcal{M}$ and $\Phi$ is an immersion, then the restricted tangent map is full rank. Hence, for an initial guess $\boldsymbol{u}_{\boldsymbol{\theta}} \in \mathcal{M}$ that is sufficiently close to $\boldsymbol{u}_0$, the displacement between $\boldsymbol{u}_{\boldsymbol{\theta}}$ and $\boldsymbol{u}_0$ can be represented to first order by a tangent perturbation in $T_{\boldsymbol{u}_0}(\mathcal{M})$, and the corresponding measurement residual lies to first order in $T_{\boldsymbol{y}}(\mathcal{Y})$. Consequently, the local quadratic loss function has a unique minimizer at $\boldsymbol{u}_0$ and gradient-based optimization with appropriate step sizes converges to this minimizer. Therefore, $m \geq d_\mathcal{M}$ gives the minimum number of measurements needed for $\boldsymbol{u}_0$ to be locally observable from an arbitrarily good initial guess. This may be a practical limit for sequential smoothers or filters, where the initial guess for each segment can become accurate after several assimilation windows or analysis steps.\par

% Global results
\subsubsection{Global behavior}
\label{sec: DoF: critical points: global}
Next, we look into the properties of critical points when $\boldsymbol{u}_{\boldsymbol{\theta}}$ and $\boldsymbol{u}_0$ need not be close. In particular, we show that $\boldsymbol{u}_0$ is the only critical point on $\mathcal{M}$ when $\Phi$ is an embedding. Starting from the gradient of Eq.~\eqref{equ: z_loss_fn} with respect to $\boldsymbol{u}_{\boldsymbol{\theta}}$,
\begin{equation}
    \frac{\partial \mathscr{J}}{\partial \boldsymbol{u}_{\boldsymbol{\theta}}}
    = \boldsymbol{r}^\top
      \frac{\partial \Phi}{\partial \boldsymbol{u}_{\boldsymbol{\theta}}},
\end{equation}
critical points arise either when $\boldsymbol{y} = \boldsymbol{y}_{\boldsymbol{\theta}}$ or when the residual lies in the left null space of $\partial \Phi / \partial \boldsymbol{u}_{\boldsymbol{\theta}}$. The column space of this Jacobian generically spans $\mathbb{R}^m$ whenever $m \leq n$, so all critical points satisfy $\boldsymbol{y} = \boldsymbol{y}_{\boldsymbol{\theta}}$.\footnote{It has been shown that one can independently perturb each observation function $h_i$ to obtain $m$ linearly independent tangent vectors $\partial h_i / \partial \boldsymbol{u}_{\boldsymbol{\theta}}$ \cite{Deyle2011}. Additional justification is required for cases with $m > n$.} This conclusion holds even if $\Phi$ is not an immersion or an embedding. Going further, when $\Phi$ is indeed an embedding, then $\boldsymbol{y} = \boldsymbol{y}_{\boldsymbol{\theta}}$ can occur only if $\boldsymbol{u}_{\boldsymbol{\theta}} = \boldsymbol{u}_0$ for $\boldsymbol{u}_{\boldsymbol{\theta}}, \boldsymbol{u}_0 \in \mathcal{M}$ since $\Phi : \mathcal{M} \to \mathcal{Y}$ is a bijection. Therefore, the only critical point on the manifold is the global minimum. To the best of the author's knowledge this is an original result.\par

Recall that $\Phi$ is generically an embedding when $m \geq 2d_\mathcal{M} + 1$ \cite{Deyle2011}. Consequently, when $m$ satisfies this bound, one might expect the state estimation problem to be well posed.\par

Regrettably, we note that the existence of a single critical point at $\boldsymbol{u}_0$ on $\mathcal{M}$ does not guarantee convergence to that point via variational state estimation. Even when the optimization begins on $\mathcal{M}$, the gradient may have components orthogonal to $T_{\boldsymbol{u}_{\boldsymbol{\theta}}}(\mathcal{M})$, pushing $\boldsymbol{u}_{\boldsymbol{\theta}}$ off the manifold, where the above analysis no longer holds and where additional minima may exist. Constraining the optimization to $\mathcal{M}$ by projecting $\boldsymbol{g}_{\boldsymbol{u}}$ onto the local tangent space---where $\boldsymbol{g}_{\boldsymbol{u}} = \partial \mathscr{J}/\partial \boldsymbol{u}_{\boldsymbol{\theta}}$ is the state space gradient---which is loosely approximated by our pseudo-projection procedure, helps to mitigate this issue. However, if $\boldsymbol{g}_{\boldsymbol{u}}$ happens to be in $N_{\boldsymbol{u}_{\boldsymbol{\theta}}} (\mathcal{M})$, then NCN steps counteract pseudo-projection and even a true manifold-constrained (i.e., Riemannian) optimization would stall. Thus, although $\boldsymbol{u}_0$ is the only critical point on $\mathcal{M}$ when $\Phi$ is an embedding, the gradient need not lie within the local tangent space, and specialized optimization strategies may be required to handle such pathologies.\par

To show that gradients can in fact point off the manifold, we invoke Whitney's strong embedding theorem, which states that the measurement manifold $\mathcal{Y} = \Phi(\mathcal{M})$, of intrinsic dimension $d_\mathcal{Y} \leq d_\mathcal{M}$, can be smoothly embedded in a Euclidean space $D \subset \mathbb{R}^m$ of dimension $d$, where for any non-linear manifold we have $d_\mathcal{Y} < d \leq \min(m, 2d_\mathcal{Y})$. Here, $D$ represents the \emph{minimal} Euclidean space that embeds $\mathcal{Y}$. The residual $\boldsymbol{r}$ necessarily lies in $D$, and because $d > d_\mathcal{Y}$ for a non-linear manifold, there must exist residuals with components in the normal space $N_{\boldsymbol{y}_{\boldsymbol{\theta}}}(\mathcal{Y})$. A visual example of this is provided on the right-hand side of Fig.~\ref{fig: DoF summary}, where the residual $\boldsymbol{r}$ does not fully reside in  $T_{\boldsymbol{y}_{\boldsymbol{\theta}}}(\mathcal{Y})$. In such instances, when $\Phi$ is an immersion or an embedding, $\boldsymbol{g}_{\boldsymbol{u}}$ necessarily contains components in $N_{\boldsymbol{u}_{\boldsymbol{\theta}}}(\mathcal{M})$, as argued next.\par

For any immersion $\Phi$, the restricted mapping
\begin{equation*}
    \left(\frac{\partial \Phi} {\partial \boldsymbol{u}_{\boldsymbol{\theta}}}\right)_{T_{\boldsymbol{u}_{\boldsymbol{\theta}}}(\mathcal{M})} :
    T_{\boldsymbol{u}_{\boldsymbol{\theta}}}(\mathcal{M}) \to
    T_{\boldsymbol{y}_{\boldsymbol{\theta}}}(\mathcal{Y})
\end{equation*}
is bijective, even though the full mapping between ambient spaces $\mathbb{R}^n \to \mathbb{R}^m$ need not be. Consequently, the gradient
\begin{equation}
    \boldsymbol{g}_{\boldsymbol{u}} =
    \left(\frac{\partial \Phi}{\partial \boldsymbol{u}_{\boldsymbol{\theta}}}\right)^\top
    \boldsymbol{r},
\end{equation}
must take any component of $\boldsymbol{r}$ that lies in $N_{\boldsymbol{y}_{\boldsymbol{\theta}}}(\mathcal{Y})$ to $N_{\boldsymbol{u}}(\mathcal{M})$, because mapping such a component into $T_{\boldsymbol{u}}(\mathcal{M})$ would contradict the bijectivity of $\partial \Phi / \partial \boldsymbol{u}_{\boldsymbol{\theta}}$ restricted to the tangent spaces. Thus, Whitney's theorem implies that residuals with normal components exist, and therefore some gradients must point off the IM when $\Phi$ is an immersion or an embedding.\par

While the residual $\boldsymbol{r}$ does not generally lie in $N_{\boldsymbol{y}}(\mathcal{Y})$, there exist points on many manifolds for which this occurs. For instance, on a circular manifold, the displacement between any pair of antipodal points is normal to the manifold. We hypothesize that \emph{some} such configurations could act as attractors in manifold-constrained optimization, posing a potential but likely uncommon pathology for variational state estimation.\par

% Tagent spaces
\subsection{Tangent spaces on \texorpdfstring{$\mathcal{M}$ and $\mathcal{Y}$}{the inertial and shadow manifolds}}
\label{sec: DoF: tangents}
Even if $\Phi$ is an embedding, the stability of variational state estimation depends on two additional factors: (1)~the numerical conditioning of the measurement map when restricted to the IM and (2)~the extent to which gradients of the loss remain aligned with the manifold. Because both the reference and observer trajectories lie on $\mathcal{M}$, these effects are governed by the local geometric structure of $\mathcal{M}$ and its image $\mathcal{Y}$. Up next, we empirically investigate the condition number of the Jacobian restricted to the tangent spaces of $\mathcal{M}$ and $\mathcal{Y}$. Since good conditioning alone does not prevent the optimizer from drifting off the manifold, we then quantify how often gradients possess non-trivial components in the normal directions. We note that these results are based on data from KS simulations and are therefore not general; however, we believe they provide useful insights and intuition.\par

% Condition number
\subsubsection{Conditioning of \texorpdfstring{$\partial \Phi/ \partial \boldsymbol{u}_{\boldsymbol{\theta}}$}{the Jacobian} restricted to \texorpdfstring{$T_{\boldsymbol{u}}(\mathcal{M}) \to T_{\boldsymbol{y}}(\mathcal{Y})$}{the tangent spaces}}
\label{sec: DoF: tangents: condition}
To restrict our analysis of $\partial \Phi / \partial \boldsymbol{u}_{\boldsymbol{\theta}}$ to the mapping between tangent spaces, we must construct a projection operator that maps $\mathbb{R}^n \to T_{\boldsymbol{u}}(\mathcal{M})$ for any state $\boldsymbol{u} \in \mathcal{M}$. Theoretically, such an operator can be obtained by differentiating the inverse of a chart with respect to $\boldsymbol{u}$, since the span of this Jacobian equals $T_{\boldsymbol{u}}(\mathcal{M})$, per Eq.~\eqref{equ: chart to tangent}. To this end, we first recall the definition of a chart:
\begin{equation*}
    \boldsymbol{z} = \psi(\boldsymbol{u}) \quad\text{and}\quad
    \boldsymbol{u} = \psi^{-1}(\boldsymbol{z}).
\end{equation*}
Although $\psi$ and its inverse are not available in closed form, we approximate these mappings using the encoder $\mathsf{E}$ and decoder $\mathsf{D}$ introduced in Sec.~\ref{sec: state estimation: KS} and detailed in Appendix~\ref{app: AE}.\par

The autoencoder's latent space $\mathcal{L}$ generally has an oversized dimension, $d_\mathcal{L} \geq d_\mathcal{M}$. In order to obtain a reduced representation that is consistent with the manifold dimension, we perform a PCA on the latent states from our long-time rollout and retain the first $d_\mathcal{M}$ principal components, storing them in $\boldsymbol{P} \in \mathbb{R}^{d_\mathcal{L} \times d_\mathcal{M}}$, as well as the mean latent vector $\boldsymbol{\ell}$. Our choice of the number of principal components is motivated by the sharp spectral drop immediately after $d_\mathcal{M}$, as can be seen in Fig.~\ref{fig: IM dimension}. This qualitative criterion is consistent with previous work using autoencoders for low-order modeling of chaotic systems \cite{Zeng2024}. Together, these elements define the affine transformation used to approximate the mappings $\mathcal{L} \to \mathcal{V}$ and $\mathcal{V} \to \mathcal{L}$. We express the chart and its inverse as
\begin{equation*}
    \boldsymbol{z} \approx \boldsymbol{P}^\top [\mathsf{E}(\boldsymbol{u}) - \boldsymbol{\ell}] \quad\text{and}\quad
    \boldsymbol{u} \approx \mathsf{D}(\boldsymbol{Pz} + \boldsymbol{\ell}).
\end{equation*}
While an atlas of charts is needed to cover $\mathcal{M}$, corresponding to a set of encoders and decoders with one pair per chart, the manifolds considered in this work are well represented by a single pair of global mappings, $\mathsf{E} : \mathcal{M} \to \mathcal{L}$ and $\mathsf{D} : \mathcal{L} \to \mathcal{M}$. The method proposed by Floryan and Graham \cite{Floryan2022} can be employed when multiple local mappings are required.\par

Given a differentiable approximation to $\psi^{-1}$, we sample $\boldsymbol{u} \in \mathcal{M}$ from the rollout and compute
\begin{equation}
    \frac{\partial \mathsf{D}(\boldsymbol{Pz} + \boldsymbol{\ell})}{\partial \boldsymbol{u}} \approx
    \frac{\partial \psi^{-1}}{\partial \boldsymbol{u}}
\end{equation}
via AD. The Jacobian $\partial \Phi / \partial \boldsymbol{u}$ is also obtained by AD. Applying a QR decomposition to our approximation of $\partial \psi^{-1} / \partial \boldsymbol{u}$ yields an orthonormal basis $\boldsymbol{Q} : \mathbb{R}^n \to T_{\boldsymbol{u}}(\mathcal{M})$ whose columns span the tangent space. We may therefore restrict the Jacobian to the mapping between tangent spaces as follows:
\begin{equation}
    \frac{\partial \Phi}{\partial \boldsymbol{u}} \boldsymbol{QQ}^\top \approx
    \left(\frac{\partial \Phi} {\partial \boldsymbol{u}}\right)_{T_{\boldsymbol{u}}(\mathcal{M})} :
    T_{\boldsymbol{u}}(\mathcal{M}) \to
    T_{\boldsymbol{y}}(\mathcal{Y}),
\end{equation}
since the null space of $\boldsymbol{QQ}^\top$ is $N_{\boldsymbol{u}}(\mathcal{M})$. We finally compute the singular value decomposition (SVD) of $(\partial \Phi / \partial \boldsymbol{u})\,\boldsymbol{Q}$ to obtain the spectrum of the restricted mapping.\footnote{The operator $\boldsymbol{Q}\boldsymbol{Q}^\top$ could be used to perform a manifold-constrained optimization. We successfully implemented a related approach with $\boldsymbol{z}$ as the control vector, i.e., by applying AD to the computational graph from $\boldsymbol{z} \to \boldsymbol{y}$. However, since the manifold is not known a priori for most reconstruction problems, we do not employ such techniques in the present work.}\par

\begin{figure}[htb!]
    \centering
    \includegraphics[height=6.5cm]{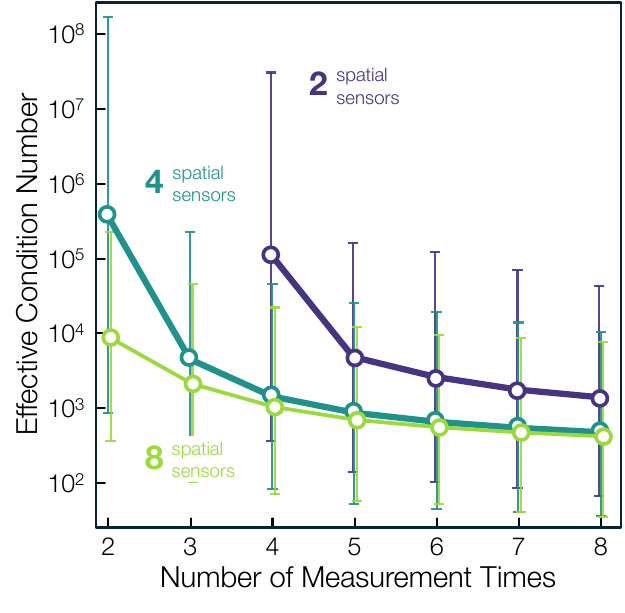}
    \caption{Mean (dots) and range (vertical lines) of condition numbers for mappings $T_{\boldsymbol{u}}(\mathcal{M}) \to T_{\boldsymbol{y}}(\mathcal{Y})$ for the $L = 22$ domain. Curves correspond to different numbers of spatial sensors $m_x$ and are plotted against the number of observation times $m_t$. All maps are full rank, consistent with the immersion criterion, and show improved conditioning with increased spatial and temporal sampling.}
    \label{fig: AE immersion}
\end{figure}

Figure~\ref{fig: AE immersion} summarizes the condition numbers obtained from the tangent space mapping for 1000 snapshots sampled from the long-time rollout in the $L = 22$ domain. Mean condition numbers are shown as solid dots, and vertical lines indicate the corresponding ranges. Results are plotted as a function of the number of measurement times $m_t$, with a separate curve for each number of spatial sensors $m_x$. Although the condition numbers are large, they remain finite and are well below the inverse of machine precision, corroborating that the mapping from $T_{\boldsymbol{u}}(\mathcal{M})$ to $T_{\boldsymbol{y}}(\mathcal{Y})$ is indeed a bijection. This behavior is consistent with the immersion criterion, whereby $\Phi \in \mathcal{P}^m$ is generically an immersion if $m \geq d_\mathcal{M}$. However, such large condition numbers imply poor numerical conditioning, meaning that some tangent directions are only weakly resolved by the measurements. Gradient components along those directions could thus be strongly attenuated, impeding optimization. As expected, we also see that conditioning improves with additional spatial and temporal observations, reflecting a more stable mapping between the IM and measurement space. That being said, improvements in the condition number inherently level off at large $m$ because it is bounded from below by the conditioning of the flow map Jacobian restricted to $T_{\boldsymbol{u}}(\mathcal{M})$, as discussed in Sec.~\ref{sec: optimization: condition}.\par

% Gradients
\subsubsection{Gradient components in \texorpdfstring{$T_{\boldsymbol{u}}(\mathcal{M})$ and $N_{\boldsymbol{u}}(\mathcal{M})$}{tangent and normal spaces}}
\label{sec: DoF: tangents: gradients}
Section~\ref{sec: DoF: critical points: global} shows that reference--observer pairs exist for which $\boldsymbol{g}_{\boldsymbol{u}} \notin T_{\boldsymbol{u}_{\boldsymbol{\theta}}}(\mathcal{M})$. Moreover, if a direction in $N_{\boldsymbol{y}_{\boldsymbol{\theta}}}(\mathcal{Y})$ intersects $\mathcal{Y}$, then gradients with $\boldsymbol{g}_{\boldsymbol{u}} \in N_{\boldsymbol{u}_{\boldsymbol{\theta}}}(\mathcal{M})$ can occur. The existence of such cases would prevent the theoretical global convergence of manifold-constrained optimization, and in practice any component of the gradient lying in $N_{\boldsymbol{u}_{\boldsymbol{\theta}}}(\mathcal{M})$ can push the initial observer state off the manifold. The frequency of these events, however, is not known a priori.\par

\begin{figure}[htb!]
    \centering
    \includegraphics[height=5.5cm]{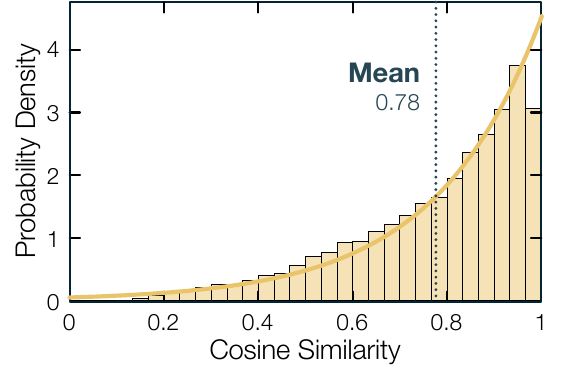}
    \caption{PDF of the cosine similarity between the full gradient and the gradient projected onto the tangent space of $\mathcal{M}$, quantifying the extent to which optimization directions point off the manifold.}
    \label{fig: grad tangent}
\end{figure}

We numerically estimate this frequency, using the basis $\boldsymbol{Q}$ to project gradients into the local tangent space and computing the cosine similarity
\begin{equation*}
    \mathrm{CS}_{\boldsymbol{g}}
    = \frac{\boldsymbol{g}_{\boldsymbol{u}}^\top \boldsymbol{Q}\boldsymbol{Q}^\top \boldsymbol{g}_{\boldsymbol{u}}}{\|\boldsymbol{g}_{\boldsymbol{u}}\|_2\,\|\boldsymbol{Q}\boldsymbol{Q}^\top \boldsymbol{g}_{\boldsymbol{u}}\|_2},
\end{equation*}
which is unity when $\boldsymbol{g}_{\boldsymbol{u}} \in T_{\boldsymbol{u}_{\boldsymbol{\theta}}}(\mathcal{M})$ and zero when $\boldsymbol{g}_{\boldsymbol{u}} \in N_{\boldsymbol{u}_{\boldsymbol{\theta}}}(\mathcal{M})$. We evaluate $\mathrm{CS}_{\boldsymbol{g}}$ for 20 random initial conditions in the $L = 22$ domain, using 1000 random initial guesses for each reference state. Gradients are computed for the $m_x = 16$, $m_t = 16$, and $L = 22$ case with $T = 20$. The resulting PDF is plotted in Fig.~\ref{fig: grad tangent}. The distribution is strongly skewed toward unity, with a mean of 0.78, indicating that gradients are usually well aligned with the tangent space, although cases with substantial normal components do occur.\par

% Effects on optimization
\subsection{Well posedness of variational state estimation}
\label{sec: DoF: posedness}
The preceding sections examined key geometric factors that influence reconstructions. We now shift to a direct assessment of how variational state estimation transitions from an ill-posed problem to a well-posed one as the immersion and embedding criteria are satisfied. These criteria do not determine the behavior of numerical optimization, per se: rank-deficient Hessians, negative curvature, and vanishing gradients (all analyzed in the next section) can obscure the observability of the reference system in practice, independent of $d_\mathcal{M}$ or the sensor configuration. Nevertheless, when optimization is stable, embedding criteria should govern the ``posedness'' of state estimation, as can be seen through the $\varepsilon^*$ metric. When $\varepsilon^*$ is small, many distinct trajectories yield nearly indistinguishable measurements, so a low loss need not imply an accurate reconstruction. To evaluate the practical onset of well-posed reconstruction, therefore, we conduct an empirical survey of reconstruction accuracy and compare these results with the theoretical criteria from Sec.~\ref{sec: DoF: critical points}.\par

\begin{figure}[htb!]
    \centering
    \includegraphics[height=6.5cm]{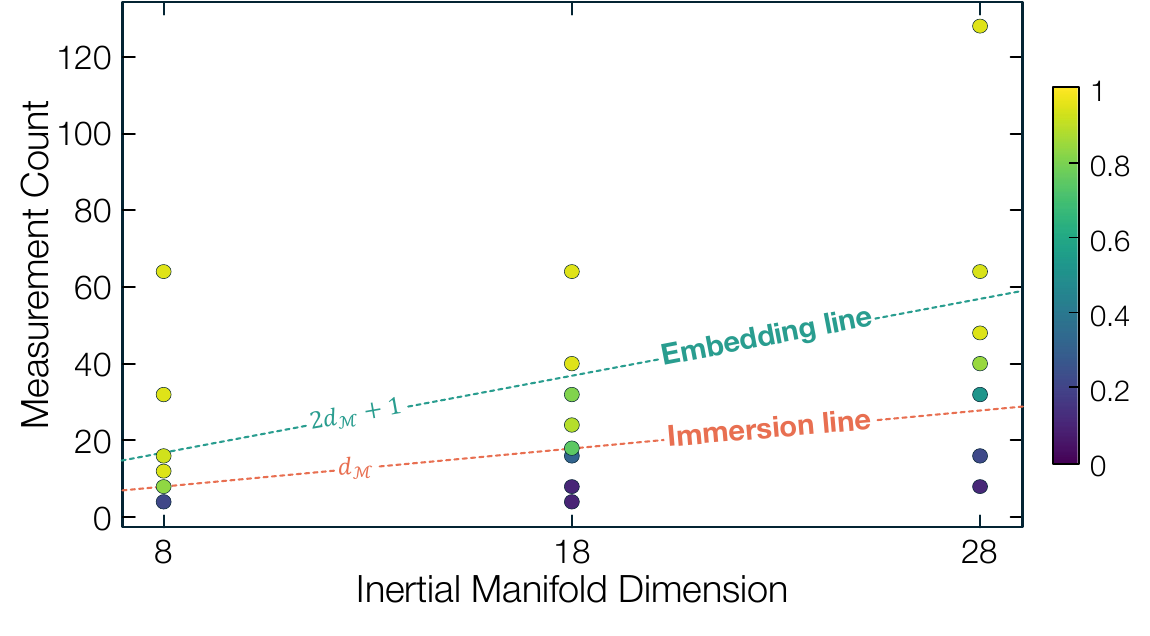}
    \caption{Summary of embedding quality across domain lengths and normalized measurement counts. Reconstructions are classified as accurate for $\mathrm{CS}_{\boldsymbol{U}} \geq 0.95$. Accuracy is low below the immersion line, high above the embedding line, and transitions smoothly between them, consistent with theory.}
    \label{fig: m vs dM}
\end{figure}

To isolate degree-of-freedom effects from the optimization dynamics, we exclude cases with poor convergence, i.e., those with a final loss above $10^{-3}$ (or a mean pointwise error over 3\%, roughly). These cases are limited by failures of optimization rather than the topological relationship between $\mathcal{M}$ and $\mathcal{Y}$. For every $(d_\mathcal{M}, m)$ point in our dataset---spanning all the sensor configurations and reference--observer pairs described in Sec.~\ref{sec: cases: generation}---we compute the probability of successful reconstruction conditioned on $\mathscr{J}<10^{-3}$. Success is defined by $\mathrm{CS}_{\boldsymbol{U}} \geq 0.95$. Figure~\ref{fig: m vs dM} plots these probabilities in the $(d_\mathcal{M}, m)$ plane, along with the immersion line $m = d_\mathcal{M}$ and the embedding line $m = 2d_\mathcal{M} + 1$. Below the immersion line, the chance of an accurate reconstruction collapses to just a few percent. Above the embedding line, the probability approaches unity. Between these bounds, the probabilities vary smoothly with $m$, reflecting a dependence on the specific reference trajectory and initial guess. The structure of Fig.~\ref{fig: m vs dM} supports the applicability of the immersion and embedding criteria. An immersion marks the onset of feasible reconstruction, and an embedding marks the regime in which a low loss reliably corresponds to an accurate reconstruction.\par

\begin{figure}[htb!]
    \centering
    \includegraphics[width=0.9\textwidth]{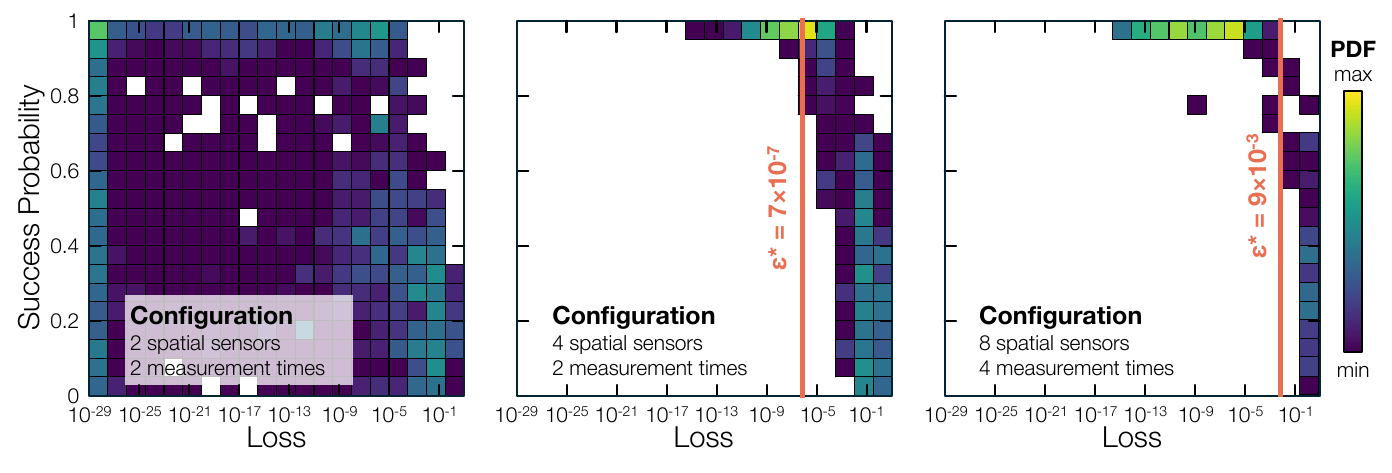}
    \caption{Joint PDF $p(\mathrm{CS}_{\boldsymbol{U}} = \tau, \mathscr{J})$ for various measurement configurations in the $L = 22$ domain. Vertical lines mark $\varepsilon^*$, which is expected to be large whenever the measurement map is an embedding.}
    \label{fig: cos_v_loss_22}
\end{figure}

Next, we look at trends in $\varepsilon^*$ and their relation to the embedding criterion to explicate the relationship between $(d_\mathcal{M}, m)$ and $p(\mathrm{CS}_{\boldsymbol{U}} \geq \tau)$. For each domain size and sensor configuration, we plot the joint density $p(\mathrm{CS}_{\boldsymbol{U}} = \tau, \mathscr{J})$ in the $(\mathscr{J}$, $\tau)$ plane. If a threshold $\varepsilon^* \in [10^{-10}, 10^{-3}]$ exists for $\tau = 0.95$ and $\delta = 0.001$, it is indicated by a vertical line. Figure~\ref{fig: cos_v_loss_22} shows these trends for the $L = 22$ domain, where $d_\mathcal{M} = 8$, using three configurations spanning $m = 4$ to $m = 32$. For $m = 4$, where $\Phi$ cannot be an embedding since $m < d_\mathcal{M}$, there is essentially no relationship between the probability of accurate reconstruction and the loss; $\varepsilon^*$ does not exist for such configurations and the reconstruction problem is hopelessly ill posed. As $m$ increases and eventually exceeds the embedding threshold, a clear correlation between $p(\mathrm{CS}_{\boldsymbol{U}})$ and $\mathscr{J}$ emerges: all the probability mass for $\mathscr{J} < \varepsilon^*$ concentrates near $\mathrm{CS}_{\boldsymbol{U}} = 1$, indicating uniformly accurate reconstructions. With $m = 32$, the $\varepsilon^*$ threshold reaches $9 \times 10^{-3}$.\par

\begin{figure}[htb!]
    \center\includegraphics[width=0.9\textwidth]{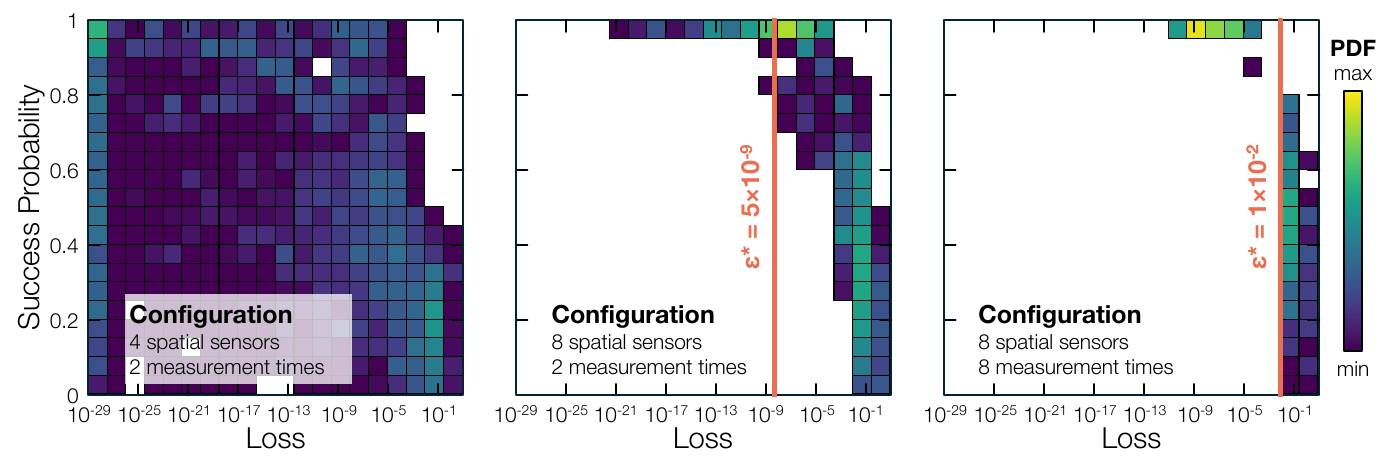}
    \caption{Joint PDF $p(\mathrm{CS}_{\boldsymbol{U}} = \tau, \mathscr{J})$ for measurement configurations in the $L = 44$ domain. Vertical lines mark $\varepsilon^*$, which is expected to be large in cases admitting an embedding.}
    \label{fig: cos_sim_vs_opt_loss_44}
\end{figure}

Figure~\ref{fig: cos_sim_vs_opt_loss_44} presents the same analysis for the $L = 44$ domain, for which $d_\mathcal{M} = 18$. Here, $m = 8$ shows no evidence of an embedding (in contrast with the transitional
behavior observed at $m = 8$ for the $L = 22$ domain). At $m = 16$, the relationship is transitional, although the threshold of $\varepsilon^* \approx 5 \times 10^{-9}$ remains small. At $m = 64$, however, the results evince a robust embedding, with $\varepsilon^* \approx 10^{-2}$, indicating that variational state estimation is theoretically well posed for such sensor arrangements.\par

\begin{figure}[htb!]
    \centering
    \includegraphics[height=5.5cm]{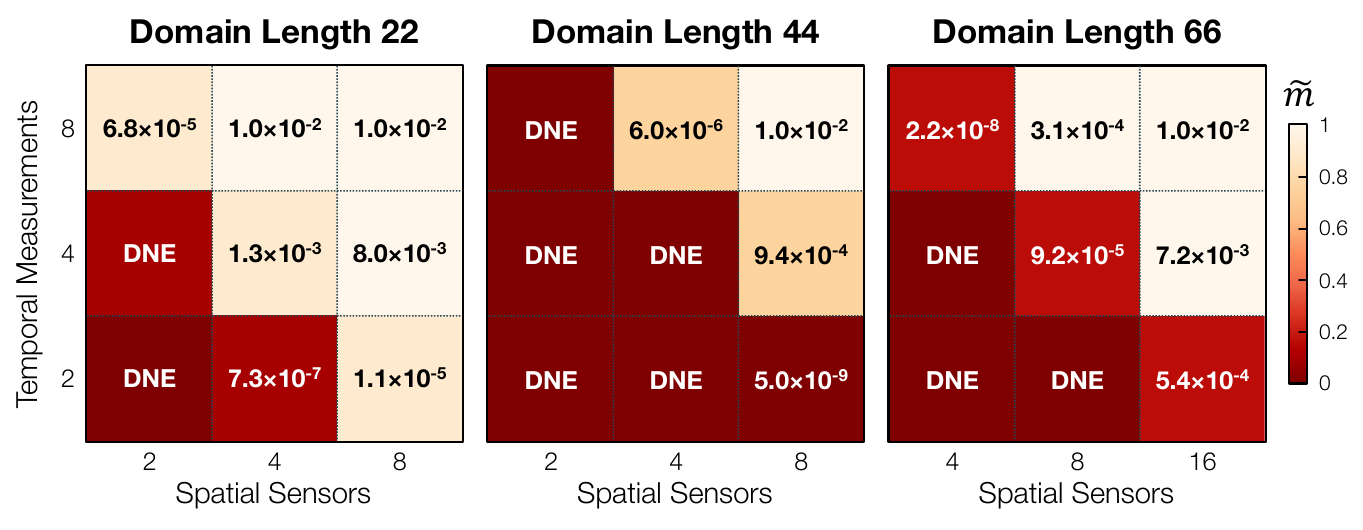}
    \caption{Summary of $\varepsilon^*$ for different domain lengths and normalized measurement counts, plotted in $(m_x, m_t)$ space. Values $\tilde{m} = 1$ mark guaranteed embeddings, $\tilde{m} \in [0,1)$ mark immersions (and possible embeddings), and $\tilde{m} = 0$ cannot support an immersion. Trends in $\varepsilon^*$ are consistent with these classifications.}
    \label{fig: global_emb_qual}
\end{figure}

These trends were computed for all domain sizes and sensor configurations and are summarized in Fig.~\ref{fig: global_emb_qual}. For each configuration, the figure reports the value of $\varepsilon^*$ for $\tau = 0.95$ and $\delta = 0.001$. If no such $\varepsilon^*$ exists above $10^{-10}$, the cell is labeled ``DNE.'' The cell colors indicate a normalized measurement count,
\begin{equation*}
    \widetilde{m} = 
    \begin{cases}
        0, & m < d_{\mathcal M}, \\
        \frac{m - (d_{\mathcal M}-1)}{\,2d_{\mathcal M}+1 - (d_{\mathcal M}-1)\,}, & d_{\mathcal M} < m < 2d_{\mathcal M}+1, \\
        1, & m \geq 2d_{\mathcal M}+1,
    \end{cases}
\end{equation*}
which equals zero below the immersion criterion, increases linearly between the immersion and embedding criteria, and saturates at unity thereafter. Larger values of $\widetilde{m}$ correlate strongly with larger $\varepsilon^*$, mirroring the trends in Fig.~\ref{fig: m vs dM}. Configurations with $\widetilde{m} = 0$, for which $\Phi$ is neither an immersion nor an embedding, either do not admit a computable value of $\varepsilon^*$ or else yield a trivial value. Once $\widetilde{m}$ reaches 1, $\varepsilon^*$ is consistently large, indicating that $\Phi$ acts as an embedding and that the state estimation problem is well posed.\par

% Implications
\subsubsection{Implications of the results}
\label{sec: DoF: posedness: implications}
It is important to note that practical DA problems involve measurement noise, so the exact injectivity of $\Phi$ is sufficient in itself for reliable state estimation. Given noisy data, the conditioning of the embedding becomes critical because small perturbations to the data should not yield disparate reconstructions. The $\varepsilon^*$ metric helps to quantify robustness to noise: larger values imply that nearby measurements are more likely to originate from nearby trajectories. The results in this section therefore identify $d_\mathcal{M}$ as a fundamental lower bound for reconstruction. If $m < d_\mathcal{M}$, then $\Phi$ cannot be an immersion, so even in the limit of vanishing noise, there exist nonzero tangent directions on $\mathcal{M}$ that produce no first-order change in the measurements. In this regime, local injectivity is lost and the inverse problem is said to be topologically ill-posed.\par

The embedding threshold should therefore be interpreted as a baseline for robust state estimation from noisy data. It does not guarantee performance in the presence of noise, but global observability in the noise-free limit is a necessary starting point for reliable reconstruction when noise is present. An important question to ask is whether performance continues to improve for $m > 2d_\mathcal{M} + 1$. Our results suggest that it does indeed: $\varepsilon^*$ generally increases with $m$, as shown in Fig.~\ref{fig: global_emb_qual}, indicating that the inverse map becomes more stable with additional measurements. This conclusion is further supported by the trends in condition number plotted in Fig.~\ref{fig: AE immersion}. However, the relationship is nonlinear and exhibits diminishing returns with $m$. Moreover, $\varepsilon^*$ is not determined by $m$ alone: it also depends on the spatial and temporal arrangement of the measurements, an effect not captured by embedding theory.\par

Although noise prevents any strict interpretation of $2d_\mathcal{M} + 1$ as the point at which state estimation becomes well posed, our results suggest that reconstruction performance should still scale with the intrinsic dimension of the dynamics. This observation points toward empirical scaling laws for practical DA problems. For example, in Navier--Stokes flows, if the effective system dimension scales with a relevant non-dimensional parameter (e.g., the Reynolds number \cite{Constantin1985, Temam2012, Cleary2025}), and if DA performance scales with that dimension, then one may be able to relate known flow parameters to the measurement density required for accurate reconstruction.\par

%%% Optimization Dynamics %%%
\section{Optimization dynamics}
\label{sec: optimization}
Numerical optimizations can still fail when $\Phi$ is an embedding, due in large part to degeneracy of the Hessian (when computed in finite precision), negative curvature, or vanishing gradients. We now analyze \emph{optimization-dynamic effects} caused by these issues and their role in determining whether a topologically well-posed reconstruction problem is numerically tractable. To start, we show the ways in which chaotic dynamics cause variational state estimation to fail. We then use analytical and empirical results to demonstrate why the Hessian creates a poorly condition loss landscape both near optimality and far away from it. These observations motivate the use of NCN, which leverages curvature information from the Hessian to rescale gradients and select effective search directions. Lastly, we present an upper bound for the loss reduction from a single NCN step, and we show numerical results that clarify when and why the optimizer becomes trapped.\par

% Equations
\subsection{Expressions for the gradient and Hessian of \texorpdfstring{$\mathscr{J}$}{the loss}}
\label{sec: optimization: equations}
To frame our discussion of optimization dynamics, we begin by recalling our loss function, which equals
\begin{equation}
        \mathscr{J} = \frac{1}{2}
        \sum_{k \in \mathcal{K}}
        \left[ \mathsf{f}^k(\boldsymbol{u}_{\boldsymbol{\theta}}) - \boldsymbol{u}_k
        \right]^\top \boldsymbol{M}_k \left[ \mathsf{f}^k(\boldsymbol{u}_{\boldsymbol{\theta}}) - \boldsymbol{u}_k \right]
    \end{equation}
    up to a constant where $\boldsymbol{M}_k\in \mathbb{R}^{n\times n}$ is a binary diagonal matrix selecting measurement positions at time index $k$,
    \begin{equation*}
    \label{eq: M_k def}
        M_{k,ii} =
        \begin{cases}
            1, & \text{if $x_i \in \mathcal{X}$ and $k\Delta t \in \mathcal{T}$}, \\
            0, & \text{otherwise},
        \end{cases}
    \end{equation*}
    and $x_i$ is the spatial position corresponding to the $i$th cell. 
    
    We differentiate it to obtain
    \begin{equation}
        \label{equ: MSE grad}
        \boldsymbol{g} =
        \sum_{k \in \mathcal{K}} \boldsymbol{g}_k =
        \sum_{k \in \mathcal{K}}
        \left( \frac{\partial \boldsymbol{u}_{\boldsymbol{\theta}}}{\partial \boldsymbol{\theta}} \right)^\top \boldsymbol{J}_k^\top
        \boldsymbol{M}_k
        \left[ \mathsf{f}^k(\boldsymbol{u}_{\boldsymbol{\theta}}) - \boldsymbol{u}_k \right] \!,
    \end{equation}
where $\mathcal{K} = \{0, \dots, K\}$ and $\boldsymbol{J}_k = \partial \mathsf{f}^k / \partial \boldsymbol{u}_{\boldsymbol{\theta}}$ is the Jacobian of the flow map at time index $k$. The vector $\boldsymbol{g}_k$ is the contribution to the gradient from time $k$, and the full gradient $\boldsymbol{g}$ is simply the sum of these ``sub-gradients.'' Finally, we note that $\partial \boldsymbol{u}_{\boldsymbol{\theta}}/ \partial \boldsymbol{\theta}$ is the inverse discrete Fourier transform, and we have $\partial^2 \boldsymbol{u}_{\boldsymbol{\theta}} / \partial \boldsymbol{\theta}^2 = \boldsymbol{0}$.\par

Equation~\eqref{equ: MSE grad} may be differentiated once more to obtain the Hessian,
\begin{equation}
\label{equ: MSE Hessian}
    \boldsymbol{H} =
    \underbrace{\left(
        \frac{\partial \boldsymbol{u}_{\boldsymbol{\theta}}}{\partial {\boldsymbol{\theta}}} \right)^\top
        \left(\sum_{k \in \mathcal{K}}
        \boldsymbol{J}_k^\top \boldsymbol{M}_k \boldsymbol{J}_k\right)
        \frac{\partial \boldsymbol{u}_{\boldsymbol{\theta}}}{\partial {\boldsymbol{\theta}}}
    }_{\sdx{\boldsymbol{H}}[\mathrm{GN}]} + 
    \underbrace{\left(
        \frac{\partial \boldsymbol{u}_{\boldsymbol{\theta}}}{\partial {\boldsymbol{\theta}}} \right)^\top
        \sum_{k \in \mathcal{K}}
        \frac{\partial (\boldsymbol{J}_k^\top)}{\partial \boldsymbol{u}_{\boldsymbol{\theta}}}
        \boldsymbol{M}_k \left[\mathsf{f}^k(\boldsymbol{u}_{\boldsymbol{\theta}}) - \boldsymbol{u}_k\right]
    }_{\sdx{\boldsymbol{H}}[\mathrm{C}]}
\end{equation}
where $\sdx{\boldsymbol{H}}[\mathrm{GN}]$ is the positive semidefinite Gauss--Newton component and $\sdx{\boldsymbol{H}}[\mathrm{C}]$ is the second-order term stemming from the curvature of the flow map. Near optimality, the residuals $\mathsf{f}^k(\boldsymbol{u}_{\boldsymbol{\theta}}) - \boldsymbol{u}_k$ vanish and the Hessian reduces to the positive semidefinite component $\boldsymbol{H} \to \sdx{\boldsymbol{H}}[\mathrm{GN}]$.\par

% Problems
\subsection{Optimization failure modes in variational state estimation}
\label{sec: optimization: failure modes}
When a reconstruction fails, it is either because $\Phi$ is not an embedding---so that $\varepsilon^{*}$ is extremely small (or not computable) and many trajectories of low loss exhibit large error---or because the optimizer fails to attain a sufficiently low loss. Section~\ref{sec: DoF: posedness} shows that $\varepsilon^{*}$ increases with $m$, but $\varepsilon^{*}$ characterizes the probability that a low loss yields an accurate reconstruction, not the probability of attaining low loss in the first place. It is therefore natural to ask whether the latter probability also increases with $m$. This would be intuitive, and results in Sec.~\ref{sec: DoF: critical points} demonstrate that the conditioning of $\Phi$ improves with added measurements, which \emph{should} increase the likelihood of successful optimization, though this hypothesis must still be verified. We address this question by comparing the probability of achieving a low loss to that of achieving an accurate reconstruction. Trends in these probabilities reveal distinct failure modes of variational state estimation, which are explained in the remainder of Sec.~\ref{sec: optimization}.\par

\begin{figure}[htb!]
    \centering
    \includegraphics[height=6.5cm]{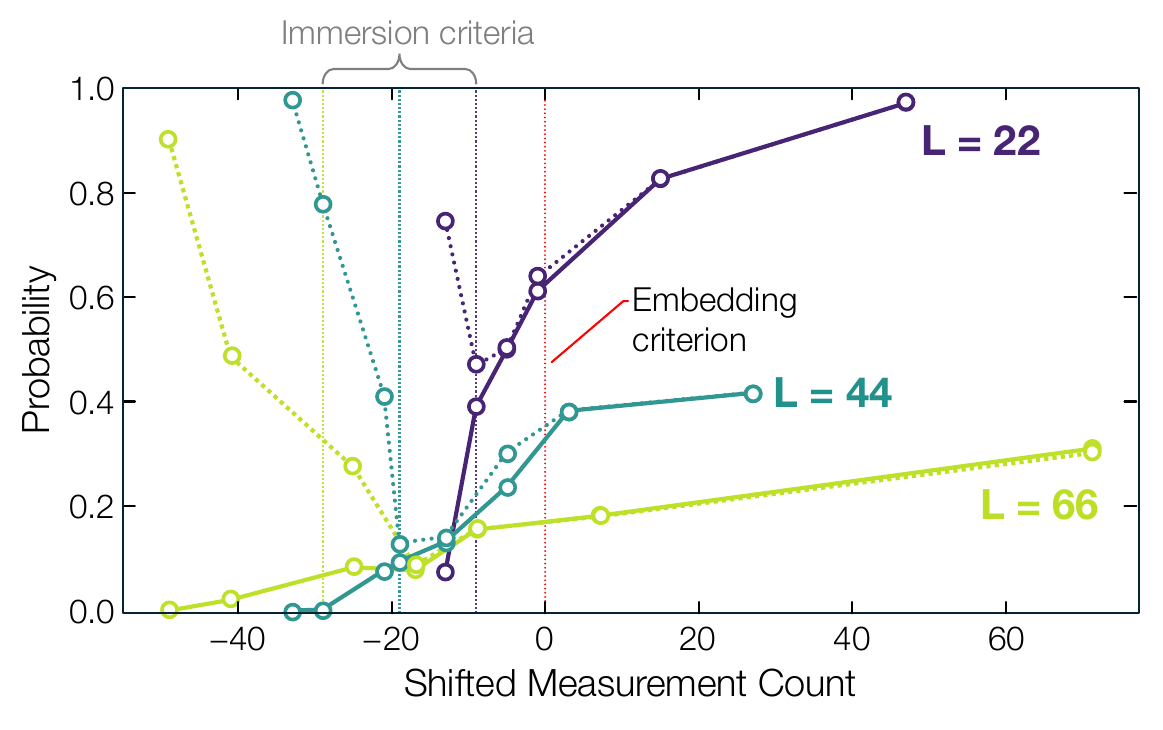}
    \caption{Probabilities of achieving high accuracy ($\mathrm{CS}_{\boldsymbol{U}} \geq 0.95$, solid lines) and low loss ($\mathscr{J} < 10^{-3}$, dashed lines) plotted against $m - (d_\mathcal{M} + 1)$ for $L \in \{22, 44, 66\}$. Vertical lines show the embedding ($x = 0$) and immersion thresholds. All cases use far-off initial conditions with $D_{ij} \in [0.8, 1]$.}
    \label{fig: global results}
\end{figure}

Figure~\ref{fig: global results} shows the probability of $\mathrm{CS}_{\boldsymbol{U}} \geq 0.95$ (solid lines) and $\mathscr{J} < 10^{-3}$ (dashed lines) for all three domains $L \in \{22, 44, 66\}$. Results are plotted against $m - (2d_\mathcal{M} + 1)$, so that $x = 0$ corresponds to the embedding criterion for all domains (marked by a red vertical line). Immersion thresholds are shown as vertical dotted lines that are color-coded by $L$. To emphasize challenging cases where negative curvature and departures from $\mathcal{M}$ are more likely, we report results for cases with an initial distance in $D_{ij} \in [0.8, 1]$.\par

The trends in Fig.~\ref{fig: global results} are consistent across all three domains. The probability of accurate reconstruction (taken as $\mathrm{CS}_{\boldsymbol{U}} \geq 0.95$) begins near zero and increases with $m$, whereas the probability of achieving a low loss exhibits a ``U'' shape: initially high, dipping to a minimum between the immersion and embedding lines, and rising again thereafter.

The initial drop in $p(\mathscr{J} < 10^{-3})$ occurs because, when $\Phi$ is not an embedding and $m$ is very low, adding measurements eliminates spurious minima from the loss landscape, making cases of low loss rarer. As $m$ continues to increase and crosses the immersion and embedding thresholds, the loss and accuracy curves converge, and low loss becomes a reliable indicator of accurate reconstructions thereafter (i.e., $\varepsilon^*$ increases). Two effects drive this transition. First, $\varepsilon^*$ rises with increasing $m$ because observation vectors $\boldsymbol{y}$ from different states on $\mathcal{M}$ become more distinct, making accurate reconstructions more likely, even at moderate loss levels. This corresponds to topological well-posedness. Second, the probability of attaining low loss itself increases because further measurements improve the conditioning of the problem. In short, gains in accuracy at low $m$ are due to degree-of-freedom effects that govern the mapping from $\mathcal{M}$ to $\mathcal{Y}$, whereas gains at high $m$ come from better performance of the optimizer. Still, the absolute probability $p(\mathscr{J} < 10^{-3})$ remains low for large-$L$ domains, even at high $m$, which underscores the need to understand why optimizations fail in the embedding regime (i.e., for topologically well-posed problems).\par

% Condition of loss landscape
\subsection{Condition and curvature of the loss landscape near optimality}
\label{sec: optimization: condition}
\begin{figure}[htb!]
    \center\includegraphics[height=6.5cm]{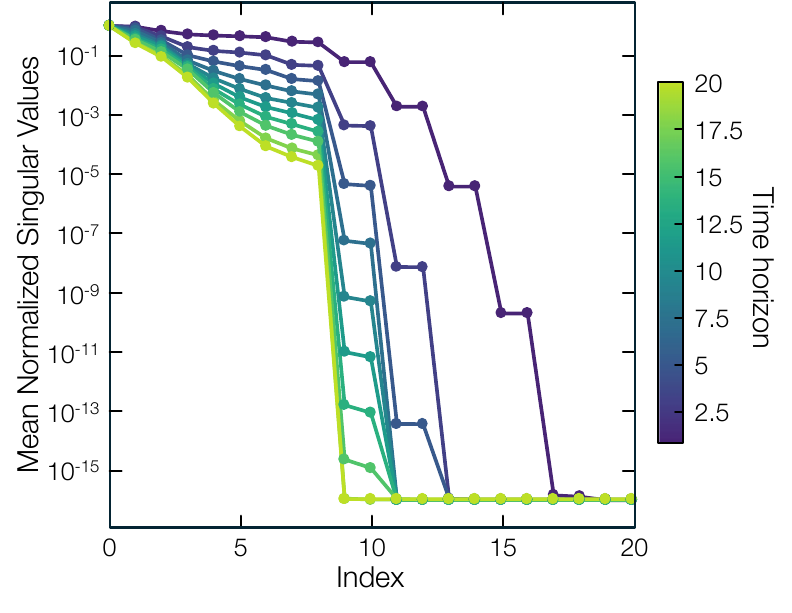}
    \caption{Normalized singular value spectra of $\boldsymbol{J}_k$ averaged over 1000 initial conditions from the $L = 22$ dataset. Line color shows the time horizon $T \in [1, 20]$.}
    \label{fig: Jac Cond}
\end{figure}

Equations~\eqref{equ: MSE grad} and \eqref{equ: MSE Hessian} reveal how the flow map Jacobian $\boldsymbol{J}_k$ (defined by $\delta \boldsymbol{u}_k = \boldsymbol{J}_k \delta \boldsymbol{u}_0$) affects gradients and curvature of the loss landscape. The singular values of $\boldsymbol{J}_k$ are directly determined by the Lyapunov spectrum, i.e., the $i$th singular value scales as $\sigma_i \sim e^{\ell_i k \Delta t}$, where $\ell_i$ is the $i$th Lyapunov exponent \cite{Oseledec1968}. As the assimilation window becomes longer, the computed Jacobian $\boldsymbol{J}_k$ becomes severely ill-conditioned: it quickly becomes singular and its nullity grows steadily with $T$. Figure~\ref{fig: Jac Cond} illustrates this behavior via normalized singular value spectra of $\boldsymbol{J}_k$ for time horizons of increasing duration, averaged over 1000 initial conditions for the $L = 22$ domain.\par

The rapid divergence of the singular values of $\boldsymbol{J}_k$ creates a fundamental trade-off in variational state estimation. Increasing $T$ initially reduces redundancy among measurements and increases $\varepsilon^*$, but the exponential amplification of perturbations progressively weakens the numerical link between early measurements and later ones. Beyond the Lyapunov time, this amplification dominates and the state estimation problem becomes badly ill-posed. A second complication arises from the unweighted MSE loss: because the operator norm $\|\boldsymbol{J}_k\|_2^2$ grows exponentially with $k$, later measurements have a disproportionate influence on the gradient, causing the optimizer to match observations at the end of the assimilation window first, as observed in prior studies on adjoint--variational state estimation \cite{Li2020, Zaki2021}. Consequently, extending the window beyond $T_\ell$, even when adding more measurements, can hinder optimization, since numerical errors accumulate, gradients with respect to $\boldsymbol{u}_0$ become unreliable, and corrupted gradient components from later times dominate the step direction.\par

Ill-conditioning of $\boldsymbol{J}_k$ also manifests in the curvature of the loss. Near optimality, the Hessian is dominated by $\sdx{\boldsymbol{H}}[\mathrm{GN}]$, whose rank is at most $m$. Because $\sdx{\boldsymbol{H}}[\mathrm{GN}]$ depends quadratically on $\boldsymbol{J}_k$, any ill-conditioning in $\boldsymbol{J}_k$ is inherited by and amplified in the Hessian. In practice, the rank of $\boldsymbol{H}$ is usually less than $m$; in cases where it manages to attain rank $m$, its effective condition number (excluding the null space) remains extremely large. These features can slow or stall the optimizer near critical points and can limit the observability of $\boldsymbol{u}_0$, even when the initial guess lies arbitrarily close to it.\par

% Hessian spectrum
\subsection{Curvature of the loss landscape away from optimality}
\label{sec: optimization: eigenvalues}
The previous section shows that for the KS systems of interest, the flow map Jacobian guarantees that $\boldsymbol{H}$ has null eigenvalues at optimality, even for short assimilation windows, thereby limiting observability. We now turn to the prevalence of negative eigenvalues away from optimality, which indicate directions of negative curvature. We show that they are ubiquitous in regions of moderate loss. To do this, we first develop a mathematical intuition for why negative eigenvalues arise, and we then present numerical evidence to support this result.\par

Far from optimality, $\sdx{\boldsymbol{H}}[\mathrm{C}]$ becomes important in Eq.~\eqref{equ: MSE Hessian}. Because the residual $\mathsf{f}^k(\boldsymbol{u}_{\boldsymbol{\theta}}) - \boldsymbol{u}_k$ has no preferred sign, $\sdx{\boldsymbol{H}}[\mathrm{C}]$ is indefinite in expectation and possesses both positive and negative eigenvalues. Meanwhile, $\sdx{\boldsymbol{H}}[\mathrm{GN}]$ has a non-trivial null space due to singularity of the flow map Jacobian and the sparsity of observations, with a nullity that almost always exceeds $n-m$ when $m < n$ (as is always the case in practical state estimation problems). To understand how $\sdx{\boldsymbol{H}}[\mathrm{GN}]$ and $\sdx{\boldsymbol{H}}[\mathrm{C}]$ contribute to the eigenvalues of $\boldsymbol{H}$, let $\lambda_i(\cdot)$ denote the $i$th ordered eigenvalue, with $\lambda_1 \geq \cdots \geq \lambda_n$. Weyl's inequality provides tight bounds on the eigenvalues of a sum of symmetric matrices. Applied to $\sdx{\boldsymbol{H}}[\mathrm{GN}]$ and
$\sdx{\boldsymbol{H}}[\mathrm{C}]$, it yields
\begin{equation}
    \lambda_i\! \left(\sdx{\boldsymbol{H}}[\mathrm{GN}]\right) +
    \lambda_n\! \left(\sdx{\boldsymbol{H}}[\mathrm{C}]\right)
    \leq
    \lambda_i\! \left({\sdx{\boldsymbol{H}}[\mathrm{GN}]} + \sdx{\boldsymbol{H}}[\mathrm{C}]\right)
    \leq
    \lambda_i\! \left({\sdx{\boldsymbol{H}}[\mathrm{GN}]}\right) +
    \lambda_1\! \left({\sdx{\boldsymbol{H}}[\mathrm{C}]}\right) \!.
\end{equation}
For many indices $i$, the Gauss--Newton term satisfies
$\lambda_i(\sdx{\boldsymbol{H}}[\mathrm{GN}]) \approx 0$. Substituting such an index into the inequality gives
\begin{equation}
    \label{equ: Weyl}
    \lambda_i\! \left(\sdx{\boldsymbol{H}}[\mathrm{GN}] + \sdx{\boldsymbol{H}}[\mathrm{C}]\right) =
    \lambda_i\! \left(\boldsymbol{H}\right)
    \in
    \left[ \lambda_n\! \left(\sdx{\boldsymbol{H}}[\mathrm{C}]\right) \!,
    \lambda_1\! \left({\sdx{\boldsymbol{H}}[\mathrm{C}]} \right) \right] \!,
\end{equation}
Since $\lambda_n\! \left(\sdx{\boldsymbol{H}}[\mathrm{C}]\right) < 0 < \lambda_1\! \left(\sdx{\boldsymbol{H}}[\mathrm{C}]\right)$, where $\lambda_1$ and $\lambda_n$ are comparable in expected magnitude, Eq.~\eqref{equ: Weyl} implies that the corresponding $\lambda_i(\boldsymbol{H})$ will take both positive and negative values. Consequently, away from optimality, where the residual is non-zero and often large, negative eigenvalues of $\boldsymbol{H}$ arise with high probability and are expected to be prevalent.\par

\begin{figure}[htb!]
    \centering
    \includegraphics[height=6.5cm]{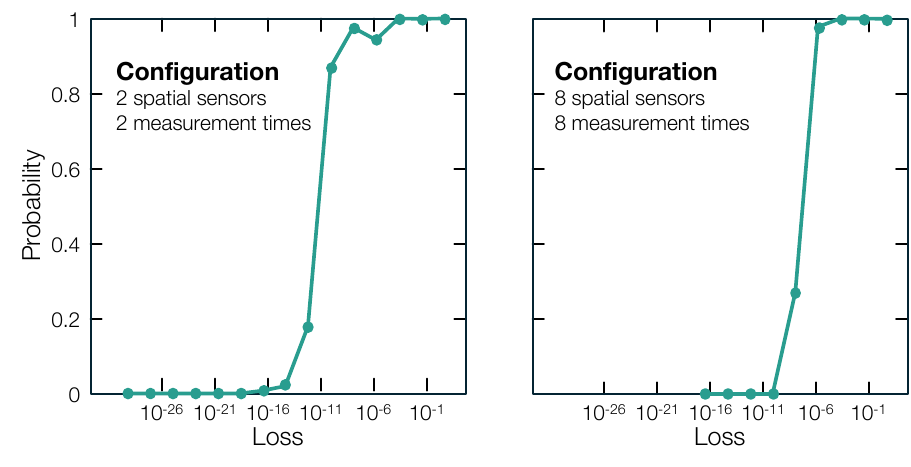}
    \caption{Probability that the terminal Hessian has a negative eigenvalue of magnitude exceeding $10^{-8}$ for the $L = 22$ dataset. Left: $m_x = 2$ and $m_t = 2$. Right: $m_x = 8$ and $m_t = 8$.}
    \label{fig: P(saddle | loss)}
\end{figure}

Figure~\ref{fig: P(saddle | loss)} provides empirical evidence for this claim. It shows the probability of $\boldsymbol{H}$ containing at least one negative eigenvalue below $-10^{-8}$, conditioned on the loss, for the $L = 22$ dataset under two sensor configurations: $m_x = 2$ and $m_t = 2$ and $m_x = 8$ and $m_t = 8$. From this plot, we see that the loss landscape almost always exhibits negative curvature in regions of moderate loss. All remaining sensor configurations across all domain sizes display the same qualitative behavior. Taken together with the severe ill-conditioning of $\boldsymbol{H}$, these observations motivate our use of the NCN optimizer.\par

% Failed convergence
\subsection{When does NCN optimization stall?}
\label{sec: optimization: convergence}
Up to this point, we have established several favorable properties of variational state estimation: when $\Phi$ is an immersion, local reconstruction is well posed; when it is an embedding, the global optimum is the only critical point on $\mathcal{M}$; negative curvature is effectively handled by NCN steps; and pseudo-projection suppresses gradient components that point off $\mathcal{M}$, thereby stabilizing the optimization. Still, NCN optimization with pseudo-projection can proceed very slowly when $\boldsymbol{g}^\top |\boldsymbol{H}|^{-1} \boldsymbol{g}$ is orders of magnitude smaller than $\mathscr{J}$. To explain why, we present an upper bound on the reduction in $\mathscr{J}$ produced by a single NCN step.\par

The upper bound in question begins with a Taylor expansion of the loss increment:
\begin{equation}
    \label{equ: sec_or_taylor_ub}
    \mathscr{J}(\boldsymbol{\theta}_{k+1})
    \leq \mathscr{J}(\boldsymbol{\theta}_k)
    + \boldsymbol{g}^\top\,(\boldsymbol{\theta}_{k+1} - \boldsymbol{\theta}_k)
    + \frac{1}{2} \left(\boldsymbol{\theta}_{k+1} - \boldsymbol{\theta}_k\right)^\top \boldsymbol{H} \left(\boldsymbol{\theta}_{k+1} - \boldsymbol{\theta}_k\right)
    + \frac{M}{6} \|\boldsymbol{\theta}_{k+1} - \boldsymbol{\theta}_k \|_2,
\end{equation}
where $M$ is the Lipschitz constant of $\boldsymbol{H}$. Substituting the NCN step, $\boldsymbol{\theta}_{k+1} = \boldsymbol{\theta}_k - \eta|\boldsymbol{H}|^{-1}\boldsymbol{g}$, yields
\begin{equation}
    \label{equ: NCN_Taylor_bound}
    \mathscr{J}(\boldsymbol{\theta}_{k+1}) - \mathscr{J}(\boldsymbol{\theta}_k) <
    -\underbrace{\eta \boldsymbol{g}^\top |\boldsymbol{H}|^{-1}\boldsymbol{g}}_{\text{first-order term}}+\underbrace{\frac{1}{2} \eta^2 \boldsymbol{g}^\top |\boldsymbol{H}|^{-1} \boldsymbol{H} |\boldsymbol{H}|^{-1} \boldsymbol{g}}_{\text{second-order term}} +
    \underbrace{\frac{M}{6} \|\eta |\boldsymbol{H}|^{-1} \boldsymbol{g}\|_2^3}_{\text{correction term}}.
\end{equation}
Here, first and second order indicates the origin of these terms in the Taylor series. The eigenvectors of $|\boldsymbol{H}|^{-1} \boldsymbol{H} |\boldsymbol{H}|^{-1}$ are the same as those of $\boldsymbol{H}$, and the eigenvalues are similar. Specifically, they are
\begin{equation}
    \lambda_i\! \left(|\boldsymbol{H}|^{-1} \boldsymbol{H} |\boldsymbol{H}|^{-1}\right) =
    \begin{cases}
        1/\lambda_i(\boldsymbol{H}), & \lambda_i(\boldsymbol{H}^{-1}) > \delta, \\
        \lambda_i(\boldsymbol{H})/\delta^2, & \lambda_i(\boldsymbol{H}^{-1})\le\delta.
    \end{cases}
\end{equation}
Since $|\boldsymbol{H}|^{-1}$ is positive definite, the first-order term in Eq.~\eqref{equ: NCN_Taylor_bound} always acts to decrease the loss so long as $\boldsymbol{g} \neq \boldsymbol{0}$ and the step size $\eta$ is sufficiently small. However, the magnitude of the reduction is governed by that of $\boldsymbol{g}^\top |\boldsymbol{H}|^{-1} \boldsymbol{g}$, which must be comparable to $\mathscr{J}$ to ensure meaningful progress. The second-order term can also reduce the loss, but only when $\boldsymbol{g}$ aligns with directions of negative curvature. If $\|\boldsymbol{g}\|_2$ is small relative to the loss, progress becomes extremely slow unless the gradient happens to point along very flat or negatively curved directions of the loss landscape.\par

\begin{figure}[htb!]
    \centering
    \includegraphics[height=6.5cm]{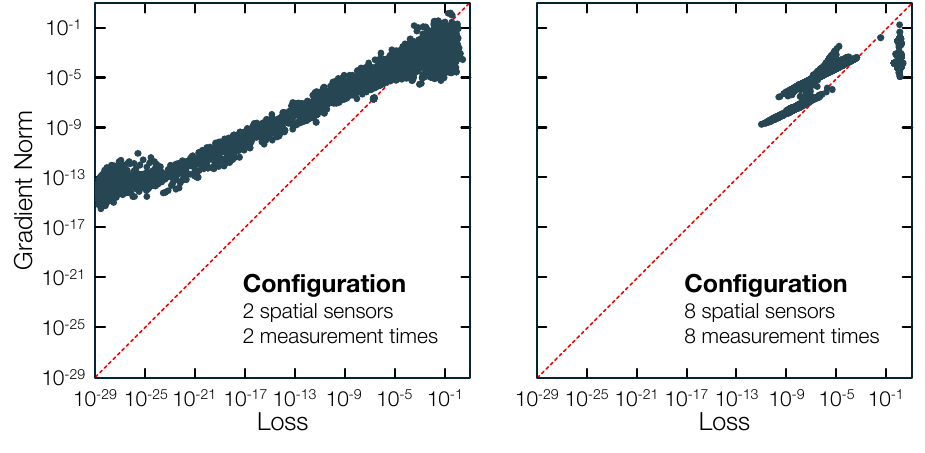}
    \caption{Gradient norm versus optimization loss for $L = 22$. Left: cases with $m_x = 2$ and $m_t = 2$. Right: cases with $m_x = 8$ and $m_t = 8$.}
    \label{fig: grad_vs_loss_L=22}
\end{figure}

\begin{figure}[htb!]
    \center\includegraphics[height=6.5cm]{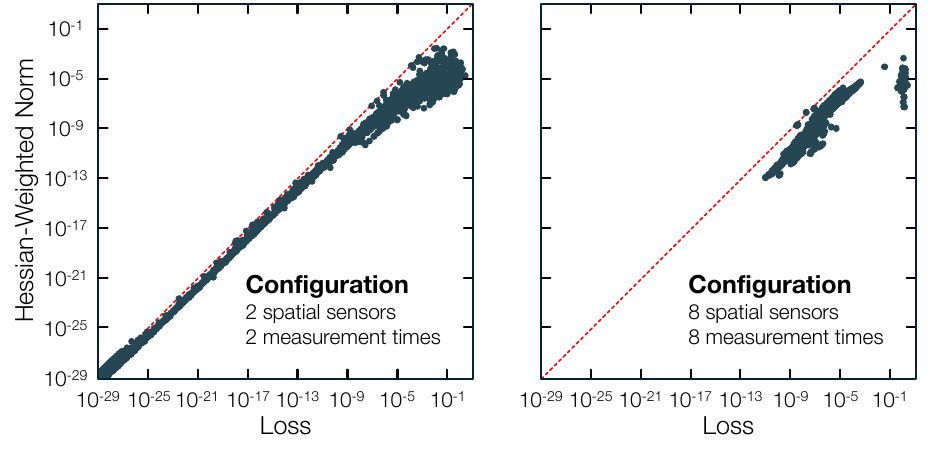}
    \caption{Hessian-weighted gradient norm, $\boldsymbol{g}^\top |\boldsymbol{H}|^{-1} \boldsymbol{g}$, versus loss for $L = 22$. Left: cases with $m_x = 2$ and $m_t = 2$. Right: cases with $m_x = 8$ and $m_t = 8$.}
    \label{fig: g^Tp_vs_loss_L=22}
\end{figure}

We now present evidence that the optimizer does not become trapped in local minima. If it were getting stuck at true critical points, we would expect little or no correlation between the loss and the gradient norm; a plot of $\|\boldsymbol{g}\|_2$ versus $\mathscr{J}$ would show no discernible trend. Instead, Fig.~\ref{fig: grad_vs_loss_L=22}, which plots the gradient norm versus the loss for all cases at $L = 22$ with the $m_x = 2$ and $m_t = 2$ and $m_x = 8$ and $m_t = 8$ sensor configurations, shows a strong, nearly linear relationship on a log--log scale. Across all sensor configurations and domain lengths tested, the minimum correlation between $\log(\mathscr{J})$ and $\log (\|\boldsymbol{g}\|_2)$ is 0.84. This behavior is consistent with a power-law relation of the form $\|\boldsymbol{g}\|_2 \sim \mathscr{J}^a$, where $a$ is a constant. Given this strong correlation, we do not attribute stalled convergence to local minima. Instead, we believe the limiting factor is the regime in which $\boldsymbol{g}^\top |\boldsymbol{H}|^{-1}\boldsymbol{g} \ll \mathscr{J}$, as illustrated in Fig.~\ref{fig: g^Tp_vs_loss_L=22}. This figure uses the same cases as Fig.~\ref{fig: grad_vs_loss_L=22}, but the $y$-axis is replaced by $\boldsymbol{g}^\top |\boldsymbol{H}|^{-1}\boldsymbol{g}$. For the overwhelming majority of points, we observe that $\boldsymbol{g}^\top |\boldsymbol{H}|^{-1}\boldsymbol{g}$ is less than $\mathscr{J}(\boldsymbol{\theta}_k)$, and for cases of high loss, it falls orders of magnitude below $\mathscr{J}$. Finally, we note that the results in Figs.~\ref{fig: grad_vs_loss_L=22} and \ref{fig: g^Tp_vs_loss_L=22} were essentially unchanged after an additional 650 NCN iterations, confirming that the optimizations in these figures are converged.\par

\begin{figure}[htb!]
    \center\includegraphics[height=6.5cm]{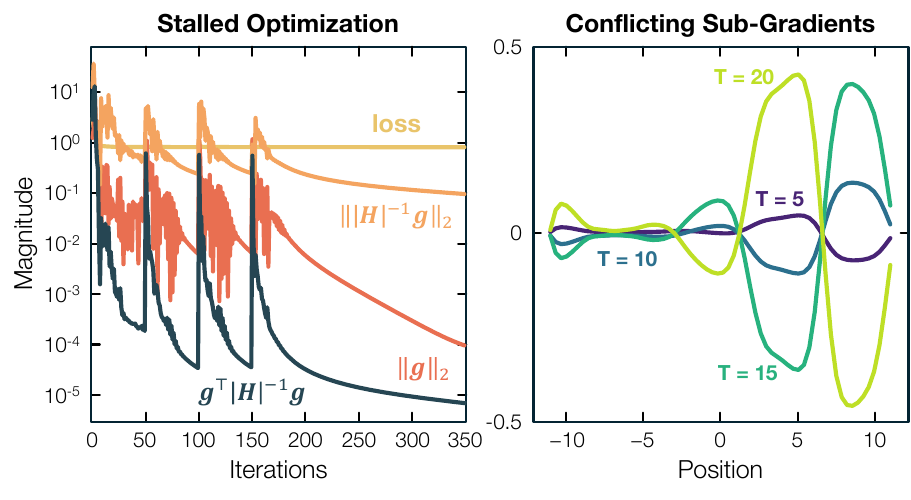}
    \caption{Representative optimization failure caused by vanishing gradients for a case from the $L = 22$ domain. Left: gradient norm $\|\boldsymbol{g}\|_2$, step $\|\,|\boldsymbol{H}|^{-1}\boldsymbol{g}\|_2$, and Hessian-weighted norm $\boldsymbol{g}^\top |\boldsymbol{H}|^{-1}\boldsymbol{g}$ versus iteration. The gradient is several orders of magnitude smaller than the loss, stalling progress. Right: sub-gradients at the final iteration visualized in state space; they nearly cancel when summed to produce the full gradient~$\boldsymbol{g}$.}
    \label{fig: vanishing_grad}
\end{figure}

To illustrate these ideas, we examine the optimizer's behavior for a failed case with $m_x = 4,\ m_t = 4,\ L = 22$. The left panel of Fig.~\ref{fig: vanishing_grad} shows the loss together with the magnitudes of $\boldsymbol{g}$, $|\boldsymbol{H}|^{-1}\boldsymbol{g}$, and $\boldsymbol{g}^\top |\boldsymbol{H}|^{-1}\boldsymbol{g}$ as functions of iteration. The loss barely decreases because $\|\boldsymbol{g}\|_2$ is several orders of magnitude smaller than $\mathscr{J}$, which in turn forces the quadratic term to be even smaller still. A plausible explanation for why the gradient becomes so small relative to the loss is the presence of conflicting sub-gradients. That is, contributions $\boldsymbol{g}_k$ from different observation times remain large at the end of the optimization, but they almost perfectly cancel out in aggregate. The right panel of Fig.~\ref{fig: vanishing_grad} illustrates this effect for the same case shown on the left; sub-gradients from all four observation times at the final iteration are mapped into state space and plotted. Although magnitudes of the individual curves are substantial, their sum, corresponding to the $\|\boldsymbol{g}\|_2$ curve in the left panel (multiplied by $\sqrt{n}$ for the conversion to state space), is nearly zero.\par

Motivated by this observation, we note that minimizing measurement residuals at different times can be viewed as a multi-task optimization problem in which sub-gradients may conflict with one another. This perspective suggests that techniques from multi-task learning, such as dynamic loss weighting \cite{Kendall2018} or gradient-conflict resolution \cite{Yu2020}, could improve the global behavior of optimizers in variational state estimation.\par

%%% Conclusions %%%
\section{Conclusions and outlook}
\label{sec: conclusions}
Variational state estimation provides a powerful framework for combining simulations and experiments to reconstruct high-fidelity trajectories of a dynamical system that are anchored to real-world observations. However, the number of measurements required for accurate reconstruction of a chaotic system is not known a priori. We address that gap using tools from dynamical systems and embedding theory. For a dissipative system whose attractor lies on an inertial manifold $\mathcal{M}$ of dimension $d_\mathcal{M}$, we show that $m \geq d_\mathcal{M}$ measurements are sufficient for local observability from an arbitrarily good initial guess, and $m \geq 2d_\mathcal{M} + 1$ are required for global observability on $\mathcal{M}$. These classical bounds determine whether the observation map $\Phi$ is an immersion or an embedding, respectively, guaranteeing the local or global existence of $\Phi^{-1}$. While such criteria are well established in the literature on state space reconstruction, we demonstrate their applicability to variational state estimation by translating the existence and conditioning of $\Phi^{-1}$ into geometric conditions for a well-posed reconstruction problem.\par

Specifically, we show that when $\Phi$ is an embedding, the global optimum is the only critical point of the loss landscape on $\mathcal{M}$. Moreover, the induced geometry is well conditioned, such that similar measurements correspond to nearby states and dissimilar measurements to distant states. These theoretical results are validated through extensive simulations of Kuramoto--Sivashinsky systems in domains of length 22, 44, and 66. In general, when the embedding criterion is satisfied, variational state estimation becomes well posed and reconstruction accuracy improves steadily with increasing $m$.\par

Despite these guarantees, embedding theory alone does not determine practical limits on observability. Even when $\Phi$ is an embedding, the loss landscape can remain severely ill conditioned due to singularity of the flow map Jacobian and sparsity of the sensor configuration. As a result, trajectories may be theoretically observable yet difficult to recover by numerical means. Two key challenges are identified. First, gradients can have large components normal to $\mathcal{M}$, causing iterates in an optimization to drift off the manifold, at which point guarantees from embedding theory no longer apply. To counteract this, we introduce a pseudo-projection step that periodically pulls the estimate back toward $\mathcal{M}$, helping to stabilize the reconstruction. Second, the Hessian is degenerate at optimality and becomes indefinite away from it. Indeed, directions of negative curvature are ubiquitous in the loss landscape at moderate loss levels. These pathological features undermine first-order, Newton, and quasi-Newton methods alike. To address this, we employ a ``non-convex Newton'' technique that explicitly handles negative curvature while preserving descent directions for indefinite Hessians. When combined with pseudo-projection, NCN enables robust state estimation once the embedding criterion is satisfied.\par

Nevertheless, optimization can still stall when the gradient norm becomes much smaller than the loss. We attribute this behavior to destructive interference among sub-gradients from different observation times, which accounts for all the failures we examined in the embedding regime. Future work will address this limitation by incorporating ideas from multi-task learning, while also extending the framework to more realistic DA problems. Important next steps include accounting for measurement noise and operator error, and determining how reconstruction performance scales with $d_\mathcal{M}$ in their presence. We will also test whether the present results extend beyond KS systems to higher-dimensional flows, including cases such as 3D turbulence where the existence of an inertial manifold is uncertain and the effective dimension of the attractor must be estimated empirically. Finally, future work should consider non-stationary measurement operators, which arise naturally in experimental settings and may alter both the embedding properties and the optimization dynamics of state estimation.\par

% Appendices
\appendix

%%% Appendex A: KSE Numerical Integration %%%
\section{Numerical simulation}
\label{app: numerics}

% Integration
\subsection{Pseudo-spectral scheme with exponential time-differencing}
\label{app: numerics: ETD}
The hyper-diffusion term in the KS equation causes Fourier coefficients associated with high-wavenumber modes to have large values of $\partial \widehat{u}_j / \partial t$, leading to rapid transients. The characteristic time scale of the $j$th Fourier mode scales as $O(j^{-4})$ for large $j$, whereas low-wavenumber modes evolve much more slowly. Differential equations that exhibit such a wide separation of scales are deemed to be ``stiff,'' and stiffness poses major challenges for classical explicit time-stepping schemes.\par

Explicit integrators require a time step that is small enough to stabilize the fastest modes, but simulations must run for long enough to resolve the slow dynamics of low-frequency modes. This combination of small $\Delta t$ and long integration times leads to a high computational cost. To overcome this issue, we employ exponential time-differencing, which analytically integrates the stiff linear terms in the KS equation and numerically integrates the non-linear term \cite{Holland2002, Petropoulos2002, Schuster2000}. This approach enables large time steps and long-time integrations without compromising stability.\par

To start, we discretize the periodic spatial domain $[-L/2, L/2]$ into $n$ uniformly spaced points,
\begin{equation}
    \boldsymbol{u}(t) =
    \left[u(-L/2, t), \dots, u(L/2-\Delta x, t)\right] \!,
\end{equation}
where $\Delta x = L/n$. The semi-discrete KS equation is written as
\begin{equation}
    \label{equ: discrete KSE}
    \frac{\mathrm{d}  \boldsymbol{u}}{\mathrm{d} t} = -\left(\sdx{\boldsymbol{D}}[(2)] + \sdx{\boldsymbol{D}}[(4)]\right) \boldsymbol{u} - \frac{1}{2} \sdx{\boldsymbol{D}}[(1)] \boldsymbol{u}^{\circ 2},
\end{equation}
with $\sdx{\boldsymbol{D}}[(i)]$ being the $i$th-order discrete derivative operator and $(\cdot)^{\circ 2}$ being the element-wise square.\par

Applying the discrete Fourier transform
\begin{equation*}
    \widehat{\boldsymbol{u}} = \mathsf{F}(\boldsymbol{u})
\end{equation*}
yields the KS equation in Fourier space,
\begin{equation}
    \frac{\mathrm{d} \boldsymbol{\widehat{u}}}{\mathrm{d} t} =
    -\left(\sdx{\widehat{\boldsymbol{D}}}[(2)] +
    \sdx{\widehat{\boldsymbol{D}}}[(4)]\right) \boldsymbol{\widehat{u}} - \frac{1}{2} \sdx{\widehat{\boldsymbol{D}}}[(1)] \,\mathsf{F}\! \left(\boldsymbol{u}^{\circ 2}\right) \!,
\end{equation}
where $\sdx{\widehat{\boldsymbol{D}}}[(i)]$ are diagonal derivative operators in Fourier space with entries
\begin{equation*}
    \widehat{D}_{jj}^{(i)} = \left(\mathrm{i} k_j\right)^i.
\end{equation*}
Here, $\mathrm{i}$ is the imaginary unit, $k_j = 2 \pi j^\prime/L$ are wavenumbers, and $j^\prime$ is the signed mode index,
\begin{equation*}
    j^\prime =
    \begin{cases}
        j, & 0 \leq j \leq n/2 \\
        j-n, & n/2 < j < n \\
    \end{cases}.
\end{equation*}
We next define the linear and non-linear terms,
\begin{subequations}
    \begin{align}
        \boldsymbol{C} &= -\sdx{\widehat{\boldsymbol{D}}}[(2)] - \sdx{\widehat{\boldsymbol{D}}}[(4)], \\
        \mathsf{N}\! \left(\boldsymbol{u}\right) &=
        -\frac{1}{2} \sdx{\widehat{\boldsymbol{D}}}[(1)]
        \,\mathsf{F}\! \left(\boldsymbol{u}^{\circ 2}\right) \!,
    \end{align}
\end{subequations}
and we apply the one-third dealiasing rule by zeroing out all modes for which $|j^\prime| > n/3$ when evaluating the non-linear term. With these elements in hand, the KS equation becomes
\begin{equation}
    \frac{\mathrm{d} \boldsymbol{\hat{u}}}{\mathrm{d} t} = \boldsymbol{C}\boldsymbol{\hat{u}} + \mathsf{N}(\boldsymbol{u}).
\end{equation}
Multiplying both sides by $e^{-\boldsymbol{C}t}$ gives
\begin{equation*}
    \frac{\mathrm{d} \boldsymbol{\hat{u}}}{\mathrm{d} t} e^{-\boldsymbol{C}t} - \boldsymbol{C} \boldsymbol{\widehat{u}} \,e^{-\boldsymbol{C}t} = \mathsf{N}(\boldsymbol{u}) \,e^{-\boldsymbol{C}t}.
\end{equation*}
We rearrange this as
\begin{equation*}
    \frac{\mathrm{d}}{\mathrm{d} t}\! \left(e^{-\boldsymbol{C}t} \boldsymbol{\widehat{u}}\right) = \mathsf{N}(\boldsymbol{u}) \,e^{-\boldsymbol{C}t},
\end{equation*}
and integrate it from $t_{k-1}$ to $t_{k}$, where the solver time step is $\Delta t = t_k - t_{k-1}$,
\begin{align*}
    \boldsymbol{\widehat{u}}(t_k) \,e^{-\boldsymbol{C}t_k} -
    \boldsymbol{\widehat{u}}(t_{k-1}) \,e^{ -\boldsymbol{C} t_{k-1}} &=
   \int_{t_{k-1}}^{t_k} \mathsf{N}[\boldsymbol{u}(t)] \,e^{-\boldsymbol{C} t} \,\mathrm{d}t, \\
   \boldsymbol{\widehat{u}}(t_k) \,e^{-\boldsymbol{C}t_k} -
    \boldsymbol{\widehat{u}}(t_{k-1}) \,e^{ -\boldsymbol{C} t_{k-1}} &=
   e^{-\boldsymbol{C}t_{k-1}} \int_{0}^{\Delta t} \mathsf{N}\! \left[\boldsymbol{u}(t_{k-1} + \tau)\right] e^{-\boldsymbol{C} \tau} \,\mathrm{d}\tau.
\end{align*}
Multiplying by $e^{\boldsymbol{C}t_k}$ results in the final expression,
\begin{equation}
    \boldsymbol{\widehat{u}}(t_n) =
    \boldsymbol{\widehat{u}}(t_{n-1}) \,e^{\boldsymbol{C}\Delta t} + e^{\boldsymbol{C}\Delta t} \int_{0}^{\Delta t} \mathsf{N}\! \left[\boldsymbol{u}(t_{n-1} + \tau)\right] e^{-\boldsymbol{C} \tau} \,\mathrm{d}\tau.
\end{equation}
Various numerical schemes can be used to approximate the integral. We employ the fourth-order exponential time-differencing Runge--Kutta method of Cox and Matthews \cite{Cox2002}.\par

% Solver validation
\subsection{Solver validation}
\label{app: numerics: validation}
\begin{figure}[htb!]
    \center\includegraphics[height=6.5cm]{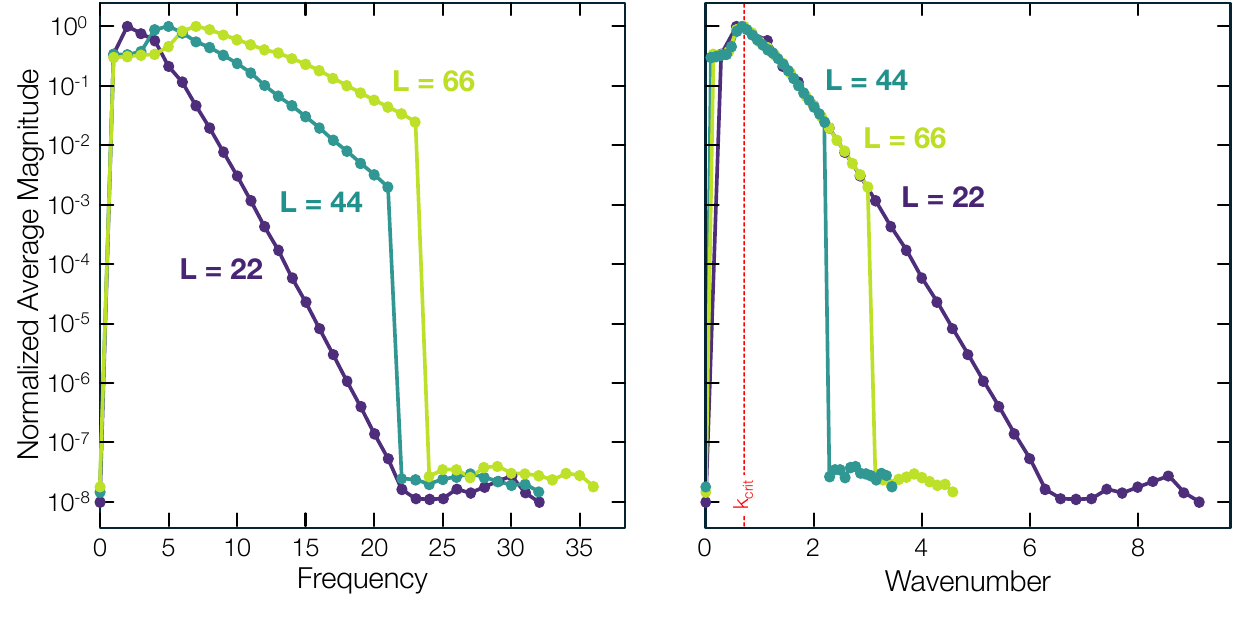}
    \caption{Average Fourier coefficient magnitude plotted versus frequency index (left) and versus wavenumber (right) for long-time trajectories in domains of length 22, 44, and 66.}
    \label{fig: freq num}
\end{figure}

Long-time KS trajectories exhibit predominantly low-frequency content, so only a modest number of spatial nodes are required for accurate simulation. In the literature, Linot et al.~\cite{Linot2020} used $64$ nodes for $L \in \{22, 44, 66\}$, while Cvitanovi{\'c} \cite{Cvitanovic2010} used 32 nodes for $L = 22$. In our numerical experiments, we adopt 64 nodes for $L \in \{22, 44\}$ and 72 nodes for $L=66$. Figure~\ref{fig: freq num} plots the average Fourier-coefficient magnitude versus mode number (left) and versus wavenumber (right) using snapshots from our dataset. These spectra confirm that the chosen spatial resolutions are sufficient to resolve the dynamically relevant frequency content of long-time solutions to the KS equation. In particular, the dominant energy lies at wavenumbers between 0 and 1, with a peak near $k_{\mathrm{crit}}$, as predicted from Eq.~\eqref{equ: KSE Fourier}.\par

\begin{figure}[htb!]
    \center\includegraphics[height=6.5cm]{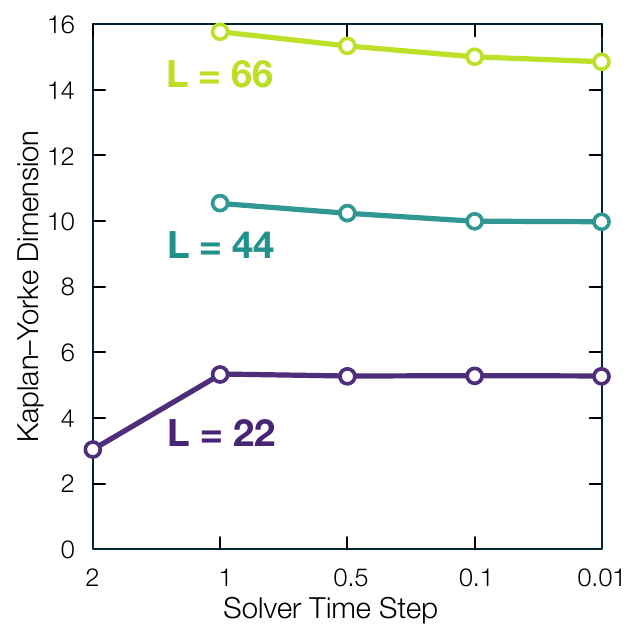}
    \caption{Kaplan--Yorke dimension $d_\mathrm{KY}$ versus solver time step $\Delta t$ for domain lengths $L = 22$, 44, and 66.}
    \label{fig: KY}
\end{figure}

Our main requirement is that the asymptotic statistical properties of our numerical solutions match those of the true dynamics. To evaluate this, Fig.~\ref{fig: KY} shows the Kaplan--Yorke dimension $d_\mathrm{KY}$ computed with our solver for $L \in \{22, 44, 66\}$ as a function of time step. The attractor dimension is nearly invariant with respect to $\Delta t$, indicating that the solver is stable and statistically consistent across a wide range of time steps. Edson et al.~\cite{Edson2019} report $d_\mathrm{KY} = 5.198$ for $L = 22$. Using time steps $\Delta t \in \{0.01, 0.1, 0.5\}$, our solver yields values of $\sim$5.23, in close agreement with Edson and co. Our solver also reproduces the expected linear growth of $d_\mathrm{KY}$ with $L$, as reported in \cite{Edson2019}. These comparisons support our supposition that our solver resolves the long-time statistics well. To balance computational cost and accuracy, we use $\Delta t = 0.1$ for $L \in \{22, 44\}$ and $\Delta t = 0.05$ for $L = 66$. Finally, we note that trajectories produced by our solver yield correct estimates of the IM dimension, per Fig.~\ref{fig: IM dimension}, further validating the fidelity of our scheme.\par

%%% Lyapunov %%%
\section{Lyapunov spectra}
\label{app: Lyapunov}
\begin{figure}[htb!]
    \center\includegraphics[height=6.5cm]{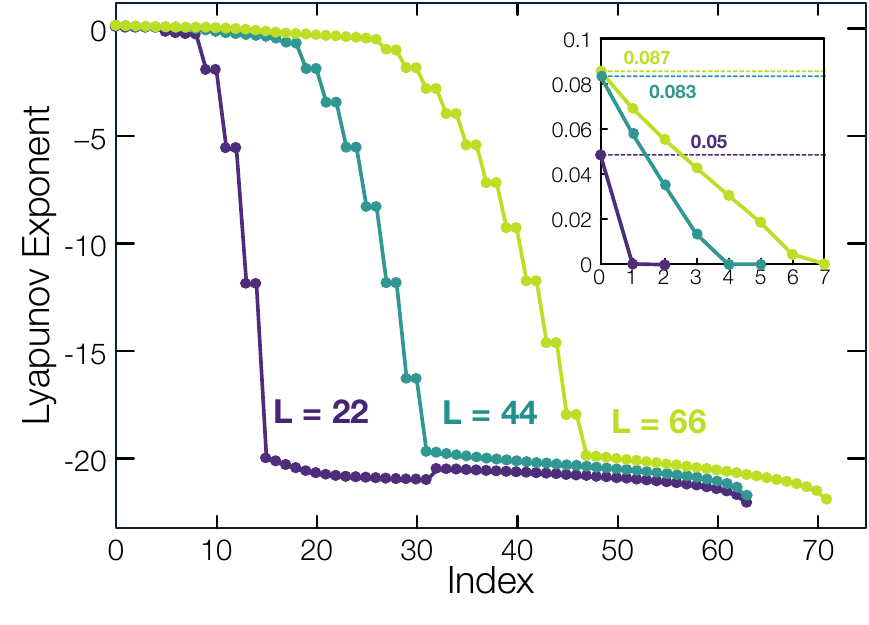}
    \caption{Lyapunov exponent spectra $\ell_i$ for $L \in \{22, 44, 66\}$ computed over a time horizon of $T = 5 \times 10^5$ with reorthogonalization every $2$ time units.}
    \label{fig: ly exp}
\end{figure}

We use the algorithm of Benettin et al.~\cite{Benettin1980} to compute the Lyapunov spectrum for the KS equation (see also Sandri~\cite{Sandri1996}). The method begins with a set of orthonormal tangent vectors $\boldsymbol{Q}_0$. They are advanced forward in time by $k$ steps using the variational equation
\begin{equation}
    \boldsymbol{V}_j = \boldsymbol{J}_k \boldsymbol{Q}_j,
\end{equation}
where $\boldsymbol{J}_k$ is the flow map Jacobian for $k$ time units of advancement. A QR decomposition is applied, $\boldsymbol{V}_j = \boldsymbol{Q}_{j+1} \sdx{\boldsymbol{R}}[(j)]$, and the process is repeated with $\boldsymbol{Q}_{j+1}$. We perform $K$ iterations of this cycle. The $i$th Lyapunov exponent $\ell_i$, which measures the average exponential growth rate of the $i$th most unstable tangent direction, is computed as
\begin{equation}
    \ell_i = \frac{1}{T} 
    \sum_{j=0}^{K-1}
    \log\! \left( |\sdx{R_{ii}}[(j)]| \right) \!,
\end{equation}
where $\sdx{R_{ii}}[(j)]$ is the $i$th diagonal entry of $\sdx{\boldsymbol{R}}[(j)]$ and $T = K \Delta t$ is the integration time. Periodic application of the QR decomposition is essential to prevent the tangent vectors from collapsing onto the dominant mode of $\boldsymbol{J}_k$. Figure~\ref{fig: ly exp} shows Lyapunov spectra that we computed for $L \in \{22, 44, 66\}$, using a total time horizon of $T = 5 \times 10^5$ and performing QR decomposition every 2 time units, i.e., $k = 2 / \Delta t$.\par

%%% Autoencoder %%%
\section{Autoencoder architecture and training}
\label{app: AE}
An autoencoder is a neural network comprising an encoder $\mathsf{E} : \mathcal{M} \to \mathcal{L}$, which for us maps discrete states on the IM, $\boldsymbol{u} \in \mathcal{M} \subset \mathbb{R}^n$, into a lower-dimensional latent space $\mathcal{L}$, and a decoder $\mathsf{D} : \mathcal{L} \to \mathcal{M}$, which approximates the inverse of $\mathsf{E}$. Their composition, $\mathsf{A} = \mathsf{D} \circ \mathsf{E}$, is trained to learn the identity on $\mathcal{M}$ such that all information in $\boldsymbol{u}$ is preserved when compressed into the latent representation. Parameters of the autoencoder are learned by minimizing the loss
\begin{equation}
    \mathscr{J} = \frac{1}{N} \sum_{i=1}^N 
    \left\|\boldsymbol{u}_i - \mathsf{D}\!\left[\mathsf{E}(\boldsymbol{u}_i)\right]\right\|_2^2,
\end{equation}
where $\boldsymbol{u}_i$ denotes the $i$th training sample, with $i = 1, 2, \dots, N$. For each domain, these samples are drawn from the corresponding long-time rollout described in Sec.~\ref{sec: cases: generation}, and the networks are trained for 2000 epochs using the Adam optimizer.\par

The autoencoders used in this work are made up of a sequence of fully connected layers. In addition, at the end of the encoder, we append a ``linear block'' composed of several fully connected linear layers, each with an output dimension equal to the latent dimension and with no biases or activation functions. This block encourages a latent space of low-rank \cite{Zeng2024, Jing2020}. The architectures employed for each domain are summarized in Table~\ref{tab: autoencoder_arch}, where $n$ is the dimension of the state space and $d_\mathcal{L}$ is the latent dimension. We set $d_\mathcal{L}$ to 20 for $L = 22$, to 30 for $L = 44$, and to 50 for $L = 66$, although the results are insensitive to this hyperparameter.\par

\begin{table}[htb!]
    \caption{Autoencoder architecture. Here, $n$ denotes the dimension of the state space, and $d_{\mathcal{L}}$ denotes the dimension of the latent space, set to $20$ for $L=22$, $30$ for $L=44$, and $50$ for $L=66$.}
    \renewcommand{\arraystretch}{1.25}
    \centering\vspace*{.3em}
    \begin{tabular}{c c c c}
        \hline\hline
        \multicolumn{1}{c}{\bf Component} &
        \multicolumn{1}{c}{\bf Input Dim. $\boldsymbol{\to}$ Output Dim.} &
        \multicolumn{1}{c}{\bf Activation} &
        \multicolumn{1}{c}{\bf Bias} \\
        \hline
        \multirow{3}{*}{Encoder}& $n \to 512$ & Swish & Yes \\
        & $512 \to 320$ & Swish & Yes \\
        & $320 \to d_{\mathcal{L}}$ & No & Yes \\[2pt]
        \hline
        \multirow{2}{*}{Linear Block} & $d_{\mathcal{L}}\to d_{\mathcal{L}}$ & No & No \\
        & $d_{\mathcal{L}} \to d_{\mathcal{L}}$& No & No \\[2pt]
        \hline
        \multirow{3}{*}{Decoder}& $d_{\mathcal{L}} \to 320$& Swish & Yes \\
        & $320 \to 512$ & Swish & Yes \\
        & $512 \to n$ & No & Yes \\[2pt]
        \hline\hline
    \end{tabular}
    \label{tab: autoencoder_arch}
\end{table}

In order to use the autoencoder for inference and to estimate $d_\mathcal{M}$, we perform a PCA in the latent space. First, we approximate the mean latent space vector as
\begin{equation}
    \boldsymbol{\ell} = \frac{1}{K}\sum_{k=1}^K \mathsf{E}(\boldsymbol{u}_k),
\end{equation}
where $k$ indicates iterations from a long-time rollout having a total of $K$ snapshots. We then construct the centered data matrix
\begin{equation}
    \boldsymbol{X} =
    [\mathsf{E}(\boldsymbol{u}_0) - \boldsymbol{\ell}, \dots,
    \mathsf{E}(\boldsymbol{u}_K) - \boldsymbol{\ell}]
\end{equation}
and compute its SVD. The number of non-trivial singular values provides an estimate of $d_\mathcal{M}$. During inference, we restrict the latent representation to the dominant subspace by projecting out directions associated with negligible singular values. To do so we build a matrix $\boldsymbol{P} \in \mathbb{R}^{d_\mathcal{L} \times d_\mathcal{M}}$ using the leading $d_\mathcal{M}$ left singular vectors of $\boldsymbol{X}$. The reduced latent coordinate is obtained as
\begin{equation}
    \boldsymbol{z} = \boldsymbol{P}^\top \left[\mathsf{E}(\boldsymbol{u}) - \boldsymbol{\ell}\right] \!,
\end{equation}
with approximate inverse
\begin{equation}
    \boldsymbol{u} = \mathsf{D}\! \left(\boldsymbol{P}\boldsymbol{z} + \boldsymbol{\ell}\right) \!.
\end{equation}\par

%%% Discrete Adjoint %%%
\section{Discrete adjoint systems}
\label{app: adjoint}
In our variational state estimation problem, we seek to minimize the objective
\begin{equation}
    \label{equ: discrete loss alt}
    \mathscr{J} = \sum_{k=0}^K \mathsf{M}_k(\boldsymbol{u}_{\boldsymbol{\theta}, k}, \boldsymbol{u}_k),
\end{equation}
where each term measures the discrepancy between the observer trajectory and the true trajectory at time index $k$. This formulation is equivalent to Eq.~\eqref{equ: discrete loss} up to a constant. In Eq.~\eqref{equ: discrete loss alt}, $\mathsf{M}_k$ is defined as
\begin{equation}
    \mathsf{M}_k(\boldsymbol{u}_{\boldsymbol{\theta}, k}, \boldsymbol{u}_k) = \frac{1}{m}(\boldsymbol{u}_{\boldsymbol{\theta}, k} - \boldsymbol{u}_k)^\top \boldsymbol{M}_k (\boldsymbol{u}_{\boldsymbol{\theta}, k} - \boldsymbol{u}_k),
\end{equation}
and $\boldsymbol{M}_k \in \mathbb{R}^{n\times n}$ is the diagonal binary matrix defined by Eq.~\eqref{eq: M_k def}, whichselects the measurement positions at each measurement time. The observer trajectory is constrained by the discrete KS dynamics,
\begin{equation}
    \boldsymbol{u}_{\boldsymbol{\theta}, k+1} = \mathsf{f}(\boldsymbol{u}_{\boldsymbol{\theta}, k}).
\end{equation}
We solve this constrained minimization problem by introducing a Lagrangian that enforces the dynamics and deriving the associated adjoint equations. Below, we present the resulting systems for computing the gradient and Hessian of the loss.\par

% Gradient
\subsection{Adjoint system for the gradient}
\label{app: adjoint: gradient}
We denote by $\mathscr{A}$ the discrete initial condition and dynamical constraints,
\begin{equation}
    \mathscr{A} = 
    \underbrace{\boldsymbol{u}_0^{\dagger}
    \left(\boldsymbol{u}_{\boldsymbol{\theta}} - \boldsymbol{\beta}\right)}_\text{initial condition} +
    \underbrace{\sum_{k=0}^{K-1}
    \boldsymbol{u}_{k+1}^\dagger
    \left[\boldsymbol{u}_{\boldsymbol{\theta}, k+1} - \mathsf{f}(\boldsymbol{u}_{\boldsymbol{\theta}, k})\right]}_\text{system dynamics},
\end{equation}
where the adjoint variables $\boldsymbol{u}_k^\dagger$ are Lagrange multipliers, arranged as row vectors, and $\boldsymbol{\beta}$ is the design parameter which determines the initial state, $\boldsymbol{u}_{\boldsymbol{\theta}} = \boldsymbol{\beta}$.
The Lagrangian is built as
\begin{equation}
    \mathscr{L} = \mathscr{J} - \mathscr{A}.
\end{equation}
Substituting $\mathscr{J}$ and $\mathscr{A}$, we get
\begin{equation}
   \mathscr{L} = \sum_{k=0}^{K}
   \mathsf{M}_k\! \left(\boldsymbol{u}_{\boldsymbol{\theta}, k}, \boldsymbol{u}_k\right) - \boldsymbol{u}_0^{\dagger} \left(\boldsymbol{u}_{\boldsymbol{\theta}} - \boldsymbol{\beta}\right) -
   \sum_{k=0}^{K-1} \boldsymbol{u}_{k+1}^{\dagger}\left[\boldsymbol u_{\boldsymbol{\theta}, k+1} - \mathsf{f}(\boldsymbol u_{\boldsymbol{\theta}, k})\right] \!.
\end{equation}
This is rearranged to obtain
\begin{equation}
    \mathscr{L} = \mathsf{M}_K(\boldsymbol u_{\boldsymbol{\theta}, K}, \boldsymbol u_K) -  \boldsymbol{u}_0^{\dagger}(\boldsymbol u_{\boldsymbol{\theta}} - \boldsymbol{\beta}) - \sum_{k=0}^{K-1}
    \left\{
        \boldsymbol{u}_{k+1}^{\dagger}
            \left[\boldsymbol u_{\boldsymbol{\theta}, k+1} - \mathsf{f}(\boldsymbol u_{\boldsymbol{\theta}, k})
        \right] - 
        \mathsf{M}_k(\boldsymbol u_{\boldsymbol{\theta}, k}, \boldsymbol u_k)
    \right\} \!.
\end{equation}
Because the dynamics are enforced during the simulation, the system constraint in $\mathscr{A}$ is always satisfied along a rollout, i.e.,
\begin{equation}
    \mathscr{A} = 0,
\end{equation}
and hence
\begin{equation}
    \mathscr{L} = \mathscr{J}.
\end{equation}
Differentiating $\mathscr{L}$ with respect to $\boldsymbol{\beta}$ gives
\begin{equation}
  \frac{\partial \mathscr{L}}{\partial \boldsymbol{\beta}} = \frac{\partial \mathsf{M}_K}{\partial \boldsymbol u_{\boldsymbol{\theta}, K}}\frac{\partial  \boldsymbol u_{\boldsymbol{\theta}, K}}{\partial  \boldsymbol \beta} + \boldsymbol{u}_0^{\dagger} - \sum_{k=0}^{K-1} \left [\boldsymbol{u}_{k+1}^{\dagger}  \left( \frac{\partial  \boldsymbol u_{\boldsymbol{\theta}, k+1}}{\partial \boldsymbol \beta} - \frac{\partial \mathsf{f}}{\partial \boldsymbol u_{\boldsymbol{\theta}, k}}\frac{\partial  \boldsymbol u_{\boldsymbol{\theta}, k}}{\partial  \boldsymbol\beta}\right) - \frac{\partial \mathsf{M}_k}{\partial \boldsymbol u_{\boldsymbol{\theta}, k}}\frac{\partial  \boldsymbol u_{\boldsymbol{\theta}, k}}{\partial  \boldsymbol \beta} \right] \!.
\end{equation}
Pulling $\boldsymbol{u}_K^{\dagger} (\partial \boldsymbol{u}_{\boldsymbol{\theta}, K}/\partial \boldsymbol{\beta})$ out of the summation, we get
\begin{equation}
    \label{equ: L grad}
    \frac{\partial \mathscr{L}}{\partial \boldsymbol{\beta}} =  \boldsymbol{u}_0^{\dagger} + \left(\frac{\partial \mathsf{M}_K}{\partial \boldsymbol u_{\boldsymbol{\theta}, K}} - \boldsymbol{u}_K^{\dagger}\right)\frac{\partial \boldsymbol u_{\boldsymbol{\theta}, K}}{\partial \boldsymbol \beta} - \sum_{k=0}^{K-1} \left( \boldsymbol {u}_k^{\dagger}  - \boldsymbol{u}_{k+1}^{\dagger}\frac{\partial \mathsf{f}}{\partial \boldsymbol u_{\boldsymbol{\theta}, k} } - \frac{\partial \mathsf{M}_k}{\partial \boldsymbol u_{\boldsymbol{\theta}, k}} \right)\frac{\partial \boldsymbol u_{\boldsymbol{\theta}, k}}{\partial \boldsymbol \beta},
\end{equation}
where
\begin{equation}
    \frac{\partial u_{\boldsymbol{\theta}, k}}{\partial \boldsymbol\beta} = \frac{\partial \boldsymbol u_{\boldsymbol{\theta}, K}}{\partial \boldsymbol u_{\boldsymbol{\theta}, K-1}}\frac{\partial \boldsymbol u_{\boldsymbol{\theta}, K-1}}{\partial \boldsymbol u_{\boldsymbol{\theta}, K-2}} \dots\frac{\partial \boldsymbol u_{\boldsymbol{\theta}, 1}}{\partial \boldsymbol\beta}.
\end{equation}
Unfortunately, direct computation of $\partial \boldsymbol{u}_{\boldsymbol{\theta}, k}/ \partial \boldsymbol{\beta} \in \mathbb{R}^{n\times n}$ is prohibitively expensive. Since $\mathscr{A} = 0$ for all choices of the adjoint variables, we may select a sequence of $\boldsymbol{u}_k^\dagger$ that annihilates the bracketed
coefficients in Eq.~\eqref{equ: L grad}. This yields the discrete adjoint recursion
\begin{subequations}
    \begin{equation}
        \boldsymbol{u}_k^\dagger = 
        \boldsymbol{u}^\dagger_{k+1} \frac{\partial \boldsymbol{u}_{\boldsymbol{\theta},k+1}}{\partial \boldsymbol{u}_{\boldsymbol{\theta}, k}}
        + \frac{\partial \mathsf{M}_k}{\partial \boldsymbol{u}_{\boldsymbol{\theta}, k}},
    \end{equation}
    for $k = K-1, \dots, 0$, with terminal condition
    \begin{equation}
        \boldsymbol{u}_{\boldsymbol{\theta}, K}^\dagger = 
        \frac{\partial \mathsf{M}_k}{\partial \boldsymbol{u}_{\boldsymbol{\theta}, K}}.
    \end{equation}
\end{subequations}
Finally, because $\boldsymbol{u}_{\boldsymbol{\theta},0} = \boldsymbol{\beta}$, we end up with
\begin{equation}
    \label{equ: loss gradient}
    \frac{\partial \mathscr{J}}{\partial \boldsymbol{u}_{\boldsymbol{\theta}}} =   \frac{\partial \mathscr{L}}{\partial \boldsymbol{u}_{\boldsymbol{\theta}}} = \boldsymbol{u}_0^\dagger.
\end{equation}\par

% Hessian
\subsection{Adjoint system for the Hessian}
\label{app: adjoint: Hessian}
Second‐order adjoints are commonly used to compute Hessian--vector products, but the modest dimension of the KS systems considered in this work allows us to form full Hessians. To derive an adjoint system for this purpose, we start by differentiating the transpose of Eq.~\eqref{equ: L grad} with respect to $\boldsymbol{\beta}$:
\begin{equation}
    \begin{split}
     \frac{\partial^2 \mathscr{L}}{\partial \boldsymbol{\beta}^2}
    &=
    \frac{\partial (\boldsymbol{u}_0^{\dagger\,\top})}{\partial \boldsymbol{\beta}}
    +
    \frac{\partial }{\partial \boldsymbol{\beta}}\left(\frac{\partial \boldsymbol{u}_{\boldsymbol{\theta}, K}}{\partial\boldsymbol{\beta}} \right)^\top\left(
        \frac{\partial \mathsf{M}_K}{\partial \boldsymbol{u}_{\boldsymbol{\theta}, K}}
        -
        \boldsymbol{u}_K^{\dagger}
    \right)^\top
    +
    \left(
        \frac{\partial \boldsymbol{u}_{\boldsymbol{\theta}, K}}{\partial \boldsymbol{\beta}}
    \right)^\top
    \left[
        \frac{\partial^2 \mathsf{M}_K}{\partial \boldsymbol{u}_{\boldsymbol{\theta}, K}^2}
        \frac{\partial \boldsymbol{u}_{\boldsymbol{\theta}, K}}{\partial \boldsymbol{\beta}}
        -
        \frac{\partial (\boldsymbol{u}_K^{\dagger\,\top})}{\partial \boldsymbol{\beta}}
    \right] \\
    &\quad
    -\sum_{k=0}^{K-1}
    \frac{\partial }{\partial \boldsymbol{\beta}}\left(\frac{\partial \boldsymbol{u}_{\boldsymbol{\theta}, k}}{\partial\boldsymbol{\beta}} \right)^\top\left(
        \boldsymbol{u}_k^{\dagger}
        -
        \boldsymbol{u}_{k+1}^{\dagger}
        \frac{\partial \mathsf{f}}{\partial \boldsymbol{u}_{\boldsymbol{\theta}, k}}
        -
        \frac{\partial \mathsf{M}_k}{\partial \boldsymbol{u}_{\boldsymbol{\theta}, k}}
    \right)^\top \\
    &\quad
    -\sum_{k=0}^{K-1}
    \left ( \frac{\partial \boldsymbol{u}_{\boldsymbol{\theta}, k}}{\partial \boldsymbol{\beta}}\right )^\top
    \left\{
    \frac{\partial (\boldsymbol{u}_k^{\dagger\,\top})}{\partial \boldsymbol{\beta}} -
    \left [\frac{\partial }{\partial \boldsymbol{u}_{\boldsymbol{\theta},  k}}\left(\frac{\partial \mathsf{f}}{\partial \boldsymbol{u}_{\boldsymbol{\theta}, k}} \right)^\top \right]\frac{\partial \boldsymbol{u}_{\boldsymbol{\theta}, k}}{\partial \boldsymbol{\beta}}\boldsymbol{u}^{\dagger\,\top}_{k+1}
     - 
    \left(\frac{\partial \mathsf{f}}{\partial \boldsymbol{u}_{\boldsymbol{\theta}, k}}\right)^\top \frac{\partial (\boldsymbol{u}^{\dagger\,\top}_{k+1})}{\partial \boldsymbol{\beta}} -
    \frac{\partial ^2\mathsf{M}_k}{\partial \boldsymbol{u}_{\boldsymbol{\theta}, k}^2}\frac{\partial \boldsymbol{u}_{\boldsymbol{\theta}, k}}{\partial \boldsymbol{\beta}}
    \right\} \!.
    \end{split}
\end{equation}
After substituting $\boldsymbol{u}_{\boldsymbol{\theta}} = \boldsymbol{\beta}$ and $\boldsymbol{u}_{\boldsymbol{\theta}, k+1} = \mathsf{f}(\boldsymbol{u}_{\boldsymbol{\theta}, k})$, we get our second-order adjoint system:
\begin{subequations}
    \begin{equation}
        \frac{\partial (\boldsymbol{u}_k^{\dagger\,\top})}{\partial \boldsymbol{u}_{\boldsymbol{\theta}}}
        = \left [\frac{\partial }{\partial \boldsymbol{u}_{\boldsymbol{\theta},  k}}\left(\frac{\partial \boldsymbol{u}_{\boldsymbol{\theta}, k+1}}{\partial \boldsymbol{u}_{\boldsymbol{\theta}, k}} \right)^\top \right]\frac{\partial \boldsymbol{u}_{\boldsymbol{\theta}, k}}{\partial \boldsymbol{u}_{\boldsymbol{\theta}}}\boldsymbol{u}^{\dagger\,\top}_{k+1}
        + \left(\frac{\partial \boldsymbol{u}_{\boldsymbol{\theta},k+1}}{\partial \boldsymbol{u}_{\boldsymbol{\theta}, k}}\right)^\top
        \frac{\partial (\boldsymbol{u}_{k+1}^{\dagger\,\top})}{\partial \boldsymbol{u}_{\boldsymbol{\theta}}}
        + \frac{\partial^2 \mathsf{M}_k}{\partial \boldsymbol{u}_{\boldsymbol{\theta}, k}^2}
        \frac{\partial \boldsymbol{u}_{\boldsymbol{\theta}, k}}{\partial \boldsymbol{u}_{\boldsymbol{\theta}}},
    \end{equation}
    for $k = K-1, \dots, 0$, with terminal condition
    \begin{equation}
        \frac{\partial (\boldsymbol{u}_{\boldsymbol{\theta}, K}^{\dagger\,\top})}{\partial \boldsymbol{u}_{\boldsymbol{\theta}}} =
        \frac{\partial^2 \mathsf{M}_k}{\partial \boldsymbol{u}_{\boldsymbol{\theta}, K}^2}
        \frac{\partial \boldsymbol{u}_{\boldsymbol{\theta}, K}}{\partial \boldsymbol{u}_{\boldsymbol{\theta}}}.
    \end{equation}
\end{subequations}
This adjoint system provides the Hessian via
\begin{equation}
    \label{equ: adjoint Hessian}
    \frac{\partial^2 \mathscr{J}}{\partial \boldsymbol{u}_{\boldsymbol{\theta}}} =         \frac{\partial^2 \mathscr{L}}{\partial \boldsymbol{u}_{\boldsymbol{\theta}}} = \frac{\partial (\boldsymbol{u}_0^{\dagger\,\top})}{\partial\boldsymbol{u}_{\boldsymbol{\theta}}}.
\end{equation}\par

%%% References %%%


\begin{thebibliography}{10}

\bibitem{Smith2014}
N.~Smith and W.~D. Arnett, \enquote{{Preparing for an explosion: hydrodynamic instabilities and turbulence in presupernovae},} Astrophys. J. \textbf{785}, 82 (2014).

\bibitem{Marati2004}
N.~Marati, C.~M. Casciola, and R.~Piva, \enquote{{Energy cascade and spatial fluxes in wall turbulence},} J. Fluid Mech. \textbf{521}, 191--215 (2004).

\bibitem{Bolotnov2010}
I.~A. Bolotnov, R.~T. Lahey~Jr, D.~A. Drew, K.~E. Jansen, and A.~A. Oberai, \enquote{{Spectral analysis of turbulence based on the DNS of a channel flow},} Comput. Fluids \textbf{39}, 640--655 (2010).

\bibitem{Robinson1991}
S.~K. Robinson, \enquote{{Coherent motions in the turbulent boundary layer},} Annu. Rev. Fluid Mech. \textbf{23}, 601--639 (1991).

\bibitem{Smits2011}
A.~J. Smits, B.~J. McKeon, and I.~Marusic, \enquote{{High--Reynolds number wall turbulence},} Annu. Rev. Fluid Mech. \textbf{43}, 353--375 (2011).

\bibitem{Holmes2012}
P.~Holmes, J.~L. Lumley, G.~Berkooz, and C.~W. Rowley, \emph{{One-dimensional “turbulence”}} (Cambridge University Press, 2012), p. 214–235, Cambridge Monographs on Mechanics.

\bibitem{Linot2023}
A.~J. Linot, K.~Zeng, and M.~D. Graham, \enquote{{Turbulence control in plane Couette flow using low-dimensional neural ODE-based models and deep reinforcement learning},} Int. J. Heat Fluid Flow \textbf{101}, 109139 (2023).

\bibitem{Duraisamy2019}
K.~Duraisamy, G.~Iaccarino, and H.~Xiao, \enquote{{Turbulence modeling in the age of data},} Annu. Rev. Fluid Mech. \textbf{51}, 357--377 (2019).

\bibitem{Argyropoulos2015}
C.~D. Argyropoulos and N.~Markatos, \enquote{{Recent advances on the numerical modelling of turbulent flows},} Appl. Math. Modell. \textbf{39}, 693--732 (2015).

\bibitem{Slotnick2014}
J.~P. Slotnick, A.~Khodadoust, J.~Alonso, D.~Darmofal, W.~Gropp, E.~Lurie, and D.~J. Mavriplis, \enquote{{CFD Vision 2030 Study: A Path to Revolutionary Computational Aerosciences},} Tech. Rep. NF1676L-18332, National Aeronautics and Space Administration (2014).

\bibitem{Cary2021}
A.~W. Cary, J.~Chawner, E.~P. Duque, W.~Gropp, W.~L. Kleb, R.~M. Kolonay, E.~Nielsen, and B.~Smith, \enquote{{Cfd vision 2030 road map: Progress and perspectives},} in \enquote{AIAA aviation 2021 forum,}  (2021), p. 2726.

\bibitem{Gronskis2013}
A.~Gronskis, D.~Heitz, and E.~M{\'e}min, \enquote{{Inflow and initial conditions for direct numerical simulation based on adjoint data assimilation},} J. Comput. Phys. \textbf{242}, 480--497 (2013).

\bibitem{Asch2016}
M.~Asch, M.~Bocquet, and M.~Nodet, \emph{{Data Assimilation: Methods, Algorithms, and Applications}} (SIAM, 2016).

\bibitem{Hayase2015}
T.~Hayase, \enquote{{Numerical simulation of real-world flows},} Fluid Dyn. Res. \textbf{47}, 051201 (2015).

\bibitem{Zaki2025a}
T.~A. Zaki and M.~Wang, \enquote{{Data assimilation and flow estimation},} in \enquote{Data Driven Analysis and Modeling of Turbulent Flows,}  (Elsevier, 2025), pp. 129--181.

\bibitem{Zaki2025b}
T.~A. Zaki, \enquote{{Turbulence from an observer perspective},} Annu. Rev. Fluid Mech. \textbf{57} (2025).

\bibitem{He2025}
C.~He, S.~Li, and Y.~Liu, \enquote{{Data assimilation: new impetus in experimental fluid dynamics},} Exp. Fluids \textbf{66}, 1--24 (2025).

\bibitem{Welch1995}
G.~Welch and G.~Bishop, \enquote{{An introduction to the Kalman filter},} Tech. Rep. TR 95-041, University of North Carolina (1995).

\bibitem{Evensen1994}
G.~Evensen, \enquote{{Sequential data assimilation with a nonlinear quasi-geostrophic model using Monte Carlo methods to forecast error statistics},} J. Geophys. Res.: Oceans \textbf{99}, 10143--10162 (1994).

\bibitem{Clark2020}
P.~Clark Di~Leoni, A.~Mazzino, and L.~Biferale, \enquote{{Synchronization to big data: Nudging the Navier-Stokes equations for data assimilation of turbulent flows},} Phys. Rev. X \textbf{10}, 011023 (2020).

\bibitem{Vela2021}
A.~Vela-Mart{\'\i}n, \enquote{{The synchronisation of intense vorticity in isotropic turbulence},} J. Fluid Mech. \textbf{913}, R8 (2021).

\bibitem{Lalescu2013}
C.~C. Lalescu, C.~Meneveau, and G.~L. Eyink, \enquote{{Synchronization of chaos in fully developed turbulence},} Phys. Rev. Lett. \textbf{110}, 084102 (2013).

\bibitem{Raissi2019}
M.~Raissi, P.~Perdikaris, and G.~E. Karniadakis, \enquote{{Physics-informed neural networks: A deep learning framework for solving forward and inverse problems involving nonlinear partial differential equations},} J. Comput. Phys. \textbf{378}, 686--707 (2019).

\bibitem{Gesemann2016}
S.~Gesemann, F.~Huhn, D.~Schanz, and A.~Schr{\"o}der, \enquote{{From noisy particle tracks to velocity, acceleration and pressure fields using B-splines and penalties},} in \enquote{18th international symposium on applications of laser and imaging techniques to fluid mechanics, Lisbon, Portugal,} , vol.~4 (2016), vol.~4.

\bibitem{Casa2013}
L.~Casa and P.~Krueger, \enquote{{Radial basis function interpolation of unstructured, three-dimensional, volumetric particle tracking velocimetry data},} Meas. Sci. Technol. \textbf{24}, 065304 (2013).

\bibitem{Wang2015}
H.~Wang, Q.~Gao, L.~Feng, R.~Wei, and J.~Wang, \enquote{{Proper orthogonal decomposition based outlier correction for PIV data},} Exp. Fluids \textbf{56}, 43 (2015).

\bibitem{Le1986}
F.-X. Le~Dimet and O.~Talagrand, \enquote{{Variational algorithms for analysis and assimilation of meteorological observations: theoretical aspects},} Tellus A: Dyn. Meteorol. Oceanogr. \textbf{38}, 97--110 (1986).

\bibitem{Wang2019}
M.~Wang, Q.~Wang, and T.~A. Zaki, \enquote{{Discrete adjoint of fractional-step incompressible Navier-Stokes solver in curvilinear coordinates and application to data assimilation},} J. Comput. Phys. \textbf{396}, 427--450 (2019).

\bibitem{Liu2008}
C.~Liu, Q.~Xiao, and B.~Wang, \enquote{An ensemble-based four-dimensional variational data assimilation scheme. part i: Technical formulation and preliminary test,} Mon. Weather Rev. \textbf{136}, 3363--3373 (2008).

\bibitem{Mons2016}
V.~Mons, J.-C. Chassaing, T.~Gomez, and P.~Sagaut, \enquote{{Reconstruction of unsteady viscous flows using data assimilation schemes},} J. Comput. Phys. \textbf{316}, 255--280 (2016).

\bibitem{Zelik2014}
S.~Zelik, \enquote{{Inertial manifolds and finite-dimensional reduction for dissipative PDEs},} Proc. R. Soc. A \textbf{144}, 1245–1327 (2014).

\bibitem{Kugiumtzis1996}
D.~Kugiumtzis, \enquote{{State space reconstruction parameters in the analysis of chaotic time series—the role of the time window length},} Physica D \textbf{95}, 13--28 (1996).

\bibitem{Deyle2011}
E.~R. Deyle and G.~Sugihara, \enquote{{Generalized theorems for nonlinear state space reconstruction},} Plos one \textbf{6}, e18295 (2011).

\bibitem{Casdagli1991}
M.~Casdagli, S.~Eubank, J.~D. Farmer, and J.~Gibson, \enquote{{State space reconstruction in the presence of noise},} Physica D \textbf{51}, 52--98 (1991).

\bibitem{Wittenberg1998}
R.~W. Wittenberg, \enquote{{Local Dynamics and Spatiotemporal Chaos. The Kuramoto--Sivashinsky Equation: A Case Study},} Ph.D. thesis, Princeton University (1998).

\bibitem{Ding2016}
X.~Ding, H.~Chat{\'e}, P.~Cvitanovi{\'c}, E.~Siminos, and K.~Takeuchi, \enquote{{Estimating the dimension of an inertial manifold from unstable periodic orbits},} Phys. Rev. Lett. \textbf{117}, 024101 (2016).

\bibitem{Linot2020}
A.~J. Linot and M.~D. Graham, \enquote{{Deep learning to discover and predict dynamics on an inertial manifold},} Phys. Rev. E \textbf{101}, 062209 (2020).

\bibitem{Linot2022}
A.~J. Linot and M.~D. Graham, \enquote{{Data-driven reduced-order modeling of spatiotemporal chaos with neural ordinary differential equations},} Chaos \textbf{32} (2022).

\bibitem{Zeng2024}
K.~Zeng, C.~E.~P. De~Jesus, A.~J. Fox, and M.~D. Graham, \enquote{{Autoencoders for discovering manifold dimension and coordinates in data from complex dynamical systems},} Mach. Learn.: Sci. Technol. \textbf{5}, 025053 (2024).

\bibitem{Yang2009}
H.-l. Yang, K.~A. Takeuchi, F.~Ginelli, H.~Chat{\'e}, and G.~Radons, \enquote{{Hyperbolicity and the effective dimension of spatially extended dissipative systems},} Phys. Rev. Lett. \textbf{102}, 074102 (2009).

\bibitem{Takeuchi2011}
K.~A. Takeuchi, H.-l. Yang, F.~Ginelli, G.~Radons, and H.~Chat{\'e}, \enquote{{Hyperbolic decoupling of tangent space and effective dimension of dissipative systems},} Phys. Rev. E \textbf{84}, 046214 (2011).

\bibitem{Foias1986}
C.~Foias, B.~Nicolaenko, G.~R. Sell, and R.~Temam, \enquote{{Inertial manifolds for the Kuramoto--Sivashinsky equation and an estimate of their lowest dimension},} J. Math. Pures Appl. \textbf{67}, 197--226 (1988).

\bibitem{Jones1996}
D.~A. Jones and E.~S. Titi, \enquote{{$C^1$ Approximations of Inertial Manifolds for Dissipative Nonlinear Equations},} Journal of Differential Equations \textbf{127}, 54--86 (1996).

\bibitem{Cox2002}
S.~M. Cox and P.~C. Matthews, \enquote{{Exponential time differencing for stiff systems},} J. Comput. Phys. \textbf{176}, 430--455 (2002).

\bibitem{Cvitanovic2010}
P.~Cvitanovi{\'c}, R.~L. Davidchack, and E.~Siminos, \enquote{{On the state space geometry of the Kuramoto--Sivashinsky flow in a periodic domain},} SIAM J. Appl. Dyn. Syst. \textbf{9}, 1--33 (2010).

\bibitem{Edson2019}
R.~A. Edson, J.~E. Bunder, T.~W. Mattner, and A.~J. Roberts, \enquote{{Lyapunov exponents of the Kuramoto--Sivashinsky PDE},} ANZIAM J. \textbf{61}, 270--285 (2019).

\bibitem{Benettin1980}
G.~Benettin, L.~Galgani, A.~Giorgilli, and J.-M. Strelcyn, \enquote{{Lyapunov characteristic exponents for smooth dynamical systems and for Hamiltonian systems; a method for computing all of them. Part 1: Theory},} Meccanica \textbf{15}, 9--20 (1980).

\bibitem{Protas2004}
B.~Protas, T.~R. Bewley, and G.~Hagen, \enquote{{A computational framework for the regularization of adjoint analysis in multiscale PDE systems},} J. Comput. Phys. \textbf{195}, 49--89 (2004).

\bibitem{Jardak2010}
M.~Jardak, I.~M. Navon, and M.~Zupanski, \enquote{{Comparison of sequential data assimilation methods for the Kuramoto--Sivashinsky equation},} Int. J. Numer. Methods Fluids \textbf{62}, 374--402 (2010).

\bibitem{Chandramoorthy2019}
N.~Chandramoorthy, P.~Fernandez, C.~Talnikar, and Q.~Wang, \enquote{{Feasibility analysis of ensemble sensitivity computation in turbulent flows},} AIAAJ \textbf{57}, 4514--4526 (2019).

\bibitem{Zaki2021}
T.~A. Zaki and M.~Wang, \enquote{{From limited observations to the state of turbulence: Fundamental difficulties of flow reconstruction},} Phys. Rev. Fluids \textbf{6}, 100501 (2021).

\bibitem{Li2020}
Y.~Li, J.~Zhang, G.~Dong, and N.~S. Abdullah, \enquote{{Small-scale reconstruction in three-dimensional Kolmogorov flows using four-dimensional variational data assimilation},} J. Fluid Mech. \textbf{885}, A9 (2020).

\bibitem{Chandramouli2020}
P.~Chandramouli, E.~M{\'e}min, and D.~Heitz, \enquote{{4D large scale variational data assimilation of a turbulent flow with a dynamics error model},} J. Comput. Phys. \textbf{412}, 109446 (2020).

\bibitem{Wang2025}
M.~Wang and T.~A. Zaki, \enquote{{Variational data assimilation in wall turbulence: from outer observations to wall stress and pressure},} J. Fluid Mech. \textbf{1008}, A26 (2025).

\bibitem{Dauphin2014}
Y.~N. Dauphin, R.~Pascanu, C.~Gulcehre, K.~Cho, S.~Ganguli, and Y.~Bengio, \enquote{{Identifying and attacking the saddle point problem in high-dimensional non-convex optimization},} Adv. Neural Inf. Process. Syst. \textbf{27} (2014).

\bibitem{Paternain2019}
S.~Paternain, A.~Mokhtari, and A.~Ribeiro, \enquote{{A Newton-based method for nonconvex optimization with fast evasion of saddle points},} SIAM J. Optim. \textbf{29}, 343--368 (2019).

\bibitem{Sun2016}
J.~Sun, Q.~Qu, and J.~Wright, \enquote{{Complete dictionary recovery over the sphere I: Overview and the geometric picture},} IEEE Trans. Inf. Theory \textbf{63}, 853--884 (2016).

\bibitem{Ge2015}
R.~Ge, F.~Huang, C.~Jin, and Y.~Yuan, \enquote{{Escaping from saddle points—online stochastic gradient for tensor decomposition},} in \enquote{Conference on learning theory,}  (PMLR, 2015), pp. 797--842.

\bibitem{Ge2016}
R.~Ge, J.~D. Lee, and T.~Ma, \enquote{{Matrix completion has no spurious local minimum},} Adv. Neural Inf. Process. Syst. \textbf{29} (2016).

\bibitem{Kawaguchi2016}
K.~Kawaguchi, \enquote{{Deep learning without poor local minima},} Adv. Neural Inf. Process. Syst. \textbf{29} (2016).

\bibitem{Bray2007}
A.~J. Bray and D.~S. Dean, \enquote{{Statistics of critical points of Gaussian fields on large-dimensional spaces},} Phys. Rev. Lett. \textbf{98}, 150201 (2007).

\bibitem{Fyodorov2007}
Y.~V. Fyodorov and I.~Williams, \enquote{{Replica symmetry breaking condition exposed by random matrix calculation of landscape complexity},} J. Stat. Phys. \textbf{129}, 1081--1116 (2007).

\bibitem{Baldi1989}
P.~Baldi and K.~Hornik, \enquote{{Neural networks and principal component analysis: Learning from examples without local minima},} Neural Networks \textbf{2}, 53--58 (1989).

\bibitem{Saxe2013}
A.~M. Saxe, J.~L. McClellans, and S.~Ganguli, \enquote{{Learning hierarchical categories in deep neural networks},} in \enquote{Proceedings of the Annual Meeting of the Cognitive Science Society,} , vol.~35 (2013), vol.~35.

\bibitem{Saarinen1993}
S.~Saarinen, R.~Bramley, and G.~Cybenko, \enquote{{Ill-conditioning in neural network training problems},} SIAM J. Sci. Comput. \textbf{14}, 693--714 (1993).

\bibitem{Greenstadt1967}
J.~Greenstadt, \enquote{{On the relative efficiencies of gradient methods},} Math. Comput. \textbf{21}, 360--367 (1967).

\bibitem{Gould1998}
N.~I. Gould and J.~Nocedal, \enquote{{The modified absolute-value factorization norm for trust-region minimization},} in \enquote{High Performance Algorithms and Software in Nonlinear Optimization,}  (Springer, 1998), pp. 225--241.

\bibitem{Gejadze2023}
I.~Gejadze, V.~Shutyaev, H.~Oubanas, and P.-O. Malaterre, \enquote{{A Bayesian-variational cyclic method for solving estimation problems characterized by non-uniqueness (equifinality)},} J. Comput. Phys. \textbf{488}, 112239 (2023).

\bibitem{Haben2011}
S.~A. Haben, A.~S. Lawless, and N.~K. Nichols, \enquote{{Conditioning and preconditioning of the variational data assimilation problem},} Comput. Fluids \textbf{46}, 252--256 (2011).

\bibitem{Ke2026}
H.~Ke, Z.~You, and Q.~Wang, \enquote{{Preconditioned adjoint data assimilation for two-dimensional decaying isotropic turbulence},} arXiv preprint arXiv:2602.14016  (2026).

\bibitem{Dembo1983}
R.~S. Dembo and T.~Steihaug, \enquote{{Truncated-Newton algorithms for large-scale unconstrained optimization},} Math. Program. \textbf{26}, 190--212 (1983).

\bibitem{Hopf1948}
E.~Hopf, \enquote{{A mathematical example displaying features of turbulence},} Commun. Pure Appl. Math. \textbf{1}, 303--322 (1948).

\bibitem{Takens1981}
F.~Takens, \enquote{{Detecting Strange Attractors in Turbulence},} in \enquote{Dynamical Systems and Turbulence, Warwick 1980,} , vol. 898 of \emph{Lecture Notes in Mathematics}, D.~A. Rand and L.-S. Young, eds. (Springer, Berlin, Heidelberg, 1981), vol. 898 of \emph{Lecture Notes in Mathematics}, pp. 366--381.

\bibitem{Noakes1991}
L.~Noakes, \enquote{{The Takens embedding theorem},} Int. J. Bifurcation Chaos \textbf{1}, 867--872 (1991).

\bibitem{Sauer1991}
T.~Sauer, J.~A. Yorke, and M.~Casdagli, \enquote{{Embedology},} J. Stat. Phys. \textbf{65}, 579--616 (1991).

\bibitem{Kugiumtzis1994}
D.~Kugiumtzis, B.~Lillekjendlie, and N.~Christophersen, \enquote{{Chaotic time series. Part I. Estimation of some invariant properties in state-space},} Int. J. Modell. Identif. Control \textbf{15} (1994).

\bibitem{Fraser1986}
A.~M. Fraser and H.~L. Swinney, \enquote{{Independent coordinates for strange attractors from mutual information},} Phys. Rev. A \textbf{33}, 1134 (1986).

\bibitem{Abarbanel1993}
H.~D. Abarbanel, R.~Brown, J.~J. Sidorowich, and L.~S. Tsimring, \enquote{{The analysis of observed chaotic data in physical systems},} Rev. Mod. Phys. \textbf{65}, 1331 (1993).

\bibitem{Rosenstein1994}
M.~T. Rosenstein, J.~J. Collins, and C.~J. De~Luca, \enquote{{Reconstruction expansion as a geometry-based framework for choosing proper delay times},} Physica D \textbf{73}, 82--98 (1994).

\bibitem{Caputo1986}
J.~Caputo, B.~Malraison, and P.~Atten, \enquote{{Determination of attractor dimension and entropy for various flows: An experimentalist’s viewpoint},} in \enquote{Dimensions and Entropies in Chaotic Systems: Quantification of Complex Behavior,}  (Springer, 1986), pp. 180--190.

\bibitem{Gibson1992}
J.~F. Gibson, J.~D. Farmer, M.~Casdagli, and S.~Eubank, \enquote{{An analytic approach to practical state space reconstruction},} Physica D \textbf{57}, 1--30 (1992).

\bibitem{bishop2012}
R.~L. Bishop and S.~I. Goldberg, \emph{{Tensor analysis on manifolds}} (Courier Corporation, 2012).

\bibitem{Floryan2022}
D.~Floryan and M.~D. Graham, \enquote{{Data-driven discovery of intrinsic dynamics},} Nat. Mach. Intell. \textbf{4}, 1113--1120 (2022).

\bibitem{Constantin1985}
P.~Constantin, C.~Foias, O.~P. Manley, and R.~Temam, \enquote{{Determining modes and fractal dimension of turbulent flows},} J. Fluid Mech. \textbf{150}, 427--440 (1985).

\bibitem{Temam2012}
R.~Temam, \emph{{Infinite-dimensional dynamical systems in mechanics and physics}}, vol.~68 (Springer Science \& Business Media, 2012).

\bibitem{Cleary2025}
A.~Cleary and J.~Page, \enquote{{Characterizing the Reynolds number dependence of the chaotic attractor in two-dimensional turbulence with dimension-minimizing autoencoders},} Phys. Rev. E \textbf{112}, 055105 (2025).

\bibitem{Oseledec1968}
V.~I. Oseledec, \enquote{{A multiplicative ergodic theorem, Lyapunov characteristic numbers for dynamical systems},} Trans. Mosc. Math. Soc. \textbf{19}, 197--231 (1968).

\bibitem{Kendall2018}
A.~Kendall, Y.~Gal, and R.~Cipolla, \enquote{Multi-task learning using uncertainty to weigh losses for scene geometry and semantics,} in \enquote{Proceedings of the IEEE conference on computer vision and pattern recognition,}  (2018), pp. 7482--7491.

\bibitem{Yu2020}
T.~Yu, S.~Kumar, A.~Gupta, S.~Levine, K.~Hausman, and C.~Finn, \enquote{Gradient surgery for multi-task learning,} NeurIPS \textbf{33}, 5824--5836 (2020).

\bibitem{Holland2002}
R.~Holland, \enquote{{Finite-difference time-domain (FDTD) analysis of magnetic diffusion},} IEEE Trans. Electromagn. Compat. \textbf{36}, 32--39 (2002).

\bibitem{Petropoulos2002}
P.~G. Petropoulos, \enquote{{Analysis of exponential time-differencing for FDTD in lossy dielectrics},} IEEE Trans. Antennas Propag. \textbf{45}, 1054--1057 (2002).

\bibitem{Schuster2000}
C.~Schuster, A.~Christ, and W.~Fichtner, \enquote{{Review of FDTD time-stepping schemes for efficient simulation of electric conductive media},} Microwave Opt. Technol. Lett. pp. 16--21 (2000).

\bibitem{Sandri1996}
M.~Sandri, \enquote{{Numerical calculation of Lyapunov exponents},} Mathematica Journal \textbf{6}, 78--84 (1996).

\bibitem{Jing2020}
L.~Jing, J.~Zbontar, and Y.~LeCun, \enquote{{Implicit rank-minimizing autoencoder},} Adv. Neural Inf. Process. Syst. \textbf{33}, 14736--14746 (2020).

\end{thebibliography}
\end{document}